\pdfoutput=1
\documentclass[final,3p]{elsarticle}
\usepackage{graphicx}
\usepackage{subfigure}
\usepackage{subfigmat}
\usepackage{amsmath,amsfonts,stmaryrd,amssymb,theorem}
\usepackage[colorlinks]{hyperref}
\usepackage{rotating}
\usepackage{subeqnarray}
\usepackage{multirow}
\usepackage[ruled]{algorithm2e}
\usepackage{algorithmic}
\usepackage{xcolor,colortbl}
\usepackage{lineno}

\newtheorem{lemma}{Lemma}

\theoremstyle{plain}
\newtheorem{remark}{Remark}
\graphicspath{{./}}

\newcommand{\be}{\begin{eqnarray}}
\newcommand{\ee}{\end{eqnarray}}
\newcommand{\ben}{\begin{eqnarray*}}
\newcommand{\een}{\end{eqnarray*}}
\newcommand{\bi}{\begin{itemize}}
\newcommand{\ei}{\end{itemize}}

\newcommand{\eph}{_{e+1/2}}
\newcommand{\enh}{_{e-1/2}}
\DeclareMathAlphabet{\mathpzc}{OT1}{pzc}{m}{it}

\journal{J. Comp. Phys.}

\begin{document}
%
%
\makeatletter
\newcommand\rlarrows{\mathop{\operator@font \rightleftarrows}\nolimits}
\makeatother
%
\newenvironment{itemizePacked}{
\begin{itemize}
  \setlength{\itemsep}{1pt}
  \setlength{\parskip}{0pt}
  \setlength{\parsep}{0pt}
}{\end{itemize}}
\newenvironment{enumeratePacked}{
\begin{enumerate}
  \setlength{\itemsep}{5pt}
  \setlength{\parskip}{0pt}
  \setlength{\parsep}{0pt}
}{\end{enumerate}}
%
\newcommand{\f}[2]{{\frac{#1}{#2}}}
\newcommand{\wt}[1]{{\widetilde{#1}}}
\newcommand{\wh}[1]{{\widehat{#1}}}
\newcommand{\wc}[1]{{\widecheck{#1}}}
\newcommand{\chem}[1]{\ensuremath{\mathrm{#1}}}
\newcommand{\lrangle}[1]{{\langle{#1}\rangle}}
\newcommand{\lrcurl}[1]{{\{{#1}\}}}
\newcommand{\ol}[1]{{\overline{#1}}}
\newcommand{\ul}[1]{{\underline{#1}}}
\newcommand{\tr}{{\scriptscriptstyle\mathsf T}}
\newcommand{\dd}{{\scriptscriptstyle\Delta}}
\newcommand{\suble}{{\scriptscriptstyle <}}
\newcommand{\eps}{{\varepsilon}}
\def\REFUP{\rm{ref}}
\def\REF{_{\REFUP}}
\def\RPP{reaction progress parameter}

\newcommand{\bfit}[1]{\textbf{\textit{#1}}}

%
%
\newcommand{\colred}[1]{{\color{red} #1}}
\newcommand{\colblue}[1]{{\color{blue} #1}}
\newcommand{\colwhite}[1]{{\color{white} #1}}
\newcommand{\colgreen}[1]{{\color{green} #1}}
\newcommand{\colbrown}[1]{{\color{Brown} #1}}
\newcommand{\colfucsia}[1]{{\color{Fuchsia} #1}}
\newcommand{\colBlue}[1]{{\color{Blue} #1}}
\newcommand{\corrections}[1]{{\color{blue}\it#1}}
\newcommand{\comment}[1]{{\color{blue}\bf#1}}
%
%
\def\chist{\chi_{Z,\up{st}}}
\def\chiq{\chi_{Z,\up{q}}}
\def\chii{\chi_{Z,\up{i}}}
%
%
\def\cldbpage{\clearpage{\pagestyle{empty}\cleardoublepage}}
%
%
\def\cA{{\cal{A}}}
\def\cB{{\cal{B}}}
\def\cC{{\cal{C}}}
\def\cO{{\cal{O}}}
\def\cD{{\cal{D}}}
\def\cE{{\cal{E}}}
\def\cF{{\cal{F}}}
\def\cH{{\cal{H}}}     
\def\cG{{\cal{G}}}     
\def\cN{{\cal{N}}}     
\def\cL{{\cal{L}}}     
\def\cS{{\cal{S}}}     
\def\cT{{\cal{T}}}
\def\cU{{\cal{U}}}
\def\cC{{\cal{C}}}     %
\def\cM{{\cal{M}}}     %
\def\cP{{\cal{P}}}     %
\def\cR{{\cal{R}}}     %
\def\cV{{\cal{V}}}     %
\def\cQ{{\cal{Q}}}     %
%
%
\def\erf{{\rm{erf}}}
%
%
\def\up{\textup}
\def\p{\partial}
\def\d{\textup d}
\def\D{\displaystyle}
\def\OmFint{\iiint\limits_{\OmF}}
\def\OmF{{\Omega_{\cal F}}}
\def\OmA{{\Omega_{\cal A}}}
\def\e{\textup{e}}
\def\i{\textup{i}}
%
%
\def\Fr{{\rm{Fr}}}
\def\Ma{{\rm{Ma}}}
\def\Re{{\rm{Re}}}
\def\Kn{{\rm{Kn}}}
\def\Rd{{\rm{Rd}}}
\def\Le{{\rm{Le}}}
\def\Da{{\rm{Da}}}
\def\Ka{{\rm{Ka}}}
\def\Nu{{\rm{Nu}}}
\def\Sc{{\rm{Sc}}}
\def\Ri{{\rm{Ri}}}
\def\Ec{{\rm{Ec}}}
\def\Tu{{\rm{Tu}}}
\def\St{{\rm{St}}}
\def\Mi{{\rm{Mi}}}
\def\Ra{{\rm{Ra}}}
%
%
 \def\avec{{\mbox{\boldmath$a$}}}
 \def\bvec{{\mbox{\boldmath$b$}}}
 \def\Bvec{{\mbox{\boldmath$B$}}}
 \def\cvec{{\mbox{\boldmath$c$}}}
 \def\dvec{{\mbox{\boldmath$d$}}}
 \def\evec{{\mbox{\boldmath$e$}}}
 \def\Fvec{{\mbox{\boldmath$F$}}}
 \def\Nvec{{\mbox{\boldmath$N$}}}
 \def\fvec{{\mbox{\boldmath$f$}}}
 \def\gvec{{\mbox{\boldmath$g$}}}
 \def\hvec{{\mbox{\boldmath$h$}}}
 \def\ivec{{\mbox{\boldmath$i$}}}
 \def\jvec{{\mbox{\boldmath$j$}}}
 \def\kvec{{\mbox{\boldmath$k$}}}
 \def\pvec{{\mbox{\boldmath$p$}}}
 \def\Pvec{{\mbox{\boldmath$P$}}}
 \def\uvec{{\mbox{\boldmath$u$}}}
 \def\Uvec{{\mbox{\boldmath$U$}}}
 \def\nvec{{\mbox{\boldmath$n$}}}
 \def\tvec{{\mbox{\boldmath$t$}}}
 \def\Rvec{{\mbox{\boldmath$R$}}}
 \def\rvec{{\mbox{\boldmath$r$}}}
 \def\svec{{\mbox{\boldmath$s$}}}
 \def\Svec{{\mbox{\boldmath$S$}}}
 \def\xvec{{\mbox{\boldmath$x$}}}
 \def\vvec{{\mbox{\boldmath$v$}}}
 \def\wvec{{\mbox{\boldmath$w$}}}
 \def\yvec{{\mbox{\boldmath$y$}}}
 \def\mvec{{\mbox{\boldmath$m$}}}
 \def\zvec{{\mbox{\boldmath$z$}}}
 \def\Xvec{{\mbox{\boldmath$X$}}}
 \def\qvec{{\mbox{\boldmath$q$}}}
 \def\0vec{{\mbox{\boldmath$0$}}}
 \def\xivec{{\mbox{\boldmath$\xi$}}}
 \def\wpvec{{\boldsymbol{\wp}}}
 \def\psivec{{\mbox{\boldmath$\psi$}}}
 \def\epsvec{{\mbox{\boldmath$\epsilon$}}}
 \def\phivec{{\mbox{\boldmath$\phi$}}}
 \def\varphivec{{\mbox{\boldmath$\varphi$}}}
 \def\zetavec{{\mbox{\boldmath$\zeta$}}}
 \def\kappavec{{\mbox{\boldmath$\kappa$}}}
 \def\varkappavec{{\pmb{\varkappa}}}
 \def\etavec{{\mbox{\boldmath$\eta$}}}
 \def\Psivec{{\boldsymbol{\Psi}}}
 \def\Wvec{{\boldsymbol{W}}}
 \def\Yvec{{\mbox{\boldmath$Y$}}}
 \def\Vvec{{\mbox{\boldmath$V$}}}
 \def\cLvec{{\boldsymbol{\cal{L}}}}
 \def\cMvec{{\boldsymbol{\cal{M}}}}
 \def\cHvec{{\boldsymbol{\cal{H}}}}
 \def\cVvec{{\boldsymbol{\cal{V}}}}
 \def\omegavec{{\mbox{\boldmath$\omega$}}}
 \def\Omegavec{{\mbox{\boldmath$\Omega$}}}
 \def\sigmavec{{\boldsymbol{\sigma}}}
 \def\Amat{{\underline{\underline{{A}}}}}
 \def\Bmat{{\underline{\underline{{B}}}}}
 \def\taumat{{\underline{\underline{{\tau}}}}}
 \def\sigmamat{{\underline{\underline{{\sigma}}}}}
 \def\Cmat{{\underline{\underline{{C}}}}}
 \def\Imat{{\underline{\underline{{I}}}}}
 \def\Smat{{\underline{\underline{{S}}}}}
 \def\Rmat{{\underline{\underline{{R}}}}}
 \def\Tmat{{\underline{\underline{{T}}}}}
 \def\Emat{{\underline{\underline{{E}}}}}
 \def\tmat{{\underline{\underline{{t}}}}}

 \def\alphavec{{\boldsymbol{\alpha}}}
 \def\betavec{{\boldsymbol{\beta}}}
 \def\tauvec{{\boldsymbol{\tau}}}
 \def\thetavec{{\boldsymbol{\theta}}}
 \def\lambdavec{{\mbox{\boldmath$\lambda$}}}
 \def\cTmat{{\underline{\underline{{\cal{T}}}}}}
 \def\cLmat{{\underline{\underline{{\cal{L}}}}}}
 \def\cMmat{{\underline{\underline{{\cal{M}}}}}}
 \def\kappamat{{\underline{\underline{{\kappa}}}}}
%
%
\def\Grad{\nabla}
\def\Div{\nabla \cdot}
\def\Lap{\nabla^2}
%

\begin{frontmatter}
\title{Entropy-Bounded Discontinuous Galerkin Scheme for Euler Equations}
\author{Yu Lv\corref{CORR1}}
\cortext[CORR1]{Corresponding author}
\ead{ylv@stanford.edu}
\author{Matthias Ihme}
\ead{mihme@stanford.edu}
\address{Department of Mechanical Engineering, Stanford University, Stanford, CA 94305, USA}
\begin{abstract}
An entropy-bounded Discontinuous Galerkin (EBDG) scheme is proposed in which the solution is regularized by constraining the entropy. The resulting scheme is able to stabilize the solution in the vicinity of discontinuities and retains the optimal accuracy for smooth solutions. The properties of the limiting operator according to the entropy-minimum principle are proofed analytically, and an optimal CFL-criterion is derived. We provide a rigorous description for locally imposing entropy constraints to capture multiple discontinuities. Significant advantages of the EBDG-scheme are the general applicability to arbitrary high-order elements and its simple implementation for two- and three-dimensional configurations. Numerical tests confirm the properties of the scheme, and particular focus is attributed to the robustness in treating discontinuities on arbitrary meshes.
\end{abstract}

\end{frontmatter}
\tableofcontents
\section{\label{INTRO}Introduction}
The stabilization of solutions near flow-field discontinuities remains an open problem to the discontinuous Galerkin (DG) community. Considerable progress has been made on the development of limiters for two-dimensional quadrilateral and triangular elements. These limiters can be categorized into three classes. Methods that limit the solution using information about the slope along certain spatial directions \cite{TVBlim, TVBlim2D} fall in the first class. The second class of limiters extends this idea by limiting based on the moments of the solution~\cite{MOMLIM1_BISWAS,MOMLIM2_KRIVO}, and schemes in which the DG-solution is projected onto a WENO~\cite{WENOLIM1_QIU,WENOLIM_ZHU_2008,WENOLIM_ZHONG} or Hermit WENO (HWENO)~\cite{HWENOLIM_LUO} representation fall in the last category. 
Although these limiters show promising results for canonical test cases on regular elements and structured mesh partitions, the following two issues related to practical applications have not been clearly answered:
\begin{itemize}
\item How can discontinuous solutions be regularized on multi-dimensional curved high-order elements?
\item How can non-physical solutions that are triggered by strong discontinuities and geometric singularities be avoided?
\end{itemize}
The present work attempts to simultaneously address both of these questions. 

Recently, positivity-preserving DG-schemes have been developed for the treatment of flow-field discontinuities, and relevant contributions are by Zhang and Shu~\cite{ZHANG_SHU_JCP2010, ZHANG_SHU_JCP2011,ZHANG_SHU_NUMERMATH2010}. The positivity preserving method provides a robust framework with provable $L_1$-stability, preventing the appearance of negative pressure and density. Resulting algorithmic modifications are minimal, and these schemes have been used in simulations of detonation systems with complex reaction chemistry~\cite{LV_IHME_PCI_2014,LV_IHME_JCP_2014}. 

Motivated by these attractive properties, the present work aims at developing an algorithm that avoids non-physical solutions on arbitrary elements and multi-dimensional spatial representations. The resulting scheme that will be developed in this work has the following properties: First, by invoking the entropy principle, solutions are constrained by a local entropy bound. Second, a general implementation on arbitrary elements is proposed without restriction to a specific quadrature rule. Third, the entropy constraint is imposed on the solutions through few algebraic operations, thereby avoiding the computationally expensive inversion of a nonlinear system. Fourth, a method for the evaluation of an optimal CFL-criterion is derived, which is applicable to general polynomial orders and arbitrary element types.

The remainder of this paper has the following structure. The governing equations and the discretization are summarized in the next two sections. The entropy-bounded DG (EBDG) formulation is presented in Sec.~\ref{SEC_ENTROPY_PRINCIPLE}, and the derivation of the CFL-constraint and the limiting operator are presented. This analysis is performed by considering a one-dimensional setting, and the generalization to multi-dimensional and arbitrary elements is presented in Sec.~\ref{SEC_GENERALIZATION}. Section~\ref{SEC_EVALUATION_ENTROPY_BOUND} is concerned with the evaluation of the entropy-bounded DG-scheme, and a detailed description of the algorithmic implementation is given in Sec.~\ref{SEC_ALG_IMPLEMENTATION}.  The EBDG-method is demonstrated by considering several test cases, and the accuracy and stability are examined in Sec.~\ref{SEC_NUMERICAL_TEST}. The paper finishes with conclusions.
\section{\label{SEC_GOVERNING_EQUATIONS}Governing equations}
We consider a system of conservation equations,
\begin{equation}
\label{EQ_GOVERNING_EQU}
\frac{\partial \mathsf{U}}{\partial t} + \nabla \cdot \mathsf{F} = 0\qquad\text{in $\Omega$}\;,
\end{equation}
where the solution variable $\mathsf{U}:\mathbb{R}\times\mathbb{R}^{N_d} \rightarrow \mathbb{R}^{N_v}$ and the flux term $\mathsf{F}:\mathbb{R}^{N_v}\rightarrow \mathbb{R}^{N_v\times N_d}$. Here, $N_d$ denotes the spatial dimension and $N_v$ is the dimension of the solution vector. For the Euler equations, $\mathsf{U}$ and $\mathsf{F}$ take the form:
\begin{subeqnarray}
\label{EQ_EULER_FLUX}
\slabel{EQ_EULER_STATE}
\mathsf{U}(x,t) &=& (\rho, \rho u, \rho e)^T\;,\\
\slabel{EQ_EULER_FLUX}
\mathsf{F}(\mathsf{U}) & = & \left(\rho u, \rho u \otimes u + p \text{I}, u (\rho e+p)\right)^T\;,
\end{subeqnarray}
where $t$ is the time, $x\in \mathbb{R}^{N_d}$ is the spatial coordinate vector, $\rho$ is the density, $u \in \mathbb{R}^{N_d}$ is the velocity vector, $e$ is the specific total energy, and $p$ is the pressure. Equation (\ref{EQ_GOVERNING_EQU}) is closed with the ideal gas law:
\begin{equation}
p = (\gamma - 1) \left(\rho e - \frac{\rho |u|^2}{2}\right)\;,
\end{equation}
in which $\gamma$ is the ratio of specific heats, which, for the present work, is set to a constant value of $\gamma=1.4$. Here and in the following, we use $|\cdot|$ to represent the Euclidean norm. With this, we define the local maximum characteristic speed as, 
\be
\nu = |u| + c\qquad\text{with}\qquad c = \sqrt{\frac{\gamma p}{\rho}}\;,
\ee
where $c$ is the speed of sound. 

Because of the presence of discontinuities in the solution of Eq.~(\ref{EQ_GOVERNING_EQU}), we seek a weak solution that satisfies physical principles. This is the so-called entropy solution. By introducing $\cU$ as a convex function of $\mathsf{U}$ with $\cU: \mathbb{R}^{N_v} \rightarrow \mathbb{R}$, Lax~\cite{LAX_ENTROPY_BOUND_1971} showed that the entropy solution of Eq.~(\ref{EQ_GOVERNING_EQU}) satisfies the following inequality:
\begin{equation}
\label{EQ_ENTROPY_GOVERNING_EQU}
\frac{\partial \mathcal{U}}{\partial t} + \nabla \cdot \mathcal{F} \leq 0\;,
\end{equation}
where $\cF: \mathbb{R}^{N_v} \rightarrow \mathbb{R}^{N_d}$ is the corresponding flux of $\mathcal{U}$. The consistency condition between Eqs.~(\ref{EQ_GOVERNING_EQU}) and~(\ref{EQ_ENTROPY_GOVERNING_EQU}) requires~\cite{LAX_ENTROPY_BOUND_1971}:
\begin{equation}
\left(\frac{\partial \mathcal{U}}{\partial \mathsf{U}} \right)^T\frac{\partial \mathsf{F}}{\partial \mathsf{U}} = \frac{\partial \mathcal{F}}{\partial \mathsf{U}}\;.
\end{equation}
The weak solution of Eq.~(\ref{EQ_GOVERNING_EQU}) that satisfies this additional condition for the pair $(\mathcal{U}, \mathcal{F})$ is called an entropy solution. With this definition,  Eq.~(\ref{EQ_ENTROPY_GOVERNING_EQU}) is commonly called entropy inequality or entropy condition, and $\mathcal{U}$ is called the entropy variable. A familiar example for gas-dynamic applications is to relate $\cU$ to the physical entropy $s$ with:
\begin{equation}
\label{ENTROPY_DEFINITION}
s = \ln(p) - \gamma \ln(\rho) + s_0\;,
\end{equation}
where $s_0$ is the reference entropy. The corresponding definition of the entropy variable and its flux in the context of the Euler system is \cite{TADMOR_1986}:
\begin{equation}
 \label{EQ_ENTROPY_DEF}
 (\mathcal{U}, \mathcal{F})=(-\rho s, -\rho s u)\;.
\end{equation}
Note that Eq.~(\ref{ENTROPY_DEFINITION}) directly provides a constraint on the positivity of pressure $p$ and density $\rho$.
\section{\label{SEC_DG_DISCRETIZATION}Discontinuous Galerkin discretization}
We consider the problem to be posed on the domain $\Omega$ with boundary $\partial \Omega$. A mesh partition is defined as $\Omega = \cup_{e=1}^{N_e} \Omega_e$, where $\Omega_e$ corresponds to a discrete element of this partition. The edge of element $\Omega_e$ is defined as $\partial \Omega_e$. In order to distinguish different sides of the edge, the superscripts ``$+$" and ``$-$" are used to denote the interior and exterior, respectively. We define a global space of test functions as
\begin{equation}
\mathpzc{V} = \oplus_{e=1}^{N_e} \mathpzc{V}_e\;,\qquad\mathpzc{V}_e = \text{span}\{\varphi_n(\Omega_e)\}_{n=1}^{N_p}\;,
\end{equation}
where $\varphi_n$ is the $n$th polynomial basis, and $N_p$ is the number of bases. On the space $\mathpzc{V}_e$ we seek an approximate solution to Eq.~(\ref{EQ_GOVERNING_EQU}) of the form:
\begin{equation}
\mathsf{U} \simeq U = \oplus_{e=1}^{N_e} U_e\;,\qquad U_e \in \mathpzc{V}_e\;,
\end{equation}
where the solution vector $U_e$ on each individual element takes the general form
\begin{equation}
U_e(x,t) = \sum_{m=1}^{N_p} \widetilde{U}_{e,m}(t) \varphi_m(x)\;,
\end{equation}
and the unknown vector of basic coefficients $\widetilde{U}_{e,m} \in \mathbb{R}^{N_v \times N_p}$ is obtained from the discretized weak solution of Eq.~(\ref{EQ_GOVERNING_EQU}):
\begin{equation}
 \label{WEAK_FORM}
 \frac{d \widetilde{U}_{e,m}}{d t}\int_{\Omega_e}\varphi_n \varphi_m d\Omega - \int_{\Omega_e} \nabla\varphi_n \cdot F(U_e)d\Omega + \int_{\partial \Omega_e}  \varphi^+_n \wh{F}(U^{+}_e, U^{-}_e, \wh{\textmd{n}}) d\Gamma = 0\;,
\end{equation}
$\forall \varphi_n$ with $n=1, \ldots,N_p$. The numerical Riemann flux $\wh{F}$ is evaluated based on the states at both sides of the interface $\partial \Omega_e$ and the outward-pointing normal vector $\wh{\textmd{n}}$. It is of interest to note that for the particular case of $N_p=1$ and $\varphi_1 = 1$, the weak form reduces to the classical first-order finite-volume (FV) discretization. It can also be seen that the DG-scheme does not rely on a specific type of basis functions. Since the following derivation is based on this mathematical property, we introduce the following lemma.
\begin{lemma}
\label{LEMMA_BASIS_SWITCH}
A polynomial $P$,
\be
\label{EQ_LM_POLYNOMIAL1}
P(x) = \sum_{m=1}^{N_p} \widetilde{P}_m \varphi_m(x)\qquad \text{for $x\in \Omega_e$}\;,
\ee
with a set of polynomial bases $\{\varphi_m(x), m = 1,\ldots,N_p\}$, can be exactly interpolated by a Lagrangian polynomial of $N_p$ points $\{y_n \in \Omega_e,~n = 1,\ldots,N_p\}$ under the condition that $\big[\varphi_m(y_n)\big]$ is non-singular: 
\be
\label{EQ_LM_POLYNOMIAL2}
P(x) = \sum_{n=1}^{N_p} P(y_n)\phi_n(x)\qquad \text{for $x\in \Omega_e$}\;.
\ee
\\
Proof: By equating Eqs.~(\ref{EQ_LM_POLYNOMIAL1}) and (\ref{EQ_LM_POLYNOMIAL2}), and comparing terms, it follows that
\be
\label{EQ_BASIS_SWITCHING_TENSOR}
\varphi_m(x) = \sum_{n=1}^{N_p} \varphi_m(y_n) \phi_n(x) \qquad \text{or}\qquad [\varphi_m(x)] =  [\varphi_m(y_n)] ~ [\phi_n(x)]\;,
\ee
where we use $[\cdot]$ to denote a tensor or a vector. Since $[\varphi_m(y_n)]$ is non-singular, Eq.~(\ref{EQ_BASIS_SWITCHING_TENSOR}) can be inverted:
\ben
[\phi_n(x)] =  [\varphi_m(y_n)]^{-1} ~ [\varphi_m(x)]\;.
\een
\end{lemma}
\begin{remark}
The significance of this lemma is that it provides a description to convert any basis set to a Lagrangian basis set with $N_p$ interpolation points, as long as they are located at general positions. To facilitate the following derivation, we choose the points $y_n$ with $ n=1,\ldots,N_p$ from the $N_q$ quadrature points\cite{CUBATURE_FORMULAS}. According to the accuracy requirement of the quadrature scheme for Eq.~(\ref{WEAK_FORM}), $N_q \geq N_p$ is alway true.
\end{remark}

\section{\label{SEC_ENTROPY_PRINCIPLE}Entropy principle and entropy-bounded discontinuous Galerkin method}
In this section, we review the entropy principle by considering a three-point FV-setting. Then, we will explore how to extend this principle to a DG-scheme, which leads to the concept of entropy boundedness. In order to enable the implementation of this concept, two important ingredients will be discussed, namely a time-step constraint and a limiting operator.  After conducting numerical analyses by considering a one-dimensional configuration, we will extend the entropy boundedness to multi-dimensional and arbitrary element types. The dimensional generality, geometric adaptability and simple implementation are major advantages of the resulting entropy-bounded DG-method.

\subsection{\label{ENTROPY_PRINCIPLE_THREE-POINT_FVM}Preliminaries and related work}
To illustrate the entropy principle, we consider a local Lax-Friedrichs flux, which can be written as:
\begin{equation}
\label{LAX_FRIEDRICHS_FLUX}
  \wh{F}(U_L, U_R, \wh{\textmd{n}}) = \frac{1}{2}\left(F(U_L) + F(U_R)\right) \cdot \wh{\textmd{n}} - \frac{1}{2}\lambda(U_R - U_L)\;,
\end{equation}
and
\ben
  \lambda \geq \max\limits_{k \in \{L, R\}}  \nu (U_k)
\een
is the dissipation coefficient. Note that this flux function satisfies consistency: $\wh{F}(U, U, \wh{\textmd{n}}) = F(U) \cdot \wh{\textmd{n}}$, conservation: $\wh{F}(U_L, U_R, \wh{\textmd{n}}) = - \wh{F}(U_R, U_L, -\wh{\textmd{n}})$, and Lipschitz-continuity. In the following, we consider the simplest case of DGP0 scheme, with $N_p = 1$, in a one-dimensional setting. This formulation is consistent with the classical three-point FV-discretization. For $x \in \Omega_e = [x_{\enh},~ x_{\eph}]$, the discretized solution to Eq.~(\ref{EQ_GOVERNING_EQU}) can be written as:
\be
 \label{1D_FVM}
 \nonumber
\widetilde{U}_e(t+\Delta t) & = &  \widetilde{U}_e -\frac{\Delta t}{h}\left(\wh{F}(\widetilde{U}_e, \widetilde{U}_{e-1}, -1) + \wh{F}(\widetilde{U}_e, \widetilde{U}_{e+1}, 1) \right)\;,
\ee
where $\widetilde{U}_e$ is the basis coefficient, which is identical to the piecewise constant approximation to the exact solution in $\Omega_e$. In the following, we introduce $\widetilde{U}_e^{\Delta t}$ to denote the solution vector $\widetilde{U}_e(t+\Delta t)$, and use the superscript $\Delta t$ to denote a temporally updated quantity at $t+\Delta t$. With the numerical flux given in Eq.~(\ref{LAX_FRIEDRICHS_FLUX}), this discretization preserves the positivity of pressure and density under the CFL-condition \cite{ZHANG_SHU_JCP2010, PERTHAME_SHU_POSI_1996}:
\begin{equation}
 \label{1D_CFL_CONDI}
 \frac{\Delta t \lambda}{h} \leq \frac{1}{2}\;.
\end{equation}
In addition, it was discussed in \cite{PERTHAME_SHU_POSI_1996} that Eq.~(\ref{1D_FVM}) satisfies the discrete minimum entropy principle proposed by Tadmor~\cite{TADMOR_1986},
\begin{equation}
\label{TADMOR_ENTROPY_PRIN}
  s(\widetilde{U}_e^{\Delta t}) \geq {s}^0_e(t) = \min\limits_{j\in \{e-1, e, e+1\}}s(\widetilde{U}_j) .
\end{equation}
To show this property, we can rewrite Eq.~(\ref{1D_FVM}) and split $\widetilde{U}_e(t + \Delta t)$ into two parts. For $x \in \Omega_e$, this is written as
\begin{subeqnarray}
\label{DISCRET_SPLIT}
\slabel{DISCRET_SPLIT_I}
\widetilde{U}_e(t + \Delta t) & = & \f{1}{2}\left(\widetilde{U}_{e,p1}^{\Delta t} + \widetilde{U}_{e,p2}^{\Delta t}  \right)\;,\\
\slabel{DISCRET_SPLIT_P1}
\widetilde{U}^{\Delta t}_{e,p1} &=& \widetilde{U}_e -\frac{\Delta t}{h}\left( F(\widetilde{U}_{e+1}) - \lambda_{\eph}\widetilde{U}_{e+1}  - F(\widetilde{U}_{e}) + \lambda_{\eph}\widetilde{U}_{e} \right)\;,\\
\slabel{DISCRET_SPLIT_P2}
\widetilde{U}^{\Delta t}_{e,p2} &=& \widetilde{U}_e + \frac{\Delta t}{h}\left(F(\widetilde{U}_{e-1}) + \lambda_{\enh}\widetilde{U}_{e-1} - F(\widetilde{U}_{e}) - \lambda_{\enh}\widetilde{U}_{e} \right)\;,
\end{subeqnarray}
where $\widetilde{U}_{e,p1}^{\Delta t}$ and $\widetilde{U}_{e,p2}^{\Delta t}$ can be viewed as the P0-approximations to the solutions of the hyperbolic systems (under the CFL constraint of Eq. (\ref{1D_CFL_CONDI}))
\begin{subeqnarray}
\label{SYSTEM_INEQ_MODIFIED}
\slabel{SYSTEM_INEQ_MODIFIED_RIGHT}
 \frac{\partial \mathsf{U}}{\partial t}  + \left(\mathsf{F}'(\mathsf{U}) -  \lambda_{\eph} \rm{I} \right) \frac{\partial \mathsf{U}}{\partial x}= 0\;,\\
\slabel{SYSTEM_INEQ_MODIFIED_LEFT}
 \frac{\partial \mathsf{U}}{\partial t}  + \left(\mathsf{F}'(\mathsf{U}) +  \lambda_{\enh}\rm{I}\right)   \frac{\partial \mathsf{U}}{\partial x } = 0\;,
\end{subeqnarray}
with the exact (Godunov) flux. If we denote the exact solutions to Eq.~(\ref{SYSTEM_INEQ_MODIFIED_RIGHT}) and (\ref{SYSTEM_INEQ_MODIFIED_LEFT}) as $\mathsf{U}_{p1}(x, t + \Delta t)$ and $\mathsf{U}_{p2}(x, t + \Delta t)$, respectively, then their P0-approximations in $\Omega_e$ yield $\widetilde{U}_{e,p1}^{\Delta t} = \frac{1}{h} \int_{x\enh}^{x\eph} \mathsf{U}_{p1}(x, t + \Delta t) dx$ and $\widetilde{U}_{e,p2}^{\Delta t} = \frac{1}{h} \int_{x\enh}^{x\eph} \mathsf{U}_{p2}(x, t + \Delta t) dx$. Both equation systems are obtained by imposing a constant shift on the characteristic speeds without modifying the characteristic variables. With these modifications, all characteristics in Eq.~(\ref{SYSTEM_INEQ_MODIFIED_RIGHT}) are right-running while those in Eq.~(\ref{SYSTEM_INEQ_MODIFIED_LEFT}) are left-running. The corresponding entropy inequalities take the form
\begin{subeqnarray}
\label{ENTROPY_INEQ_MODIFIED}
\slabel{ENTROPY_INEQ_MODIFIED_RIGHT}
\frac{\partial \mathcal{U}}{\partial t} + \frac{\p}{\p x} \left(\mathcal{F} - \lambda_{\eph}\mathcal{U}  \right) \leq 0\;,\\
\slabel{ENTROPY_INEQ_MODIFIED_LEFT}
\frac{\partial \mathcal{U}}{\partial t} + \frac{\p}{\p x} \left(\mathcal{F} + \lambda_{\enh}\mathcal{U}  \right)\leq 0\;.
\end{subeqnarray}

Without loss of generality, we now consider Eq. (\ref{ENTROPY_INEQ_MODIFIED_RIGHT}) and integrate over $[t,t+\Delta t]\times[x\enh,x\eph]$, resulting in the following expression:
\begin{equation}
\begin{split}
\int\limits_{x\enh}^{x\eph}\mathcal{U}\big(\mathsf{U}_{p1}(x, t+{\Delta t})\big) dx - 
\int\limits_{x\enh}^{x\eph}  \mathcal{U}(\widetilde{U}_e)dx + 
\int\limits_t^{t+\Delta t}\left(\mathcal{F}(\mathsf{U}(x_{\eph},t)) -   \lambda_{\eph}\mathcal{U}(\mathsf{U}(x_{\eph},t)) \right) dt \\
 -  \int\limits_t^{t+\Delta t}\left(\mathcal{F}(\mathsf{U}(x_{\enh},t)) -  \lambda_{\enh}\mathcal{U}(\mathsf{U}(x_{\enh},t)) \right) dt \leq 0\;.
\end{split}
\end{equation}
Recognizing that all characteristics are right-running, the temporal integral can be evaluated exact since $\mathsf{U}(x_{\enh},t) = \widetilde{U}_e$ and $\mathsf{U}(x_{\eph},t) = \widetilde{U}_{e+1}$ under the condition of Eq.~(\ref{1D_CFL_CONDI}). Then by utilizing the convexity of $\mathcal{U}$ with respect to $\mathsf{U}$, the following estimate is obtained:  
\ben
\mathcal{U}(\widetilde{U}_{e, p1}^{\Delta t}) & = & \mathcal{U}\left(\frac{1}{h}\int_{x\enh}^{x\eph}\mathsf{U}_{p1}(x, t+\Delta t)dx \right) \leq \frac{1}{h}\int_{x\enh}^{x\eph}\mathcal{U}(\mathsf{U}_{p1}(x, t+\Delta t)) dx\;\\
& \leq &  \mathcal{U}(\widetilde{U}_e) + \frac{\Delta t}{h}\left(\mathcal{F}(\widetilde{U}_{e}) - \lambda_{\eph} \mathcal{U}(\widetilde{U}_{e}) \right) - \frac{\Delta t}{h}\left( \mathcal{F}(\widetilde{U}_{e+1}) -  \lambda_{\eph} \mathcal{U}(\widetilde{U}_{e+1}) \right)  \;.
\een
With the definition of $(\cU, \cF)$, given in Eq.~(\ref{EQ_ENTROPY_DEF}), it follows
\be
s(\widetilde{U}_{e,p1}^{\Delta t}) \geq \frac{\rho_e}{\rho_{e,p1}^{\Delta t}} \left[1- \frac{\Delta t}{h}\left(\lambda_{\eph} - u_e\right)\right] s(\widetilde{U}_e) + \frac{\rho_{e+1}}{\rho_{e,p1}^{\Delta t}} \frac{\Delta t}{h}\left(\lambda_{\eph} - u_{e+1} \right) s(\widetilde{U}_{e+1})\;.
\ee
The constraint (\ref{1D_CFL_CONDI}) ensures that the coefficients in front of $s(\widetilde{U}_{e})$ and $s(\widetilde{U}_{e+1})$ are positive and sum to unity according to Eq. (\ref{DISCRET_SPLIT_P1}). From these arguments directly follows:
\begin{equation}
s(\widetilde{U}_{e,p1}^{\Delta t}) \geq \min\{s(\widetilde{U}_e), s(\widetilde{U}_{e+1})\}\;,
\end{equation}
and
\begin{equation}
s(\widetilde{U}_{e,p2}^{\Delta t}) \geq \min\{s(\widetilde{U}_{e-1}), s(\widetilde{U}_{e})\}\;.
\end{equation}
Combining these two relations with the quasi-concavity of the entropy $s$ (Lemma 2.1 of \cite{ZHANG_SHU_NUMERMATH2010}), the discrete minimum entropy principle of Eq. (\ref{TADMOR_ENTROPY_PRIN}) is obtained. 

The result above is obtained for a one-dimensional setting. However, in the following we derive a rotational version of this entropy principle for multi-dimensional cases by following the idea of~\cite{PERTHAME_SHU_POSI_1996}. By considering an arbitrary face with a normal $\wh{\textmd{n}}$, we can define a tangential vector $\wh{\textmd{t}}$. By projecting the velocity in Cartesian coordinates onto this local coordinate $(\wh{\textmd{n}},~ \wh{\textmd{t}})$ with the following mapping,
\ben
u ~\rightarrow ~(u_n,~u_t )^T\;,
\een
with 
\ben
u_n = u \cdot \wh{\textmd{n}}, \qquad~ u_t = u \cdot \wh{\textmd{t}}\;,
\een
with which, the conservative variables and flux terms can be written as follows
\begin{subeqnarray}
\label{3D_AUG_SYS}
\slabel{3D_AUG_SYS_U}
\mathsf{U} &=& (\rho,~\rho u_n,~\rho u_t,~\rho e_{n},~\rho e_{t})^T,\;\\
\slabel{3D_AUG_SYS_F}
\mathsf{F} &=& (\rho u_{n},~\rho u_{n}^2 + p,~\rho  u_{t} u_{n},~ u_{n}(\rho e_{n} + p ),~\rho u_{n}e_{t})^T\;,
\end{subeqnarray}
and
\ben
e_{n} = \frac{p}{\rho(\gamma - 1)} + \frac{1}{2}u_{n}^2, \qquad~e_{t} = \frac{1}{2} u_{t}^2.
\een
For this augmented system, it can be seen that the variations of density, pressure and entropy are all governed by a 1D reduced system that is parallel to $\wh{\textmd{n}}$:
\begin{subeqnarray}
\label{3D_RED_SYS}
\slabel{3D_RED_SYS_U}
\mathsf{U} &=& (\rho,~\rho u_n,~\rho e_n)^T,\;\\
\slabel{3D_RED_SYS_F}
\mathsf{F} &=& (\rho u_n,~\rho u_n^2 + p,~ u_n(\rho e_n + p ))^T,\;
\end{subeqnarray}
and its discretized version is identical to Eq. (\ref{1D_FVM}). Thus, the solution preserves the positivity of density and pressure and satisfies the entropy principle under the CFL-constraint, Eq. (\ref{1D_CFL_CONDI}), with $\lambda \geq \max\limits_{j \in \{e-1, e, e+1\}} \nu(U_j)$.  

We now conclude the above analysis with the following lemma as a critical element for the subsequent derivation.

\begin{lemma}
\label{LEMMA_THREEPOINT_ROT}
For a three-point system defined on ${\mathbb{R}}^{N_d}$, the solution along an arbitrary direction {$\wh{\textmd{n}}$},
\be
\label{ROT_FVM_DISC}
\widetilde{U}_e^{\Delta t} =  \widetilde{U}_e + \frac{\Delta t}{h}\left(\wh{F}(\widetilde{U}_e, \widetilde{U}_{e-1}, - \wh{\textmd{n}}) + \wh{F}(\widetilde{U}_e, \widetilde{U}_{e+1}, \wh{\textmd{n}})\right)\;, 
\ee
with the flux function $\wh{F}$ specified in Eq. (\ref{LAX_FRIEDRICHS_FLUX}), preserves the positivity of density and pressure, and satisfies the entropy principle:
\be
\label{ROT_ENTROPY_PRIN}
s(\widetilde{U}_e^{\Delta t}) \geq \min\limits_{j \in \{e-1, e, e+1\}} s(\widetilde{U}_j)\;,
\ee
under the CFL condition:
\ben
\frac{\Delta t \lambda}{h } \leq \frac{1}{2},~\qquad \lambda \geq \max\limits_{j \in \{e-1, e, e+1\}} \nu(\widetilde{U}_j),
\een
\end{lemma}
This three-point system is consistent with that used by Zhang and Shu~\cite{ZHANG_SHU_JCP2010, ZHANG_SHU_NUMERMATH2010}. The difference is that we introduce a local entropy bound $s_e^0$ at time $t$ instead of a global entropy bound that is derived from the initial conditions $\min\limits_{x\in\Omega} s(x, 0)$. Although a local Lax-Friedrichs flux was used for illustrative purposes, other Riemann solvers that preserve positivity and entropy stability are equally suitable, for example, the Roe-type solver with entropy fix \cite{ROE_ENTROPY_FIX_TADMOR}, the kinetic-type solver \cite{KINE_ENTROPY_SOLVER_1994} and the exact Godunov solver.

\subsection{\label{SEC_ENTROPY_BOUNDED_DG}Entropy-bounded DG-scheme}

To robustly capture shocks while retaining the high-order benefit of the DG-scheme, sub-cell shock resolution is required \cite{ZJWANG_REVIEW_2013}. We now extend the discussion by considering a high-order DG-solution with sub-cell representation. In each DG-cell, the whole solution is approximated by a function space. However, there is no guarantee that the high-order ($N_p > 1$) solution obeys the physical entropy principle. This is the reason that DG suffers from numerical instability in the vicinity of discontinuities. To suppress these instabilities, one approach is to consider imposing constraints based on the behavior of the entropy solution. The positivity-preserving DG-method \cite{ZHANG_SHU_JCP2010, ZHANG_SHU_JCP2011} is a successful example for this approach. Based on the entropy principle, Eq.~(\ref{TADMOR_ENTROPY_PRIN}), Zhang and Shu \cite{ZHANG_SHU_NUMERMATH2010} extended their implementation to an entropy-based constraint. Here, we propose a general framework that is based on the entropy principle, and major differences and advantages have been highlighted in Sec. \ref{INTRO}. 

We define the constraint for the high-order DG-scheme as follows:
\be
\label{ENTROPY_PRINCIPLE_FOR_DG}
\forall x \in \Omega_e,~~s(U_e^{\Delta t}(x)) \geq \min\{s(U(y))|~ y \in \Omega_e \cup \partial \Omega_e^-\} \equiv s_e^0(t)\;.
\ee
In this equation, the right-hand-side sets an entropy bound for an element-local solution in $\Omega_e$; with this, we refer to a DG-solution as \emph{entropy-bounded} if it satisfies this principle. $s_e^0(t)$ is a local estimate for the true entropy minimum in $\Omega_e$, $|s_e^0(t) - \min\limits_{x \in \Omega_e} s\left(\mathsf{U}(x)\right)| \sim \mathcal{O}(h^k)$, where $k$ is the local order of accuracy. Besides that, $s_e^0(t)$ is bounded if the entropy is bounded at the domain boundaries, $s_e^0(t)\geq \min\limits_{x\in \Omega}s(U(x, t=0)) = s^0$, where $s^0$ is the minimum entropy at the initial condition. By imposing this constraint, we expect that the sub-cell DG-solution is regularized, avoiding the appearance of non-physical solutions. This idea is illustrated in Fig. \ref{EB_IDEA_DEMO}. At time level $t$, $s_e^0(t)$ is calculated and used to set a reference bound for the solution at the next step, $U_e^{\Delta t}$. If $U_e^{\Delta t}$ yields entropy undershoot with respect to $s_e^0(t)$, it will be modified to satisfy the constraint of Eq.~(\ref{ENTROPY_PRINCIPLE_FOR_DG}). In order to implement this regularization for a high-order DG scheme, the following aspects require addressing:
\begin{itemize}
\item[i] To impose Eq.~($\ref{ENTROPY_PRINCIPLE_FOR_DG}$) on the DG-solution, we introduce a limiting operator $\cL$. The regularized solution, denoted by $^\mathcal{L}U_e^{\Delta t}$, requires that $s(^\mathcal{L}U_e^{\Delta t}(x)) \geq  s_e^0(t) $ $\forall x \in \Omega_e$. In the following, we relax this condition, and impose Eq.~(\ref{ENTROPY_PRINCIPLE_FOR_DG}) only on the set of quadrature points, $\mathcal{D}$, that are involved in solving the weak form in Eq. (\ref{WEAK_FORM}). 

\begin{figure}
\centering
\includegraphics[width=0.5\textwidth, clip=, keepaspectratio]{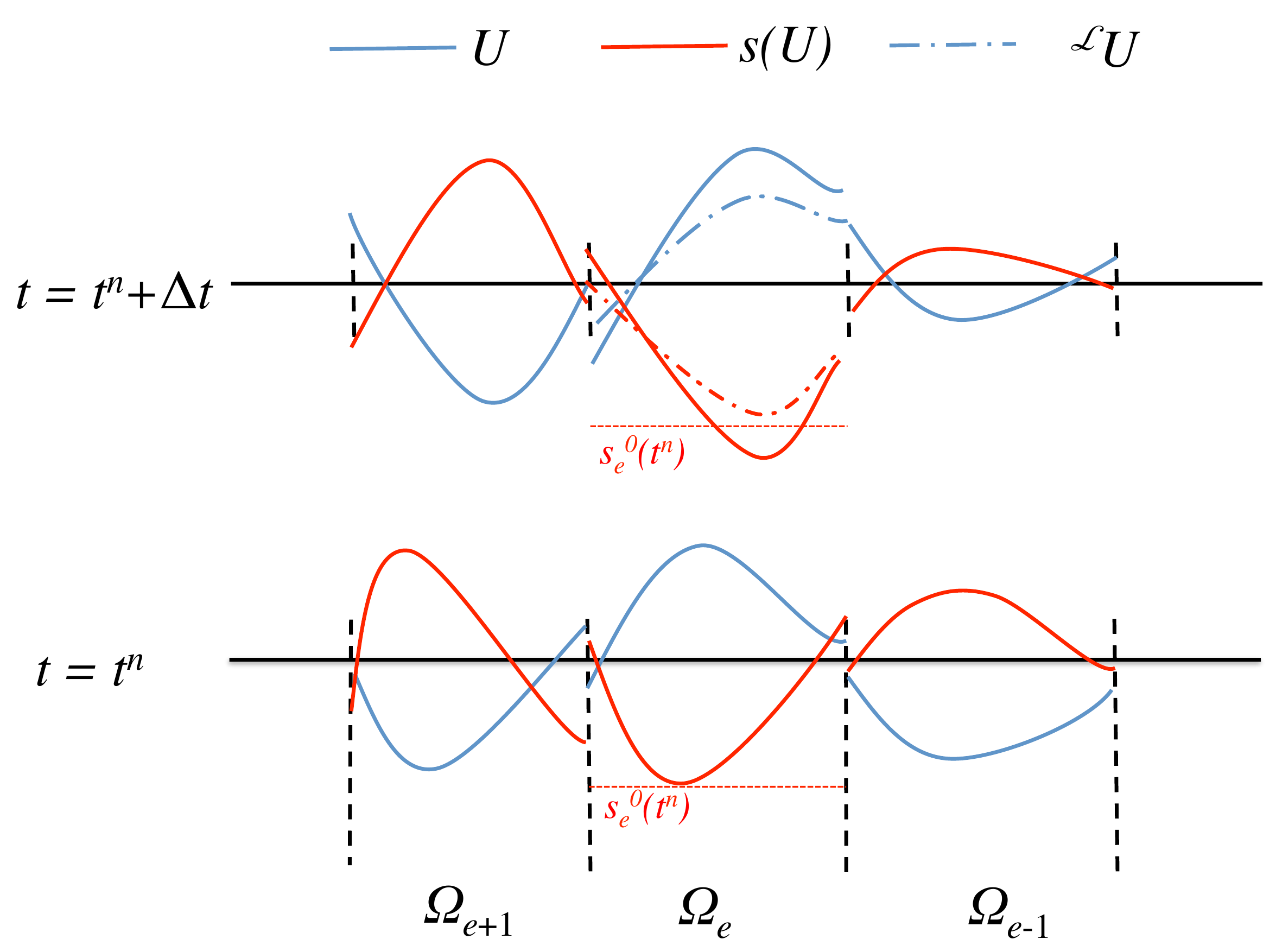}
    \caption{\label{EB_IDEA_DEMO}(Color online) Schematic of entropy-bounding of the EBDG scheme.}
\end{figure}

\item[ii] Guaranteeing that the constraint (\ref{ENTROPY_PRINCIPLE_FOR_DG}) is always imposed requires the existence of the operator $\cL$. A sufficient condition for this is that the element-averaged solution is entropy-bounded, $s(\ol{U}_e^{\Delta t}) \geq s_e^0(t)$. Enforcing this condition relies on the selection of a proper CFL-condition, and this analysis will be developed in Sec.~\ref{DERIVE_CFL_1DDG} for a one-dimensional system. Subsequently, this analysis is then extended in Sec.~\ref{SUBSEC_CFL_CONSTRAINT_ARB_ELEMENT}  to general multi-dimensional elements.

\item[iii] Algorithmic details on the implementation of the operator $\cL$ constraining the element-local DG-solution are discussed in Sec.~\ref{1D_LIMITER_DESIGN}.

\item[iv] The evaluation of the lower bound $s_e^0(t)$ that is necessary to constrain the entropy solution is given in Sec.~\ref{SEC_EVALUATION_ENTROPY_BOUND}.
\end{itemize}

\subsection{\label{DERIVE_CFL_1DDG}CFL-constraint for one-dimensional entropy-bounded DG}
The objective now is to extend the analysis for DGP0 to a DG-scheme with high-order polynomial representations. Consider a one-dimensional domain in which the element $\Omega_e$ is centered at $x_e$, and a quadrature rule with weights $w_q$ and $\sum_{q=1}^{N_q}w_q = 1$. These quadrature weights are evaluated at the quadrature points $x_q \in [x\eph,~x\enh]$. The discretized cell-averaged solution $\ol{U}_e$ is defined as:
\be
\ol{U}_e=\f{1}{h}\int_{\Omega_e} U_e d x,
\ee
(for the P0-case discussed above, $\ol{U}_e = \wt{U}_e$), which can be further expanded by a quadrature rule with sufficient accuracy:
\begin{equation}
\begin{split}
\overline{U}_e 
\label{1DDG_CELL_AVE}
 = & \sum\limits_{q=1}^{N_q}w_qU_e(x_q)\;,\\
 = & \sum_{q=1}^{N_q}\left({w_q} - \theta_l\phi_q(x\enh) - 
     \theta_r \phi_q(x\eph)\right)U_e(x_q) + \theta_l U_e(x\enh)+ \theta_rU_e(x\eph)\;,\\
 = & \sum_{q=1}^{N_q}\theta_q U_e(x_q)+ \theta_l U_e(x\enh)+ \theta_rU_e(x\eph)\;,
\end{split}
\end{equation}
where the first line utilizes the exactness of the quadrature rule, the second line utilizes Lemma~\ref{LEMMA_BASIS_SWITCH}, and the third line defines $\theta_q =w_q - \theta_l\phi_q(x\enh) - \theta_r \phi_q(x\eph)$. Under the condition that $\theta_{r,l} > 0$ and $\theta_q \geq 0$, the last line of Eq.~(\ref{1DDG_CELL_AVE}) is a convex combination. Since the quadrature weights $w_q$ are positive, the existence of $\theta_{r,l}$ is guaranteed through the condition $w_q\geq \theta_l\phi_q(x\enh) + \theta_r \phi_q(x\eph)$. If $\phi_q(x_{e\pm {1}/{2}})>0$, $\theta_{r,l}$ is constrained as $(0, \min_q{w_q/\max\{\phi_q(x_{e\pm {1}/{2}})\}}]$. If some of  $\phi_q(x_{e\pm {1}/{2}})$ are negative, they are not essential in setting the upper bound for $\theta_{r,l}$.

\begin{remark}
In the following, $\theta_{r,l}$ will be related to a CFL-constraint. To obtain an optimal CFL-number, the largest value of  $\theta_{r,l}$ needs to be found. This can be formulated as a maximization problem subject to the constraints, $\theta_{r,l} > 0$ and $\theta_q \geq 0$.
\end{remark}

For illustration, we fully discretize Eq.~(\ref{WEAK_FORM}) using a forward Euler time integration scheme and insert the results from Eq. (\ref{1DDG_CELL_AVE}). The element-averaged solution in $\Omega_e$ is then updated as:
\begin{eqnarray}
\label{DG_MEAN_SOLUTION}
\overline{U}_e^{\Delta t} 
& = &\overline{U}_e  - \frac{\Delta t}{h}\left(\wh{F}(U_e(x\enh), U_{e-1}(x\enh), - 1) + \wh{F}(U_e(x\eph), U_{e+1}(x\eph), 1)\right)\;,\nonumber\\
& = &  \sum_{q=1}^{N_q}\theta_q U_e(x_q)  + 
\nonumber\\
&& \theta_l U_e(x\enh) - \frac{\Delta t}{h}\left(\wh{F}(U_e(x\enh), U_{e-1}(x\enh), -1) + \wh{F}(U_e(x\enh), U_e^*, 1)\right) + \nonumber
\\
&& \theta_r U_e(x\eph) - \frac{\Delta t}{h}\left(\wh{F}(U_e(x\eph), U_e^*, -1) + \wh{F}(U_e(x\eph),U_{e+1}(x\eph), 1)\right)\;,
\end{eqnarray}
where
\be
\label{EXP_U_STAR}
U_e^* = \frac{1}{2}\left(U_e(x\enh) + U_e(x\eph)\right)- \frac{1}{2\lambda} \left( F(U_e(x\eph)) -  F(U_e(x\enh)) \right)
\ee
is introduced to simplify subsequent analyses. Note that $U_e^*$ ensures the validity of the second equality in Eq. (\ref{DG_MEAN_SOLUTION}) with a $\lambda$ that is defined in the following lemma. We can see that Eq. (\ref{DG_MEAN_SOLUTION}) contains two three-point systems discussed in Sec.~\ref{ENTROPY_PRINCIPLE_THREE-POINT_FVM}. To guarantee that $\overline{U}_e^{\Delta t}$ is entropy-bounded, it is necessary that these systems conform to the entropy principle of Eq.~(\ref{TADMOR_ENTROPY_PRIN}). This leads to the following lemma. 

\begin{lemma} 
\label{ENTROPY_THREE_ELEMENT_DG}
For a one-dimensional DG-system, the element-averaged solution satisfies the entropy principle
\begin{equation}
  \label{MEAN_ENTROPY_BOUNDED_1D_DG}
  s(\overline{U}_e^{\Delta t})   \geq s^0_e(t) = \min\{s(U(y))|~ y\in \Omega_e \cup \partial\Omega_e^- \}\;,
\end{equation}
under the CFL-constraint
\begin{equation}
\label{CFL_CONDITION_1D_DG}
\frac{\Delta t \lambda}{h}  \leq \frac{1}{2} \min \left\{\theta_l, \theta_r\right\},
\end{equation}
where $\lambda$ is the maximum wave speed that is evaluated over the set of point-wise solutions
\be
\lambda \geq \max_{U}~\nu(U)\qquad\text{with}\qquad  U \in \{U_{e-1}(x\enh), U_{e}(x_{e\pm1/2}), U_{e+1}(x\eph)\}\;,
\ee
given the conditions $\theta_{r,l} > 0 $ and $\theta_q \geq 0$.
\\
Proof: First, we need to show that $U_e^*$ satisfies the discretized entropy principle of Eq.~(\ref{TADMOR_ENTROPY_PRIN}). This is indeed the case since $U_e^*$ is essentially the left-hand-side of Eq. (\ref{DISCRET_SPLIT_P1}) with time step taking $\Delta t = h/(2\lambda)$, corresponding to the upper bound of the CFL-constraint. Since $U_e^*$ is an entropy solution, $\nu(U_e^*)$ is bounded by $\lambda$, which makes the Lax-Friedrichs flux $\wh{F}$ involving $U_e^*$ in Eq. (\ref{DG_MEAN_SOLUTION}) valid according to the definition in Eq.~(\ref{LAX_FRIEDRICHS_FLUX}). Therefore, we have 
\ben
s(U_e^*) &\geq& \min\left\{s(U_e(x\enh)), s(U_e(x\eph))\right\}\;,
\een

Second, we reformulate Eq. (\ref{DG_MEAN_SOLUTION}) as
\be
\label{DG_MEAN_SOLUTION_REF}
\overline{U}_e^{\Delta t} =\sum_{q=1}^{N_q}\theta_q U_e(x_q) + \theta_l \overline{U}_{e,p1}^{\Delta t} + \theta_r \overline{U}_{e,p2}^{\Delta t}\;,
\ee
in which $\overline{U}_{e,p1}^{\Delta t}$ and $\overline{U}_{e,p2}^{\Delta t}$ are the two updated solutions of the three-point system. Their definitions are readily obtained by comparing Eqs. (\ref{DG_MEAN_SOLUTION_REF}) and (\ref{DG_MEAN_SOLUTION}). The given constraints, $\theta_{r,l} > 0 $ and $\theta_q \geq 0$, guarantee that the form of the convex combination in Eq. (\ref{DG_MEAN_SOLUTION_REF}) always holds. According to  Lemma \ref{LEMMA_THREEPOINT_ROT}, it follows
\ben
s(\overline{U}_{e,p1}^{\Delta t}) &\geq& \min\left\{s(U_e^*), s(U_e(x\enh)),   s(U_{e-1}(x\enh)) \right\}\;,\\
s(\overline{U}_{e,p2}^{\Delta t}) &\geq& \min\left\{s(U_e^*), s(U_e(x\eph)),   s(U_{e+1}(x\eph))\right\}\;,
\een
under the given CFL-constraint, Eq. (\ref{CFL_CONDITION_1D_DG}). Combining this with the quasi-concavity of entropy~\cite{ZHANG_SHU_NUMERMATH2010}, it follows
\ben
s(\overline{U}_e^{\Delta t}) &\geq& \min\left\{s(U_e(x_q)),s(\overline{U}_{e,p1}^{\Delta t}), s(\overline{U}_{e,p2}^{\Delta t}) \right\}\;,\\
 & \geq & \min\{s(U(y))|~ y\in \Omega_e \cup \partial \Omega_e^- \}\;.
\een
\end{lemma}
\begin{remark}
Equation (\ref{CFL_CONDITION_1D_DG}) ensures the positivity of $p(\overline{U}_e^{\Delta t})$ and $\rho(\overline{U}_e^{\Delta t}).$
\end{remark}
In this context, we emphasize that the CFL-constraint (\ref{CFL_CONDITION_1D_DG}) provides a description for the entropy boundedness and does not conflict with the general CFL-constraint for linear stability, CFL$^{\rm L}$. To distinguish both constraints, here we use CFL$^{\rm EB}$ to denote the CFL-number for guaranteeing the entropy boundedness. In general, the time step has to be selected to satisfy both criteria. Equation~(\ref{CFL_CONDITION_1D_DG}) shows that CFL$^{\rm EB}$ depends on the value $\min \{\theta_l, \theta_r\}$, and a rigorous evaluation for this will be given below. 

Although we considers the specific case of a forward Euler time discretization scheme, all the derivation and conclusions are directly applicable to any explicit Runge-Kutta (RK) methods with positive coefficients, since the RK-solution is a convex combination of solutions obtained from several forward Euler sub-steps. In practice, RK-methods are preferred as DG time-integration schemes due to their compatible stability properties \cite{COCKBURN_SHU_JSC01}. 

\subsection{\label{1D_LIMITER_DESIGN}Construction of a limiting operator $\mathcal{L}$}
Following Lemma~\ref{ENTROPY_THREE_ELEMENT_DG}, the entropy constraint is imposed on the set of quadrature points, $x \in \mathcal{D}\subset \Omega_e$. For the one-dimensional case, $\cD$ is:
\be
\mathcal{D} =  \{x_{e\pm1/2}, x_q, q=1,\ldots,N_q ~(N_q \geq N_p) \}\;.
\ee
In the following, we are concerned with the construction of a limiting operator $\cL$, such that 
\be
\label{EQ_ENTROPY_1D_DG_APPL}
\forall x \in \mathcal{D}, \qquad s(^{\cL}U_e^{\Delta t}(x)) \geq s^0_e(t) \;, 
\ee
Since the operator $\cL$ is applied at the end of each sub-iteration, we will omit the superscript $\Delta t$ in the subsequent analysis. According to the entropy definition (\ref{ENTROPY_DEFINITION}), Eq.~(\ref{EQ_ENTROPY_1D_DG_APPL}) can be written as:
\begin{equation}
\label{OPERATOR_CONSTRAINT}
p(^{\cL}U_e(x)) \geq \exp(s^0_e)\rho^\gamma (^{\cL}U_e(x))\qquad\forall x\in \mathcal{D}\;.
\end{equation}
To define the operator $\cL$, we follow the work of Zhang and Shu~\cite{ZHANG_SHU_JCP2010,ZHANG_SHU_NUMERMATH2010}, and 
introduce a linear scaling:
\begin{equation}
\label{PROJECTED_SOLUTION}
^\mathcal{L}U_e = U_e + \varepsilon (\overline{U}_e  - U_e)\;.
\end{equation}
The parameter $\varepsilon$ is then determined by substituting Eq.~(\ref{PROJECTED_SOLUTION}) into Eq.~(\ref{OPERATOR_CONSTRAINT}) and by applying Jensen's inequality:
\begin{eqnarray}
\label{EQ_EPSILON_RELATION}
p\left((1- \varepsilon) U_e +  \varepsilon \overline{U}_e\right) & \geq &(1-\varepsilon) p(U_e) + \varepsilon p(\overline{U}_e) \;,\nonumber\\
 & \geq & \exp({s}^0_e) \left[(1-\varepsilon) \rho^\gamma(U_e) + \varepsilon \rho^\gamma(\overline{U}_e) \right]\;,\\
& \geq & \exp({s}^0_e) \rho^\gamma \left((1- \varepsilon) U_e +  \varepsilon \overline{U}_e\right)\;.\nonumber
\end{eqnarray}
Solving for $\varepsilon$ gives
\begin{equation}
\label{EVALUATE_VAREPSILON}
 \varepsilon = \frac{\tau}{\tau - [p(\overline{U}_e) - \exp({s}^0_e)\rho^\gamma(\overline{U}_e)]}\qquad\text{with}\qquad
 \tau = \min\left\{0,~\min_{x\in \mathcal{D}}\{ p(U_e(x)) - \exp({s}^0_e) \rho^\gamma(U_e(x))\}\right\}\;,
\end{equation}
which is subject to the conditions
\begin{equation}
\label{CONSTRAINTS_ON_MEAN}
\rho(\overline{U}_e) > 0\;,\qquad
p(\overline{U}_e) >\exp(s^0_e) \rho^\gamma(\overline{U}_e)\;.
\end{equation}
These conditions are automatically guaranteed through the CFL$^{\rm EB}$-constraint of Lemma \ref{ENTROPY_THREE_ELEMENT_DG}. While the positivity condition for pressure is embedded in Eq.~(\ref{EQ_EPSILON_RELATION}), the positivity of density must be imposed for all $x \in \mathcal{D}$ before $\mathcal{L}$ is applied, and the methodology for this is presented in \cite{ZHANG_SHU_JCP2011}.

Compared to the limiting operator, presented in \cite{ZHANG_SHU_NUMERMATH2010}, the herein proposed method is substantially simplified. Specifically, the step for imposing the positivity of pressure is avoided; in addition, $\varepsilon$ is obtained from an algebraic relation, and does not require a computationally expensive Newton iteration. It is also noted that the operator $\cL$ contains the positivity-preserving limiter as a special case, which is obtained by setting $s_e^0 \to -\infty$.
\subsection{Numerical analysis of the limiting operator $\mathcal{L}$}
In this section, numerical properties of the limiting operator are examined. 
\paragraph{Conservation} 
Integrating Eq. (\ref{PROJECTED_SOLUTION}) over $\Omega_e$, 
\begin{equation}
\nonumber
\int_{\Omega_e} {}^\mathcal{L} U_e dx 
= (1-\varepsilon) \int_{\Omega_e} U_e dx + \varepsilon \overline{U}_e \int_{\Omega_i} dx 
=   \int_{\Omega_e}  U_e dx\;,
\end{equation}
confirms that the limiting operator preserves the conservation properties of the solution vector.
\paragraph{Stability} Since the positivity of density and pressure is preserved at the quadrature points, $\mathcal{L}$ is $L_1$-stable, which was shown in~\cite{ZHANG_SHU_JCP2010, ZHANG_SHU_NUMERMATH2010}. Here, we extend this stability analysis and evaluate the $L_2$-stability. By considering a periodic domain and taking the $L_2$-norm of Eq.~(\ref{PROJECTED_SOLUTION}) we obtain:
\begin{subeqnarray}
||^\mathcal{L}U_e||^2 & = & \int_{\Omega_e} [U_e + \varepsilon ( \overline{U}_e - U_e)]^2 dx\;,\\
& = & (1-\varepsilon)^2 \int_{\Omega_e} U_e^2 dx - \varepsilon (\varepsilon-2)  \int_{\Omega_e} \overline{U}_e^2 dx\;,\\
& \leq & (1-\varepsilon)^2 \int_{\Omega_e} U_e^2 dx - \varepsilon (\varepsilon-2) \int_{\Omega_e} U_e^2 dx\;,\\
& \leq & ||U_e||^2\;.
\end{subeqnarray}
After integrating over the entire domain, we obtain 
\ben
\nonumber
 ||^\mathcal{L}U||^2_{\Omega}  \leq ||U||^2_{\Omega}\;,
\een
which shows that $\mathcal{L}$ does not affect the stability of the DG-discretization. Further, since $\cL$ constrains pressure and density, $\lambda$ in Eq.~(\ref{CFL_CONDITION_1D_DG}) provides a robust CFL-criterion, without the need for arbitrarily reducing $\Delta t$ to increase the stability region. 

\paragraph{Accuracy} In regions where the solution is smooth, we assume that the weak solution before limiting has optimal accuracy:
\ben
||U - \mathsf{U} ||_\Omega \leq C_1 h^{p+1}\;,
\een
and that undershoots in entropy remain small, so that $\varepsilon \sim \cO(h ^p)$. Thus, the error is estimated as follows:
\begin{subeqnarray}
||^\mathcal{L} U - \mathsf{U}||_\Omega^2 & = & \sum_e ||^\mathcal{L}U_e - \mathsf{U}_e||^2\;,\\
& = & \sum_e ||\varepsilon (\overline{U}_e - \mathsf{U}_e) + (1-\varepsilon) (\mathsf{U}_e-U_e) ||^2\;,\\
& \leq & \sum_e \left( 2\varepsilon^2||(\overline{U}_e-\mathsf{U}_e) ||^2 + 2(1-\varepsilon)^2 ||U_e-\mathsf{U}_e ||^2 \right)\;,\\
& \leq & C_2 h^{2p+2}\;,
\end{subeqnarray}
where for simplicity, we introduce $\mathsf{U}_e$ to denote the element-wise representations to $\mathsf{U}$. Here we use the fact that $\overline{U}_e$ is locally a first-order approximation to $\mathsf{U}_e$, $\overline{U}_e = \mathsf{U}_e + C^e_3 \cO(h)$. 

In the vicinity of a discontinuity, the DG-solution looses its regularity so that the convergence rate reduces to first-order: $||U - \mathsf{U} ||_\Omega \leq C_4 h.$

Triggered by spurious sub-cell solutions, the entropy undershoot can be very large, so that $\varepsilon \sim \cO(1) $. By repeating the above argument, we obtain an estimate for the accuracy of the discontinuous solution:
\ben
||^\mathcal{L}U - \mathsf{U} ||_\Omega \leq C_5 h.
\een
The accuracy arguments given here are substantiated through numerical tests in Sec. \ref{SEC_NUMERICAL_TEST}.
 
\section{\label{SEC_GENERALIZATION}Generalization to multi-dimension and arbitrary elements}
The entropy-bounded DG scheme that was presented for one-dimensional systems in the previous section can be generalized to arbitrary elements in multi-dimensions. This extension is the subject of the following analysis.

Since EBDG does not rely on a specific quadrature rule, any quadrature method can be used as long as it accurately integrates the problem and ensures the positivity of the quadrature weights. The limiting procedure requires the definition of a new set of quadrature points $\mathcal{D}$ for the general multi-dimensional setting. The selection of these points is given in the next section. The extension to arbitrary elements requires special consideration of the CFL$^{\rm EB}$ number.

\subsection{\label{SEC_GENERALIZATION_MULT_EXT}Generalization to multi-dimension and arbitrary elements}
To present a general formulation for multi-dimensional configurations, we first introduce necessary notations to describe general elements with curvatures. For this, we define a geometric mapping function $\Phi: \mathbb{R}^{N_d} \rightarrow \mathbb{R}^{N_d}$ on a reference element $\Omega_e^{\rm{r}}$, such that $x = \Phi(r)$ maps any point $r \in \Omega_e^{\rm{r}}$ onto $x \in \Omega_e$, and $\mathcal{J}  = [\partial x/\partial r ]$ is the geometric Jacobian. With these specifications, we can write the discretized state vector as:
\begin{equation}
U_e(x, t) = \sum_{m=1}^{N_p} ~\widetilde{U}_{e,m}(t)\varphi_m(r)\;,\qquad x = x(r) \in \Omega_e\;,\qquad \forall r \in \Omega_e^{\rm{r}}
\end{equation}
The mapping function is commonly parameterized by a polynomial function $x(r) = \sum_{m=1}^{N_g} \widetilde{x}_m\varphi^g_m(r) $, where $\varphi^g_m(r)$ is a Lagrangian interpolation and $N_g$ is the number of geometric bases used to represent $\Omega_e$. Since the reference element is regular, we can use a subspace of $r$ to parameterize the element edges. Therefore, to parameterize the $k$th edge of $\Omega_e$ we define $g_k = \mathcal{P}_k(r) \in \mathbb{R}^{N_d-1}$ such that $\forall r \in \partial \Omega_{e, k}^{\rm{r}}$, $r = \mathcal{P}^{-1}_{k}(g_k)$, in which $\mathcal{P}^{-1}_k$ is the pseudo-inverse of $\mathcal{P}_k$. For the physical element, the edge can be represented as:
\begin{equation}
\partial \Omega_{e,k} = \left\{x \in {\Omega}_e ~|~ x  = \Phi(r),~ r = \mathcal{P}^{-1}_k (g_k) \in \partial \Omega_{e,k}^{\rm{r}}\right\}\;.
\end{equation}

The integral in Eq.~(\ref{WEAK_FORM}) is evaluated using multi-dimensional quadrature rules. Considering the complexity of the dimensionality, here we follow the quadrature convention that is $\sum_{v=1}^{N_q} w_v = V_e^{\rm{r}}$ (the volume of $\Omega_e^{\rm{r}}$) and $\sum_{q=1}^{N_q^k} w_{k,q} = S_{e,k}^{\rm{r}}$ (the area of $\partial \Omega_{e,k}$). With these preliminaries, we can evaluated any volume integral in Eq.~(\ref{WEAK_FORM}) as:
\begin{equation}
\int_{\Omega_e} f(x) dx = \int_{\Omega_e^{\rm{r}}} |\mathcal{J}(r)| f(x(r)) dr = 
                          \sum_{v=1}^{N_q}|\mathcal{J}(r_v)| f(x(r_v))w_v\;.
\end{equation}
The surface integral of a scalar function on $\partial \Omega_{e,k}$ can be written as:
\begin{eqnarray}
\int_{\partial \Omega_{e,k}} f(x) d\Gamma &=& \int_{\partial \Omega_{e,k}^{\rm{r}}} f(x(g_k)) |\mathcal{J}_k^{\partial}| dg_k =  \sum_{q=1}^{N_{q,k}^{\partial}}|\mathcal{J}_k^\partial(g_{k,q})|f(x(g_{k,q})) w_{k,q}\;,
\end{eqnarray}
and the surface integral of a vector-valued function is evaluated as:
\begin{eqnarray}
\int_{\partial \Omega_{e,k}}f(x) \cdot \wh{\textmd{n}} d\Gamma &=& \int_{\partial \Omega_{e,k}^{\rm{r}}} f(x(g_k)) \cdot \mathcal{J}^\partial_k dg_k = \sum_{q=1}^{N^{\partial}_{q,k}} |\mathcal{J}_k^\partial(g_{k,q})|f(x(g_{k,q})) \cdot \wh{\textmd{n}}(g_{k,q})w_{k, q}\;,
\end{eqnarray}
where $\mathcal{J}^\partial_k$ is the surface Jacobian, and $\widehat{\textmd{n}}$ refers to the unit vector ${\mathcal{J}^\partial_k}/{|\mathcal{J}^\partial_k|}$. Note that for a general element, the quadrature expression might be subject to a tiny numerical error bounded by $\cO(h ^{2p+1})$, given the accuracy requirement for integrating Eq.~(\ref{WEAK_FORM}). 


With the above notation, we are now able to define the set of quadrature points $\mathcal{D}$ for general curved elements:
\begin{equation}
\mathcal{D} = \bigcup_{k=1}^{N_\partial} \{g_{k,q},~q=1,\ldots,N^\partial_{q,k}\} \bigcup \{r_v,~v=1,\ldots,N_q\}\;,
\end{equation}
where $N_\partial$ is the number of element edges (which is equal to the number of neighbor elements). In this context, it is noted that $\mathcal{D}$ includes all quadrature points that are involved in the integration. With this specification of $\mathcal{D}$, the limiting operator $\mathcal{L}$, developed in Sec.~\ref{1D_LIMITER_DESIGN}, can be directly extended to arbitrary elements on multi-dimensional configurations. In the following, a CFL$^{\rm EB}$-constraint is derived that extends the results of  Lemma~\ref{ENTROPY_THREE_ELEMENT_DG}, thereby ensuring the existence of the general limiter $\mathcal{L}$.

\subsection{\label{SUBSEC_CFL_CONSTRAINT_ARB_ELEMENT}CFL-constraint}
\begin{figure}[!htb!]
\centering
\includegraphics[width=0.7\textwidth, clip=, keepaspectratio]{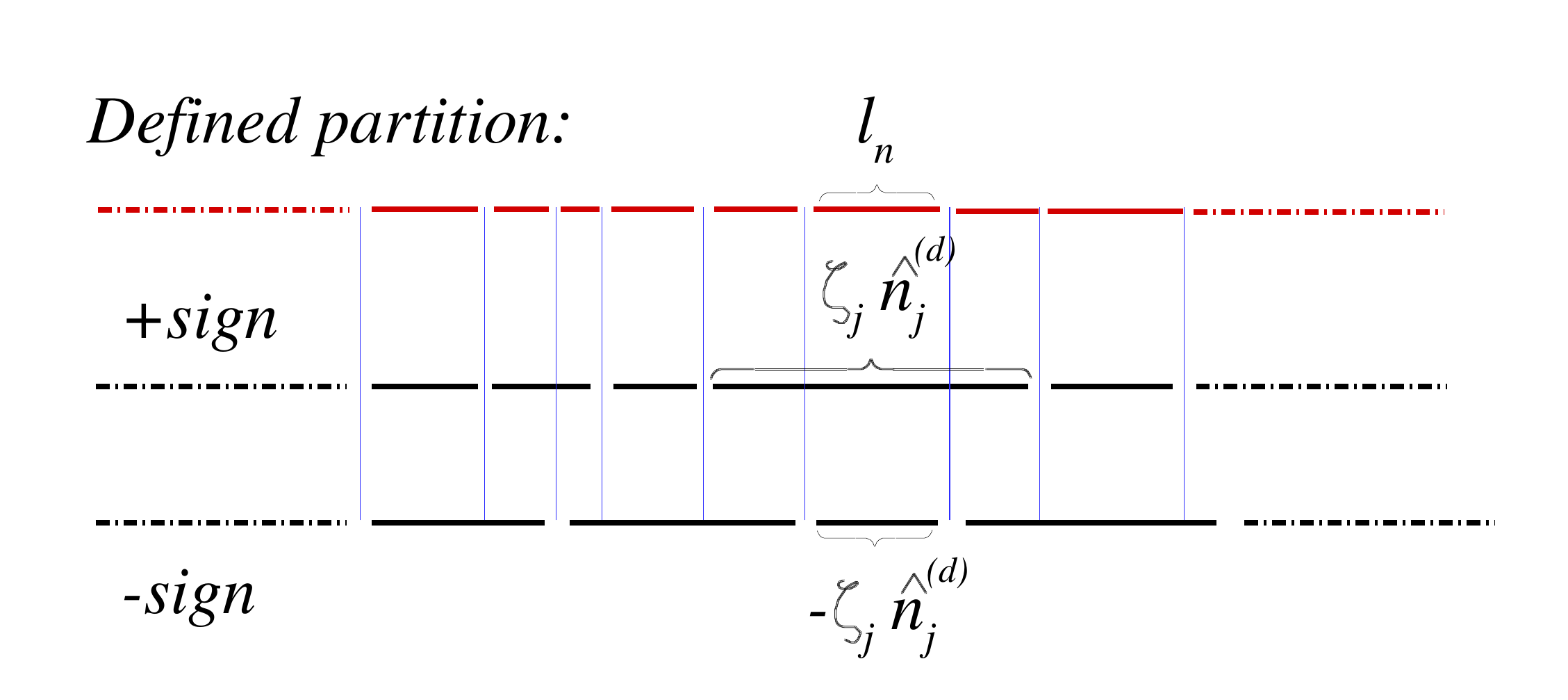}
    \caption{\label{PARTITION_IDEA_DEMO} Illustration of the partition introduced in Lemma~\ref{LEMMA_ENTROPY_MULTID}.}
\end{figure}

Following the same approach as for the one-dimensional derivation in Sec. \ref{DERIVE_CFL_1DDG}, the element-averaged solution of $U_e^{\Delta t}$ is evaluated as:
\begin{subeqnarray}
 \label{DG_MEAN_GENERAL_ELEMENT}
\overline{U}_e^{\Delta t} & = & \overline{U}_e - \frac{\Delta t}{V_e} \sum_{k=1}^{N_\partial}\int_{\partial \Omega_{e,k}} \wh{F}\left(U_e^{+}, U_e^{-}, \wh{\textmd{n}}\right) d\Gamma\;,\\
& = & \sum_{v=1}^{N_q} \frac{|\mathcal{J}(r_{v})| w_v}{V_e}U_e(r_v) - \sum_{k=1}^{N_\partial}\sum_{q=1}^{N_{q,k}^\partial}  \frac{\Delta t |\mathcal{J}^\partial_k(g_{k,q})|w_{k,q}}{V_e} \wh{F}\left(U_e^+(r(g_{k,q})), U_e^-(r(g_{k,q})), \wh{\textmd{n}}(g_{k,q})\right) \;,\\
\nonumber
& = &  \sum_{v=1}^{N_q}\theta_v U_e(r_v)
+ \sum_{k=1}^{N_\partial}\sum_{q=1}^{N^\partial_{q,k}} \left[\theta_{k,q} U_e^+(r(g_{k,q})) \right. \\
 \slabel{DG_MEAN_GENERAL_ELEMENT_C}
&& \left. - {\Delta t \zeta_{k,q}} \left(\wh{F}\left(U_e^+(r(g_{k,q})), U_e^-(r(g_{k,q})), \wh{\textmd{n}}(g_{k,q})\right)  + \wh{F}\left(U_e^+(r(g_{k,q})), U_e^*, -\wh{\textmd{n}}(g_{k,q})\right) \right) \right] \;,
\end{subeqnarray}
where 
\be
\theta_v & = & \frac{|J(r_v)| w_v}{V_e} -  \sum_{k=1}^{N_\partial}\sum_{q=1}^{N^\partial_{q,k}} \theta_{k,q}\phi_v\left(r(g_{k,q})\right)\;
\ee
is introduced to decompose the volumetric quadrature to obtain $U_e^+(r(g_{k,q}))$. For notational simplification, we defined
\be
\zeta_{k,q} & = & \frac{|\mathcal{J}_k^\partial(g_{k,q})|w_{k,q}}{V_e}\;,
\ee
so such that $S_e=\sum_{k=1}^{N_\partial}\sum_{q=1}^{N^\partial_{q,k}}\zeta_{k,q}$ is equal to the surface area of $\Omega_e$, and $\sum_{k=1}^{N_\partial}\sum_{q=1}^{N^\partial_{q,k}}\zeta_{k,q} \wh{\textmd{n}}(g_{k,q}) =0$ since $\Omega_e$ has a closed surface. To apply the results from the three-point system to the multi-dimensional configuration, we introduce the auxiliary variable $U_e^*$:
\be
\label{MULTID_U_STAR}
U_e^* =  \sum_{k=1}^{N_\partial}\sum_{q=1}^{N^\partial_{q,k}}\frac{\zeta_{k,q}}{S_e} \left[ U_e^+(r(g_{k,q})) -\frac{1}{\lambda^*} F(U_e^+(r(g_{k,q}))) \cdot \wh{\textmd{n}}(g_{k,q})\right]\;.
\ee
It can be shown that $U_e^*$ is essentially the solution to the following equation:
\be
  \sum_{k=1}^{N_\partial}\sum_{q=1}^{N^\partial_{q,k}} \zeta_{k,q} 
  \wh{F}\left(U_e^+(r(g_{k,q})), U_e^*, -\wh{\textmd{n}}(g_{k,q})\right) = 0\;,
\ee
subject to a preselected dissipation coefficient $\lambda^*$, so that the equality in Eq.~(\ref{DG_MEAN_GENERAL_ELEMENT_C}) holds true. Here, we evaluate $\lambda^*$ from the following relation: 
\be
  \label{LAMBDA_STAR}
  \lambda^* = \tau \max\left\{\nu(U)~|~U \in \{U_e^+(r(g_{k,q}))\}\right\},\qquad \tau  = \max\left\{\sqrt{N_d}, \sqrt{2+\gamma(\gamma-1)}\right\}
\ee
and the rationale for this selection is provided later in Remark~\ref{LAMBDA_RATIONALE}. To prove that $\overline{U}_e^{\Delta t}$ is entropy bounded, we present the following lemma. 
\begin{lemma}
\label{LEMMA_ENTROPY_MULTID}
$U_e^*$ in Eq. (\ref{MULTID_U_STAR}) satisfies $s(U_e^*) \geq \min \left\{s(U)~|~ U \in \{U_e^+(r(g_{k,q}))\}\right\}$.
\\
Proof: For notational simplification, we combine the indices $k$ and $q$ into a single index $j$, and we denote the total number of surface quadrature points on $\p\Omega_e$ by $N_{\rm{tot}}$, $N_{\rm{tot}} = \sum_{k=1}^{N_\partial} N^\partial_{q,k}$. Considering $\sum_{j=1}^{N_{\rm{tot}}} \zeta_j\wh{\textmd{n}}_j^{(d)} = 0$ and $\zeta_j >0$, the $d$th components of the surface-normal vectors have different signs. To denote each component, we introduce the superscript $(d)$. By sorting $\wh{\textmd{n}}_j^{(d)}$ so that the first $N^>_{\rm{tot}}$ vector components are positive. The following statement is true for any $d$:
\be
\sum\limits_{j=1}^{N_{\rm{tot}}^>} \zeta_j\wh{\textmd{n}}_j^{(d)} = -\sum\limits_{j=N_{\rm{tot}}^> + 1}^{N_{\rm{tot}}} \zeta_j\wh{\textmd{n}}_j^{(d)} = \sum_{n=1}^{N_{par}} l_n\;,
\ee
where $l_n$ introduces a partition as illustrated in Fig. \ref{PARTITION_IDEA_DEMO} and $N_{par}$ is the dimension of this partition. With this, we are able introduce a variable mapping,
\ben
U_{n}^{s+} = U_e^+(r_j),~~&\text{if}&~\sum\limits_{i=1}^{j-1} \zeta_i\wh{\textmd{n}}_i^{(d)}< \sum\limits_{i=1}^n l_i \leq \sum\limits_{i=1}^{j}\zeta_i\wh{\textmd{n}}_i^{(d)} \;,\\
U_{n}^{s-} = U_e^+(r_j),~~&\text{if}&~-\sum\limits_{i=N_{\rm{tot}}^> + 1}^{j-1} \zeta_i\wh{\textmd{n}}_i^{(d)}< \sum\limits_{i=1}^n l_i \leq -\sum\limits_{i=N_{\rm{tot}}^> + 1}^{j}\zeta_i\wh{\textmd{n}}_i^{(d)} \;.
\een

With this, Eq. (\ref{MULTID_U_STAR}) is equivalent to:
\ben
U_e^* &=& \frac{1}{S_e}\left(\sum_{j=1}^{N_{\rm{tot}}} \zeta_j U_e^+(r_j) - \sum_{j=1}^{N_{\rm{tot}}} \frac{\zeta_j}{\lambda^*}F(U_e^+(r_j)) \cdot \wh{\textmd{n}}_j \right)\;,\\
 & = & \sum_{j=1}^{N_{\rm{tot}}} \frac{\zeta_j}{S_e} \left(1 -  \sum_{d=1}^{N_d} \frac{|\wh{\textmd{n}}_{j}^{(d)}|}{\sqrt{N_d}}\right) U_e^+(r_j) + \sum_{d=1}^{N_d}  \sum_{j=1}^{N_{\rm{tot}}} \frac{1}{S_e}\left(\frac{\zeta_j|\wh{\textmd{n}}_{j}^{(d)}|}{\sqrt{N_d}} U_e^+(r_j)- \frac{\zeta_j\wh{\textmd{n}}_j^{(d)}}{\lambda^*} F^{(d)}(U_e^+(r_j))\right) \;,\\
 & = & \sum_{j=1}^{N_{\rm{tot}}} \frac{\zeta_j}{S_e} \left(1 -  \sum_{d=1}^{N_d} \frac{|\wh{\textmd{n}}_{j}^{(d)}|}{\sqrt{N_d}}\right) U_e^+(r_j) + \sum_{d=1}^{N_d} \frac{1}{S_e}\left(\frac{2}{\sqrt{N_d}} \sum_{n=1}^{N_{par}}l_n U_{d,n}^{**}\right) \;,\\
\een
where we introduce
\ben
 U_{d,n}^{**} &=& \frac{1}{2}(U_{n}^{s+} + U_{n}^{s-}) - \frac{\sqrt{N_d}}{2\lambda^*}\left(F^{(d)}(U_{n}^{s+}) - F^{(d)}(U_{n}^{s-})\right)\;,
\een
which takes the same form as the left-hand-side of Eq.~(\ref{EXP_U_STAR}). Note that $U_{d,n}^{**}$ is essentially expressed in a one-dimensional setting along $x^{(d)}$. Therefore, one can follow the same argument used for Eq. (\ref{EXP_U_STAR}) in Lemma \ref{ENTROPY_THREE_ELEMENT_DG} to verify that 
\ben
s(U_{d,n}^{**} ) \geq \min \left\{s(U_{n}^{s\pm})\right\} \geq \min \left\{s(U)~|~ U \in \{U_e^+(r_j)\},~j = 1, \ldots, N_{\rm{tot}}\right\}\;,
\een
with the given form of $\lambda^*$ in Eq. (\ref{LAMBDA_STAR}). As given above, $U_e^*$ is a convex combination; by using the quasi-concavity of entropy~\cite{ZHANG_SHU_NUMERMATH2010}, we conclude that
\ben
s(U_e^*) > \min \left\{s(U)~|~ U \in \{U_e^+(r_j)\},~j = 1, \ldots, N_{\rm{tot}}\right\}\;.
\een
\end{lemma}

\begin{remark}
\label{LAMBDA_RATIONALE}
Note that the maximum characteristic speed of $U_{d,n}^{**}$ is bounded, $\nu(U_{d,n}^{**}) \leq \nu(U_e^+(r))$. According to the combination law given in Appendix~\ref{APP_COMB_RULE}, we have $\nu(U_e^*)  \leq \sqrt{2+ \gamma(\gamma-1)} \max \left\{\nu(U_e^+(r))\right\} \leq \lambda^*$. With this, we are now able to show that the maximum characteristic speed of $U_e^*$ is bounded by the chosen value of $\lambda^*$, so that the Lax-Friedrichs flux $\wh{F}\left(U_e^+(r), U_e^*, -\wh{\textmd{n}}\right)$ in Eq. (\ref{DG_MEAN_GENERAL_ELEMENT}) is valid according to the definition of Eq. (\ref{LAX_FRIEDRICHS_FLUX}). 
\end{remark}
To enforce the entropy boundedness, the decomposition of $\overline{U}_e^{\Delta t}$ in Eq. (\ref{DG_MEAN_GENERAL_ELEMENT}) is required to be convex. This can be satisfied under the following condition:
\be
\label{DG_ENTROPY_CONVEX_CONS}
\begin{cases}
\theta_v \geq 0\;,~~\forall v = 1,\ldots,N_q,\\
\theta_{k,q} > 0\;,~~\forall (k,~ q),~q = 1,\ldots, N^\partial_k,~k = 1,\ldots, N_\partial\;,
\end{cases}
\ee
With this, the entropy boundedness of $\ol{U}_e^{\Delta t}$ is shown by the following lemma.
\begin{lemma}
\label{ENTROPY_GENERAL_ELEMENT_DG}
For a general DG element, the element-averaged solution is entropy bounded,
\begin{equation}
s(\overline{U}^{\Delta t}_e) \geq s_e^0(t) = \min\{s\left(U(y)\right)|~ y \in \Omega_e \cup \partial \Omega_e^-\},
\end{equation}
under the condition that Eq. (\ref{DG_ENTROPY_CONVEX_CONS}) holds and that the following constraint is fulfilled:
\begin{equation}
\label{CFL_GENERAL}
\Delta t \lambda \leq \frac{1}{2}\min~\left\{\frac{\theta_{k,q}}{\zeta_{k,q}}  \right\},~~\forall (k,~ q),~q = 1,\cdots, N^\partial_k,~k = 1,\cdots, N_\partial\;,
\end{equation}
where $\lambda \geq \max \{\nu(U)~|~ U \in \{U^{\pm}_e(r(g_{k,q}))\}$  and $\lambda \geq \lambda^*$.
\\
Proof: The proof follows Lemma~\ref{ENTROPY_THREE_ELEMENT_DG}, utilizing Lemma~\ref{LEMMA_THREEPOINT_ROT} and the quasi-concavity of entropy. 
\end{lemma}

Note that Lemma~\ref{ENTROPY_GENERAL_ELEMENT_DG} does not rely on any assumption regarding the dimensionality or shape of the finite element, and is therefore general. Another observation is that Eq.~(\ref{CFL_GENERAL}) essentially provides an estimate for CFL$^{\rm EB}$ that is only a function of the geometry of the element. For practical applications, we require the right-hand-side of Eq.~(\ref{CFL_GENERAL}) to be as large as possible so that larger time steps can be taken. This can be achieved by solving a convex optimization problem:
\begin{align}
\label{OPT_FOR_TIME_STEP}
&\textrm{maximize}~ \left(\min~\left\{\frac{\theta_{k,q}}{\zeta_{k,q}}  \right\} \right)\;,\\
\nonumber
&\textrm{subject~to}~ \textrm{Eq. (\ref{DG_ENTROPY_CONVEX_CONS})}
\end{align}
where $\theta_{k,q},~\zeta_{k,q}$ are properties of the geometry alone. This problem can be solved for each individual element as a pre-processing step prior to the simulation. Another way to interpret the expression is to identity a length scale from the right-hand-side of Eq. (\ref{CFL_GENERAL}), for which the CFL$^{\rm EB}$ number can be explicitly defined. For this, $L_e = \min V_e / |\mathcal{J}^\partial_k(g_{k,q})|$ is used as a characteristic length for $\Omega_e$. Hence, 
\ben
 \min~\left\{\frac{\theta_{k,q}}{\zeta_{k,q}}  \right\} \geq L_e \min ~\left\{\frac{\theta_{k,q}}{w_{k,q}} \right\}
\een
and an alternative expression to Eq.~(\ref{OPT_FOR_TIME_STEP}) is
\begin{align}
&\textrm{maximize}~ \min~\left\{\frac{\theta_{k,q}}{w_{k,q}}  \right\}\;,\\
\nonumber
&\textrm{subject~to}~ \textrm{Eq. (\ref{DG_ENTROPY_CONVEX_CONS})}\;,
\end{align}
where the optimal solution is the value of $\text{CFL}^{\text{EB}}$. With this, the CFL-constraint can be written as
\be
\label{CFL_COND_FINAL}
\frac{\Delta t \lambda}{L_e} \leq \frac{1}{2}\text{CFL}^{\text{EB}}\;,
\ee
which is used in the following numerical tests. The factor of 1/2 is a consequence of the Riemann flux formulation. For some of the most relevant element types with regular shapes, the value of $\text{CFL}^{\rm EB}$ has been calculated and listed in Table~\ref{OPT_CFL_NUM_LIST} for different polynomial orders. In practice, we found that the bound in Eq.~(\ref{CFL_COND_FINAL}) leads to a conservative estimate for the time step. Considering the computation of efficiency and the constraint for the linear stability by \cite{COCKBURN_SHU_JSC01}, this condition is relaxed and we consider $0.8\text{CFL}^{\rm EB}$ for the following numerical experiments.

\begin{table}[!htb!]
\centering
\caption{\label{OPT_CFL_NUM_LIST}Summary of quadrature orders and optimal CFL numbers for different types of elements. Quadrature rule (QR) applied: Line, Quadrilateral and Brick: tensor-product Gauss-Legendre; Triangle: Dunavant \cite{DUNAVANT_QUADRATURE_1985}; Tetrahedron: Zhang, et al.~\cite{ZHANG_QUADRATURE_2009}.  (note that Dunavant's triangle rule includes negative weights for 3rd-and 7th-order quadrature, therefore, only quadrature rules with positive weights are used with one extra order.)}
    \begin{tabular}{| c || c | c | c | c | c | c|}
    \hline
Element & Order & QR on $\partial \Omega_e$  & QR on $\Omega_e$ & $\text{CFL}^{\rm EB}$ \\ \hline\hline
 \multirow{4}{*}{\includegraphics[width=0.1\textwidth, clip=, keepaspectratio]{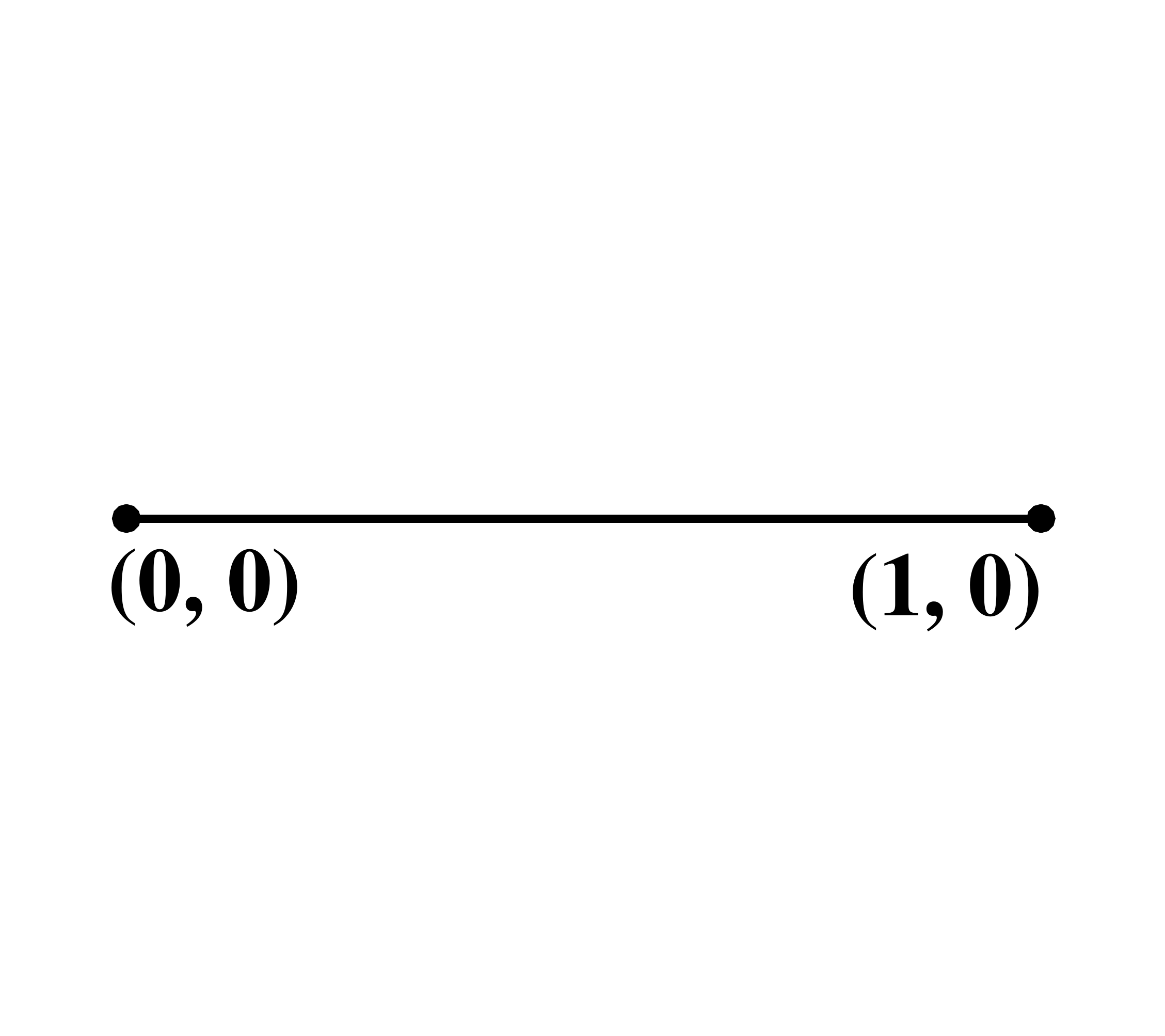}} &  $p=1$       & /  &  3 & 0.5 \\ \cline{2-5}
& $p=2$                & /  &  5 & 0.167 \\ \cline{2-5}
& $p=3$                & /  &  7 & 0.123 \\ \cline{2-5}
& $p=4$                & /  &  9 & 0.073 \\ \cline{2-5} \hline
 \multirow{4}{*}{\includegraphics[width=0.1\textwidth, clip=, keepaspectratio]{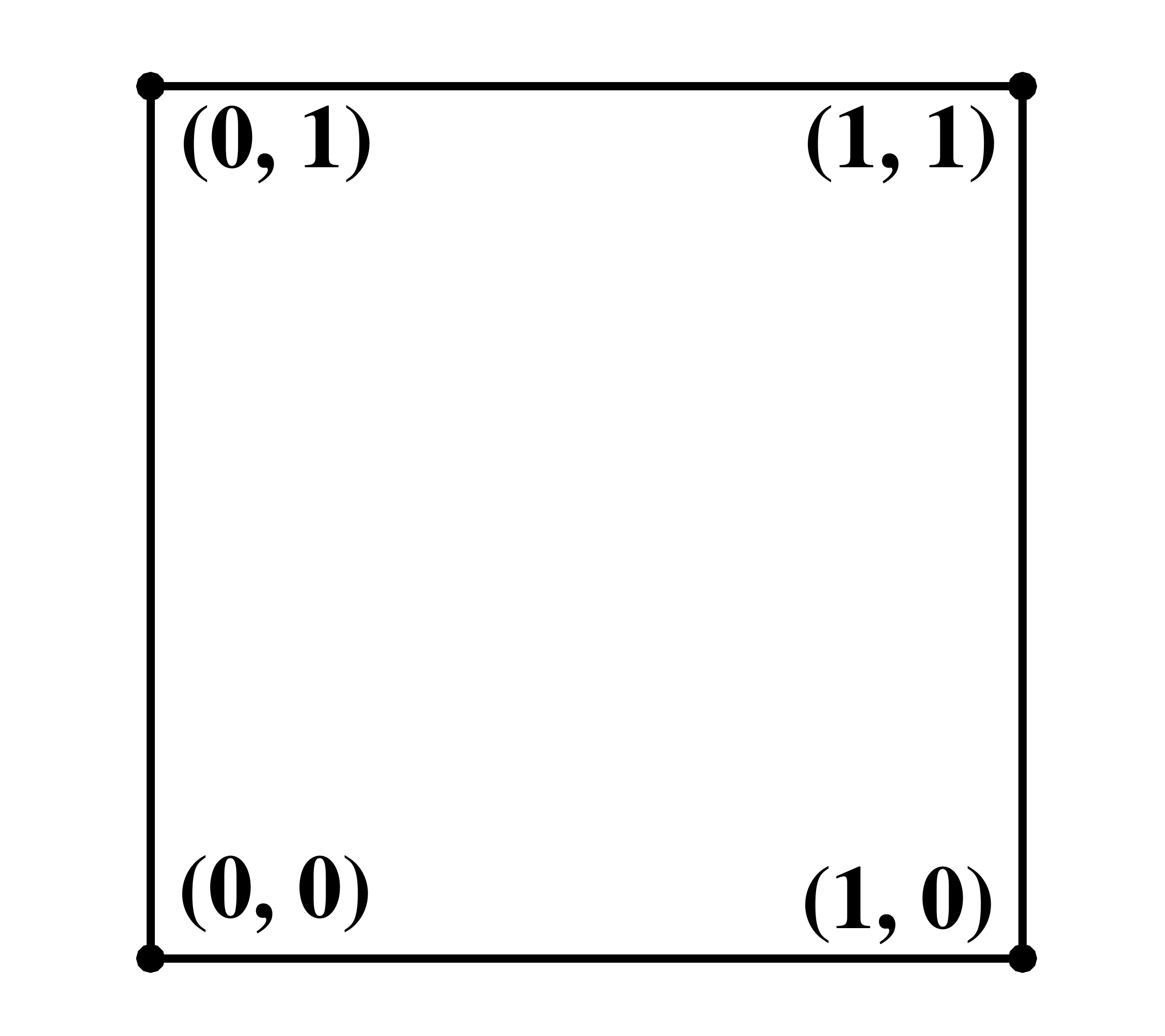}} & $p=1$   & 3 &  3 & 0.25 \\ \cline{2-5}
 &$p=2$                & 5 &  5 & 0.083 \\ \cline{2-5}
 &$p=3$                & 7 &  7 & 0.062 \\ \cline{2-5}
 &$p=4$                & 9 &  9 & 0.036 \\ \cline{2-5} \hline
 \multirow{4}{*}{\includegraphics[width=0.1\textwidth, clip=, keepaspectratio]{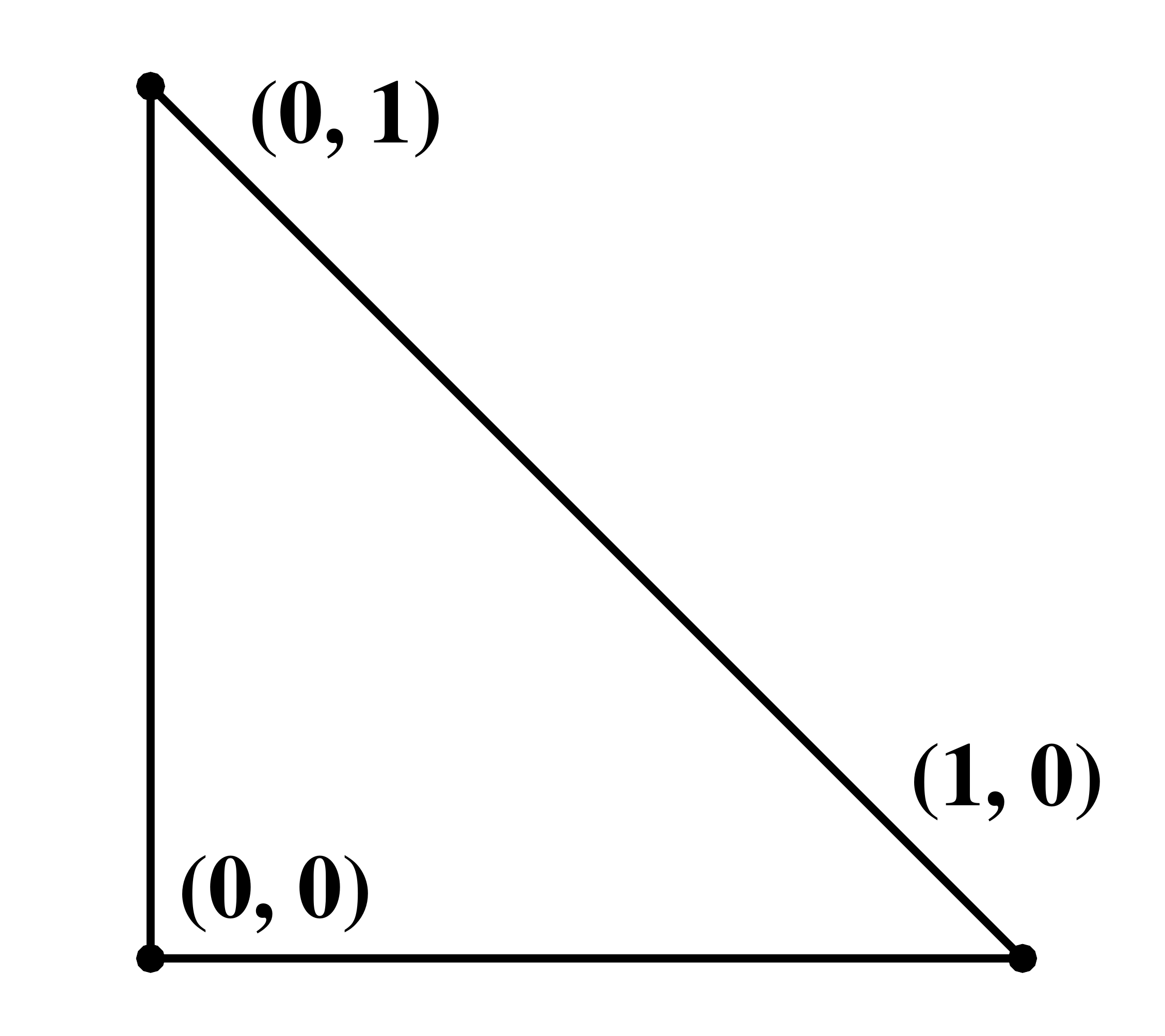}} & $p=1$   & 3 &  4 & 0.135 \\ \cline{2-5}
 &$p=2$                & 5 &  5 & 0.067 \\ \cline{2-5}
 &$p=3$                & 7 &  8 & 0.058 \\ \cline{2-5}
 &$p=4$                & 9 &  9 & 0.033 \\ \cline{2-5} \hline
  \multirow{4}{*}{\includegraphics[width=0.1\textwidth, clip=, keepaspectratio]{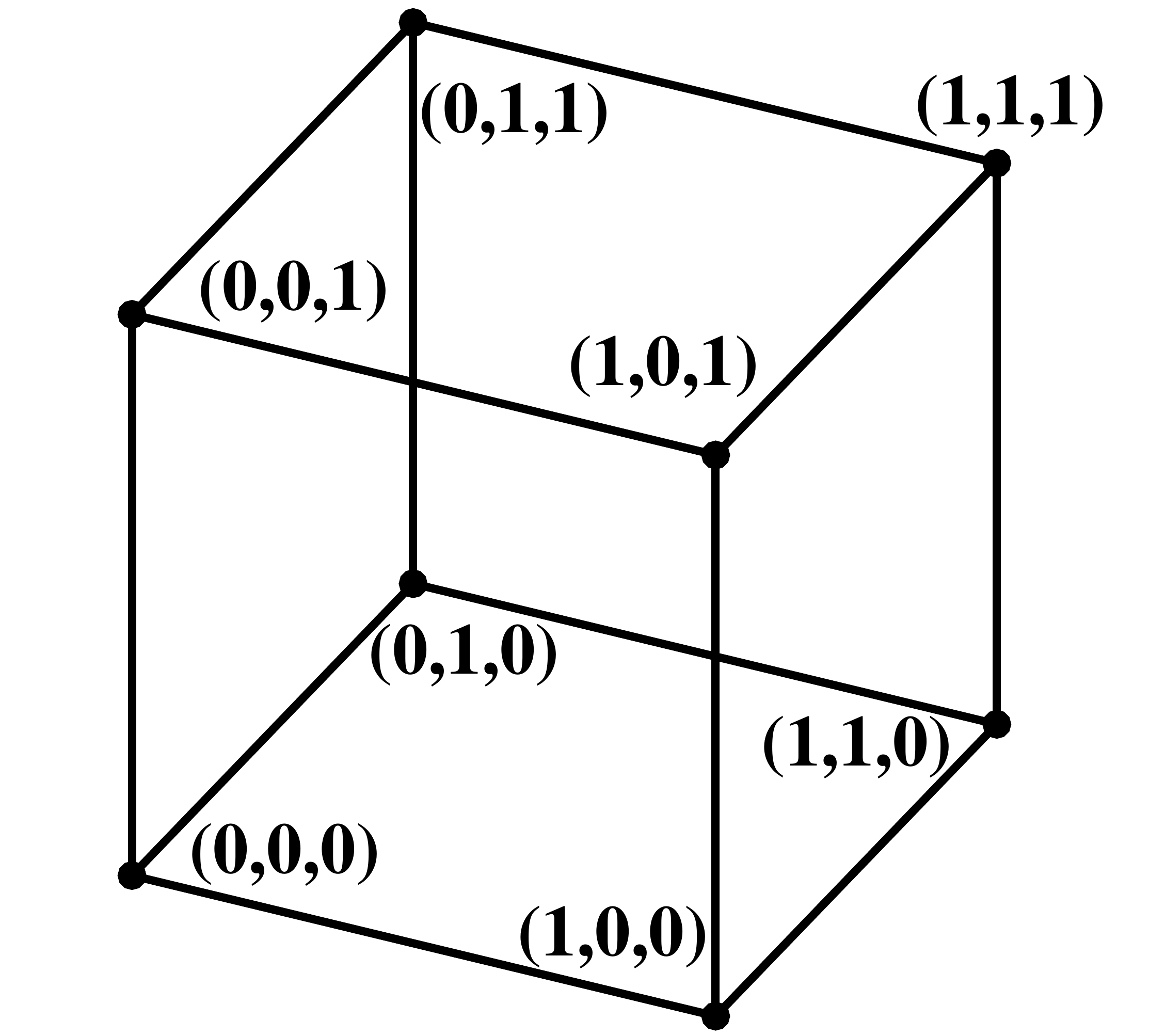}} & $p=1$      & 3 &  3 & 0.167 \\ \cline{2-5}
 &$p=2$                & 5 &  5 & 0.056 \\ \cline{2-5}
 &$p=3$                & 7 &  7 & 0.041 \\ \cline{2-5}
 &$p=4$                & 9 &  9 & 0.024 \\ \cline{2-5} \hline
 \multirow{4}{*}{\includegraphics[width=0.1\textwidth, clip=, keepaspectratio]{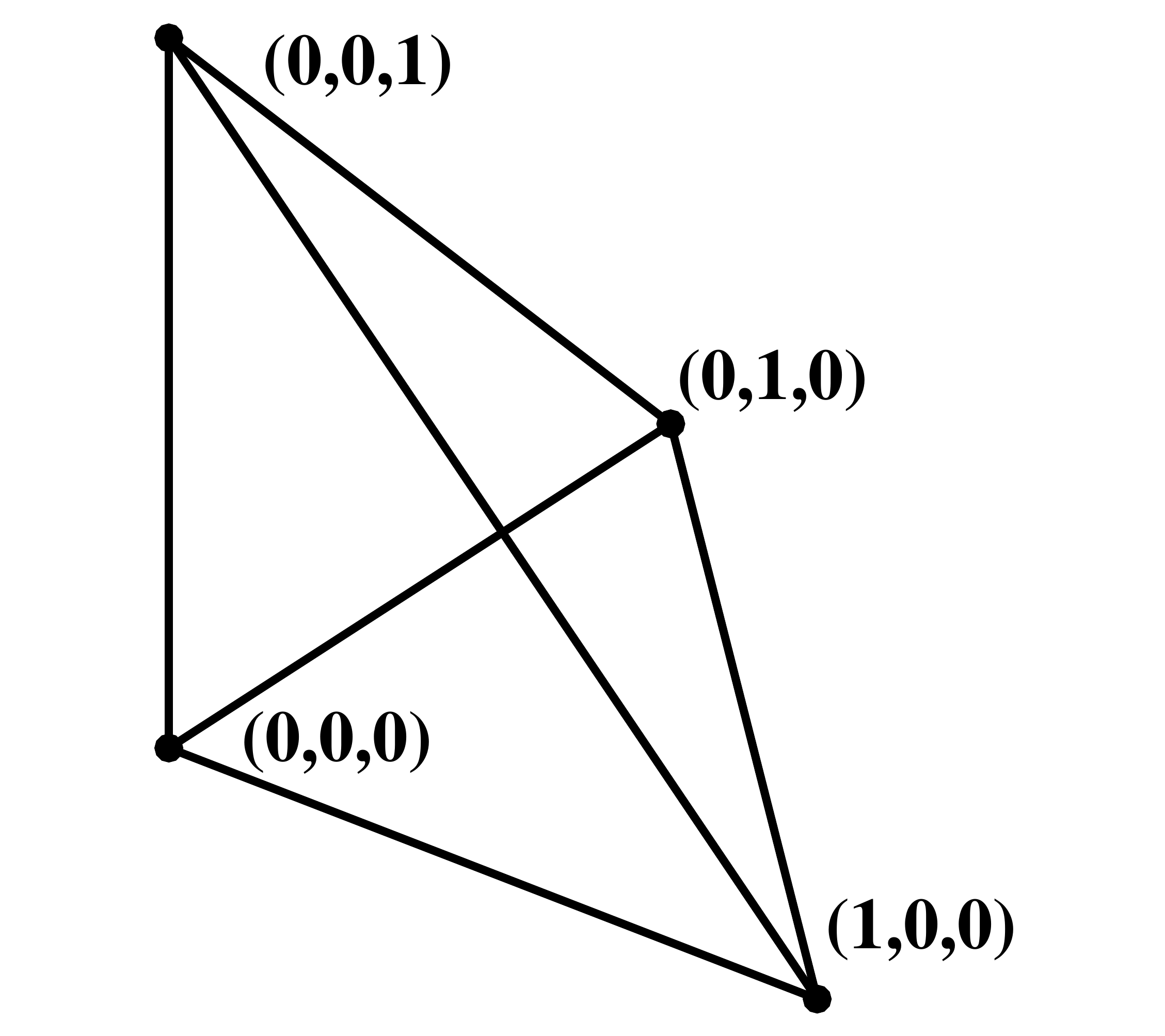}} & $p=1$  & 4 &  3 & 0.066 \\ \cline{2-5}
 &$p=2$                & 5 &  5 & 0.035 \\ \cline{2-5}
 &$p=3$                & 8 &  7 & 0.015 \\ \cline{2-5}
 &$p=4$                & 9 &  9 & 0.013 \\ \hline
    \end{tabular}
\end{table}

\section{\label{SEC_EVALUATION_ENTROPY_BOUND}Evaluation of entropy bound}
\label{Evaluation of Entropy Bound}
In this section, we propose an approach for evaluating the entropy bound $s_e^0(t)$ to answer the fourth implementation problem listed in Sec. \ref{SEC_ENTROPY_BOUNDED_DG}. Obviously, the most accurate way for evaluating a lower bound of entropy is to use Newton's method. However, this approach can significantly impair the efficiency, since searching the minimum on a multi-dimensional high-order element is intractable in terms of computational cost. To overcome this issue, we propose the following two approaches:

\begin{itemize}
\item{\emph{User-defined global bound}}. The first strategy is to let the user specify a global entropy bound, which is then kept constant and used everywhere in the computational domain. Although this approach is simple and robust, it is not optimal.  It is suitable for certain problems with a well-defined entropy bound. As example, for a supersonic flow over an airfoil, the free stream entropy can be used to impose this bound. However, for more complex cases with multiple discontinuities that include several entropy jumps, such a constant bound is not able to enforce the constraint for all elements. Note that this approach recovers the positivity constraint of Zhang and Shu \cite{ZHANG_SHU_JCP2010, ZHANG_SHU_JCP2011} in the limit of $s_e^0(t) \rightarrow -\infty$.


\item{\emph{Estimate of local entropy bound}}. This strategy imposes an entropy bound for each element and dynamically updates $s_e^0(t)$ during the simulation. Instead of relying on a sophisticated search algorithm, $s_e^0(t)$ can be approximately evaluated by reusing available information on quadrature points. According to the definition of $s_e^0(t)$, Eq. (\ref{ENTROPY_PRINCIPLE_FOR_DG}), we also need to consider the set of quadrature points on $\partial \Omega_e^-$, denoted as $\mathcal{D}^-$. Therefore, the estimation for $s_e^0(t)$ is obtained according to the following formulation, 
\begin{equation}
\label{ESITMATE_SENOB_T_F1}
^\mathcal{E}{s}_e^0(t) =  \min\left\{\min_{x\in \mathcal{D}^- }{s(U(x))},~s_{m} - \frac{\min\limits_{x \in \mathcal{D}, x \neq x_m} \{|x_m- x|\}}{|x_m - x_n|} \big(s_{n}-s_{m}\big)\right\}
\end{equation}
where we introduce $x_m$ and $x_n$ to denote the locations of the minimum and maximum entropy values, respectively, 
\ben
s_{m}  = s(U_e(x_m)) = \min\limits_{x \in \mathcal{D}}s(U_e(x))\;,\\
s_{n} = s(U_e(x_n)) = \max\limits_{x \in \mathcal{D}}s(U_e(x))\;.
\een
Although this estimate is simple and inexpensive, one has to realize that any extrapolation in the vicinity of discontinuities becomes dangerous due to the spurious behavior of the sub-cell solution. However, it can be resolved by referring to the entropy bounds around $\Omega_e$ at the last time step,
\begin{equation}
\label{ESITMATE_SENOB_T_F2}
s_e^0(t) = \max \left\{^\mathcal{E}{s}_e^0(t),~ \min_{k \in \mathcal{N}_e \cup \{e\}} s_k^0(t-\Delta t)\right\}\;,
\end{equation}
where $\mathcal{N}_e$ refers to the set of the indices of all neighbor elements of $\Omega_e$ that share a common edge.
\end{itemize}

For practical tests, we found that the above strategy can be applied in a combined way. Specifically, Eq. (\ref{ESITMATE_SENOB_T_F1}) is used for initializing the simulation, and Eq. (\ref{ESITMATE_SENOB_T_F2}) is then applied during the subsequent simulation. 

\section{\label{SEC_ALG_IMPLEMENTATION}Algorithmic implementation}
Algorithm \ref{BRILLIANT_YU} provides a description of the implementation details of the EBDG scheme.

\begin{algorithm}[H]
\label{BRILLIANT_YU}
\SetAlgoLined
 \textbf{Pre-computation of CFL condition}:
 For each element, solve Eq. (\ref{OPT_FOR_TIME_STEP}); alternatively, take CFL$^{\text{EB}}$ from Table~\ref{OPT_CFL_NUM_LIST} and compute $L_e$(recommended for simplicity)

 \textbf{Initialization}: Initialize solution vector $U(x, 0) = U_0$

 \While{$t \leq t_{\rm{end}}$}{
 \For {each element} {
 Estimate entropy bound $s_e^{0}(t)$ according to Eqs. (\ref{ESITMATE_SENOB_T_F1}) and (\ref{ESITMATE_SENOB_T_F2}) \\
  Find $\lambda$ and estimate time step size $\Delta t$ according to Eq. (\ref{CFL_GENERAL})\\
  }
  Find the minimum permissible time step $\Delta t_{\min}$ over all elements\;
  \For {each stage $k$ of a Runge-Kutta integration scheme} {
  \For {each element} {
  step 1: Update solution vector $U^{k+1} = U^k + \Delta t_{\min} R^k$ ($R$ refers to the residual)\\
  step 2: Apply $\mathcal{L}$ on $U^{k+1}$ with $s_e^0(t)$ according to Eqs. (\ref{PROJECTED_SOLUTION}) and (\ref{EVALUATE_VAREPSILON})\\
  }
  }
  Advance time $t = t + \Delta t_{\min}$ 
 }

 \caption{Implementation of EBDG scheme.}
\end{algorithm}

\section{\label{SEC_NUMERICAL_TEST}Results and numerical test cases}
In the following, EBDG is applied to a series of test cases to demonstrate the performance of this method. We begin by considering one-dimensional configurations to confirm the high-order accuracy and essential convergence properties. This is followed by two- and three-dimensional cases  with specific emphasis on applications to unstructured meshes and general curved elements.

\subsection{One-dimensional smooth solution}
The first case considers a one-dimensional periodic domain $x\in[0, 1]$ with smooth initial conditions:
\ben
\rho(x,0) & = & 1 + 0.1\sin(2\pi x)\;,\\
u(x,0) & = & 1\;,\\
p(x,0) & = & 1\;.
\een
The accuracy is examined by considering different spatial resolutions and polynomial orders. For each polynomial order, the CFL number is assigned to 0.8CFL$^{\rm EB}$, in which CFL$^{\rm EB}$ is taken from Table \ref{OPT_CFL_NUM_LIST}. Initially, $s_0$ is set to 0.874, corresponding to the minimum entropy value of the initial condition. The SSPRK33 time-integration scheme~\cite{SSP_REF_2001} is used, and the convergence rate is given in Table~\ref{1D_ACCURATE_SSPRK34}. Although the EBDG-scheme remains stable, it can be seen that the solutions do not reach the optimal rates for DGP3 and DGP4. The reason for this is that the stability region is not sufficiently large for both cases. To demonstrate this, we switch the time-integration scheme to a standard RK45. As can be seen from Table~\ref{1D_ACCURATE_RK45}, the optimal convergence rates for all cases are achieved, demonstrating that the optimal convergence for smooth solutions is preserved by the EBDG-scheme. In the following, the standard RK45 is used for all other cases.
\begin{table}[htp]
	\begin{center}
		\begin{tabular}{|c||cc|cc|cc|cc|}\hline
		 \multirow {2}{*}{$h$} & \multicolumn{2}{c|}{DGP1} & \multicolumn{2}{c|}{DGP2} & \multicolumn{2}{c|}{DGP3} & \multicolumn{2}{c|}{DGP4}\\
		  &  $L_2$-error   & rate & $L_2$-error   & rate  & $L_2$-error   & rate & $L_2$-error   & rate \\
		\hline\hline
		1/10 &				3.074\e{-3} & {-} & 1.274\e{-4} & {-} & 4.716\e{-6} & {-} & 2.036\e{-7} & {-}\\
		1/20 &				6.508\e{-4} & 2.240  & 1.513\e{-5} & 3.073  & 3.073\e{-7} & 3.940  & 1.980\e{-8} & 3.362 \\
		1/40 &				1.535\e{-4} & 2.084  & 1.891\e{-6} & 3.000  & 2.182\e{-8} & 3.816  & 2.454\e{-9} & 3.013 \\
		1/80 &				3.775\e{-5} & 2.024  & 2.364\e{-7} & 3.000  & 1.880\e{-9} & 3.537  & 3.130\e{-10} & 2.971 \\
		1/160 &				9.398\e{-6} & 2.006  & 2.955\e{-8} & 3.000  & 2.001\e{-10} & 3.232  & 3.924\e{-11} & 2.995 \\
		1/320 &				2.347\e{-6} & 2.002  & 3.694\e{-9} & 3.000  & 2.401\e{-11} & 3.059  & 4.922\e{-12} & 2.995 \\
		\hline 
		\end{tabular}
	\end{center}
	\caption{\label{1D_ACCURATE_SSPRK34}Convergence test of 1D advection with SSPRK33, showing degradation of convergence order for DGP3 and DGP4 (here we use density to evaluate the error).}
\end{table}
\begin{table}[htp]
	\begin{center}
		\begin{tabular}{|c||cc|cc|cc|cc|}
		\hline 
		 \multirow {2}{*}{$h$} & \multicolumn{2}{c|}{DGP1} & \multicolumn{2}{c|}{DGP2} & \multicolumn{2}{c|}{DGP3} & \multicolumn{2}{c|}{DGP4}\\
		  & $L_2$-error   & rate & $L_2$-error   & rate &  $L_2$-error   & rate &  $L_2$-error   & rate \\		\hline\hline
		1/10 & 				3.494\e{-3} & {-} & 2.140\e{-4} & {-} & 4.650\e{-6} & {-} & 1.438\e{-7} & {-}\\
		1/20 &				7.231\e{-4} & 2.273  & 1.513\e{-5} & 3.823  & 2.920\e{-7} & 3.993  & 4.517\e{-9} & 4.992 \\
		1/40 & 				1.630\e{-4} & 2.150  & 1.891\e{-6} & 3.000  & 1.826\e{-8} & 3.999  & 1.419\e{-10} & 4.992 \\
		1/80 &				3.790\e{-5} & 2.105  & 2.364\e{-7} & 3.000  & 1.141\e{-9} & 4.000  & 4.444\e{-12} & 4.997 \\
		1/160 &				9.398\e{-6} & 2.012  & 2.955\e{-8} & 3.000  & 7.134\e{-11} & 4.000  & 1.497\e{-13} & 4.892 \\
		1/320 &				2.347\e{-6} & 2.002  & 3.694\e{-9} & 3.000  & 4.463\e{-12} & 3.999  & 8.930\e{-14} & 7.453\e{-1}\\
		\hline 
		\end{tabular}
	\end{center}
	\caption{\label{1D_ACCURATE_RK45}Convergence test of 1D advection with standard RK45 (here we use density to evaluate the error).}
\end{table}

\subsection{\label{1D_MOV_SHOCK}One-dimensional moving shock wave}
\begin{figure}[!tb!]
     \begin{subfigmatrix}{3}
        \subfigure[${\rm{Ma}}=2, h = 1/100$.]{\includegraphics[width=0.32\textwidth, clip=, keepaspectratio]{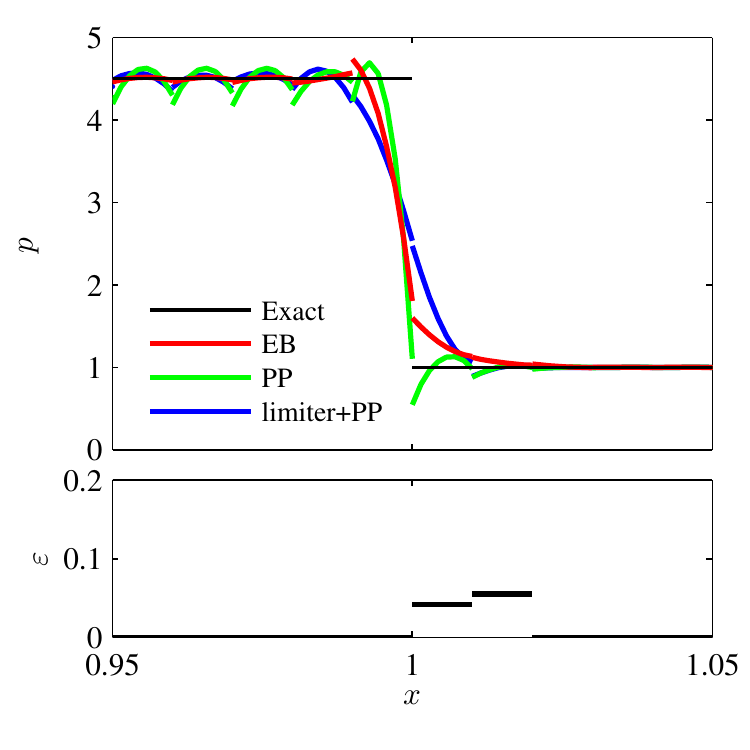}}
        \subfigure[${\rm{Ma}}=5, h = 1/100$.]{\includegraphics[width=0.32\textwidth, clip=, keepaspectratio]{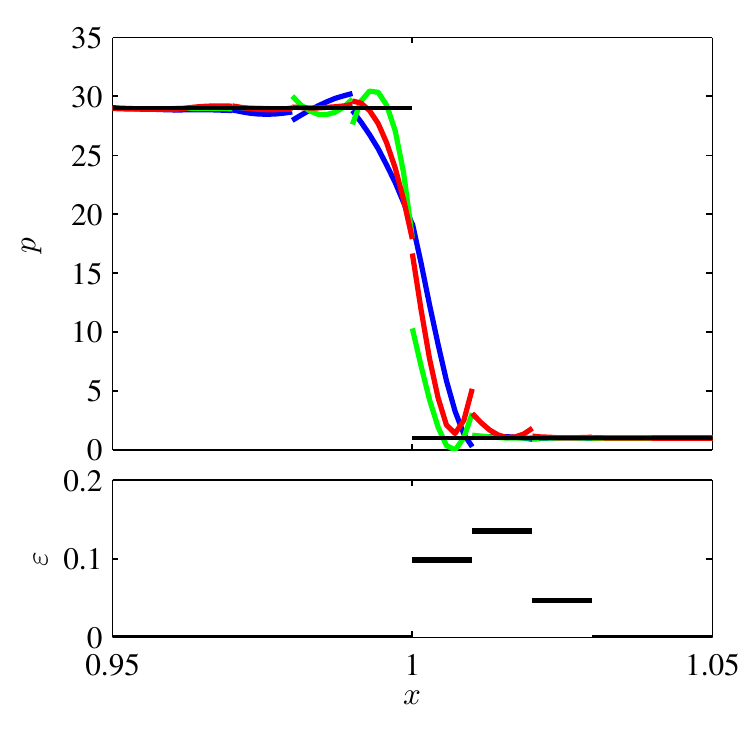}}
        \subfigure[${\rm{Ma}}=100, h = 1/100$.]{\includegraphics[width=0.32\textwidth, clip=, keepaspectratio]{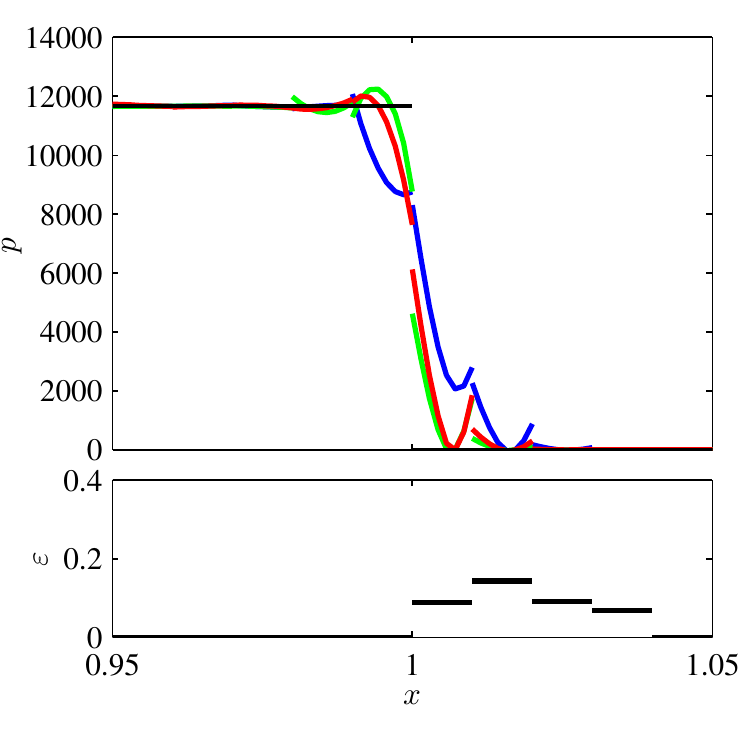}}
        \subfigure[${\rm{Ma}}=2, h = 1/200$.]{\includegraphics[width=0.32\textwidth, clip=, keepaspectratio]{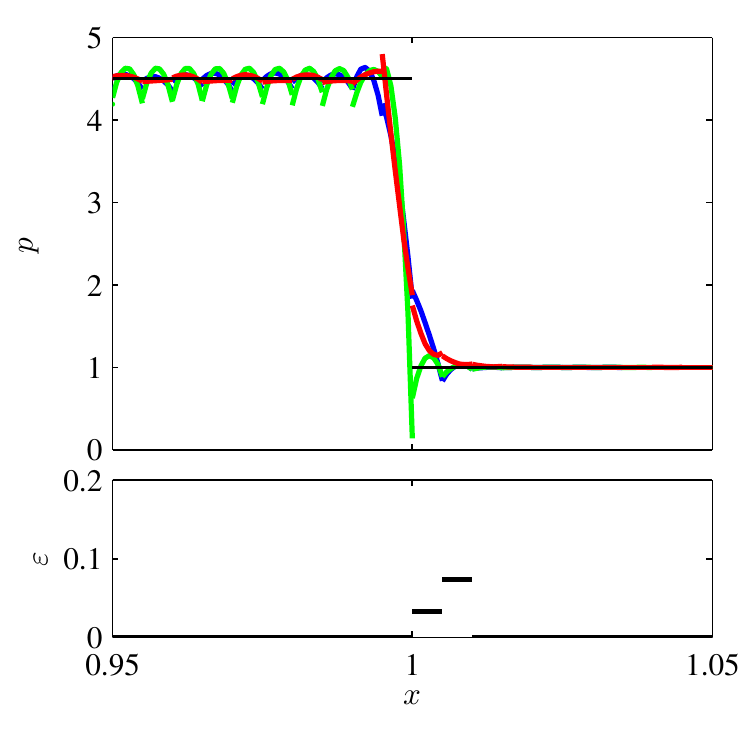}}
        \subfigure[${\rm{Ma}}=5, h = 1/200$.]{\includegraphics[width=0.32\textwidth, clip=, keepaspectratio]{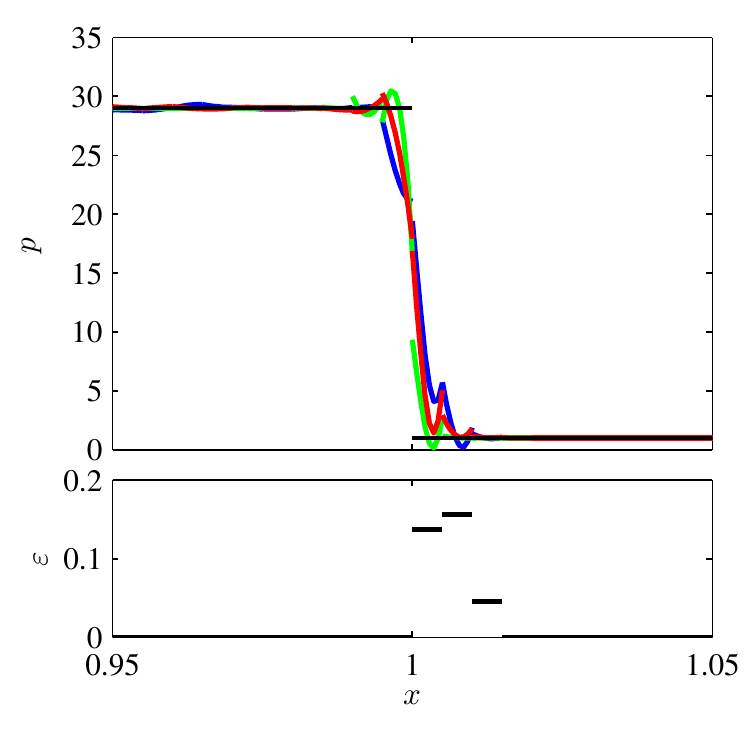}}
        \subfigure[${\rm{Ma}}=100, h = 1/200$.]{\includegraphics[width=0.32\textwidth, clip=, keepaspectratio]{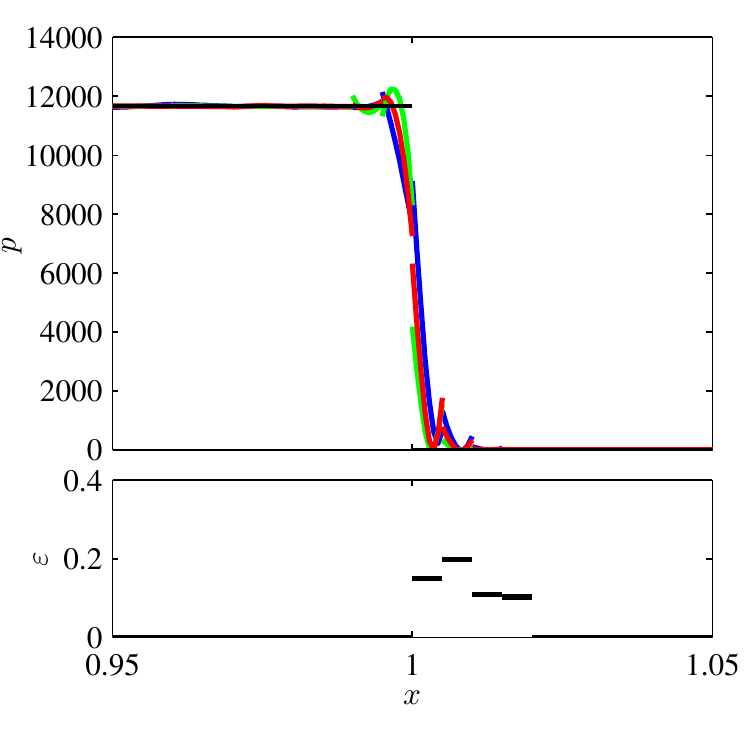}}
    \end{subfigmatrix}
    \caption{\label{1D_SHOCK_CAPTURE}(Color online) DGP2 simulation of moving shock wave for different Mach numbers (Abbreviation: EB$-$entropy bounding; PP$-$positivity preserving~\cite{ZHANG_SHU_JCP2010}; Limiter+PP$-$WENO limiter~\cite{WENOLIM_ZHONG}) with positivity preserving \cite{ZHANG_SHU_JCP2010}.}
\end{figure}

A moving shock-wave in a one-dimensional domain is considered as a test-case for evaluating the robustness and performance of EBDG for shock-capturing. A domain with $x\in[-0.1, 1.1]$ is considered, in which the initial shock front is located at $x=0$. The domain is initialized in $x<0$ with the following 
pre-shock state:
\ben
\rho &=& 1.4\;,\\
u &=& 0\;,\\
p &=& 1\;.
\een
Shocks are specified with different Mach numbers $({\rm{Ma}} = u_s/c)$, and Ma $=\{2, 5, 100\}$ are considered in this case. For all cases considered, the initial value for the entropy, $s_0$, is set to a value of 0.620, corresponding of the minimum value in the initial condition. The simulation ends when the exact solution of the shock front reaches the location at $x=1$. Results are illustrated in Fig.~\ref{1D_SHOCK_CAPTURE}, showing that the entropy boundedness guarantees the robustness and consistent performance over a wide range of shock strengths. Entropy bounding ($\varepsilon \neq 0 $) is only activated in elements that are occupied by flow discontinuities. Compared to the positivity-preserving method, entropy bounding entirely avoids unphysical undershoots in pressure, and provides an improved suppression of oscillations in the post-shock region. Compared to limiting, the entropy bounding shows better robustness in describing shocks at different conditions, introducing lower dissipation in the vicinity of discontinuities. 

\subsection{Two-dimensional flow over a cylinder}
In this section, we verify the convergence order of the EBDG-scheme for high-order curved elements by considering a two-dimensional flow over a round cylinder. The radius of the cylinder is $R=1$ and the far-field boundary is a concentric circle with $R=20$. The condition in the free stream is given as:
\ben
\rho_\infty &=& 1.4\;,\\
u_\infty &=& 5.32\;,\\
v_\infty &=& 0.0\;,\\
p_\infty &=& 1\;.
\een
The corresponding Mach number is 0.38 and characteristic boundary conditions are imposed at the far-field. The entire domain is initialized with free-stream conditions and $s_0 = 0.620$. We compare results on quadrilateral and triangular meshes at three levels of refinement. High-order elements are generated using cubic polynomials to accommodate the curvature of the geometry. The CFL number is set to the CFL$^{\rm EB}$ number from Table \ref{1D_ACCURATE_SSPRK34} for the corresponding shape and polynomial order, multiplied by a factor of 0.8.

A main issue in these simulations is the occurrence of numerical instabilities that are initiated at the leading edge of the cylinder. As a result of this instability, DGP2 and DGP3 without any entropy-bounding diverge (the code blows up) after few iterations. Previously, limiters have been used in this case for stabilizing the transient solutions \cite{HWENOLIM_LUO}. However, for high-order polynomials, it is difficult to develop limiters to achieve the optimal convergence rate without a nontrivial implementation. In contrast, EBDG provides a considerably simpler implementation for enabling high-order simulations for such complex geometric configurations.

\begin{table}[!h!]
	\begin{center}
		\begin{tabular}{|c||cc|cc|cc|}
		\hline 
        \multirow {2}{*}{Mesh} & \multicolumn{2}{c|}{DGP1} & \multicolumn{2}{c|}{DGP2} & \multicolumn{2}{c|}{DGP3} \\
		  & $L_2$-error   & rate & $L_2$-error   & rate&  $L_2$-error   & rate\\
		\hline\hline
        \multicolumn{7}{|>{\columncolor[gray]{0.7}}c|}{Quadrilateral Elements}       \\ \hline

         Level 1 &	7.272\e{-2} & {-} & 1.694\e{-2} & {-} & 3.816\e{-3} & {-}\\
		 Level 2 &	1.318\e{-2} & 2.464  & 7.219\e{-4} & 4.552  & 1.827\e{-4} & 4.384 \\
		 Level 3 &	2.441\e{-3} & 2.433  & 6.029\e{-5} & 3.582  & 1.036\e{-5} & 4.141 \\
		\hline
		\hline
        \multicolumn{7}{|>{\columncolor[gray]{0.7}}c|}{Triangular Elements}       \\ \hline
         Level 1 &	1.137\e{-1} & {-} & 2.590\e{-2} & {-} & 4.086\e{-3} & {-}\\
		  Level 2 &	1.865\e{-2} & 2.608  & 8.899\e{-4} & 4.863  & 1.291\e{-4} & 4.984 \\
		  Level 3 &	3.391\e{-3} & 2.459  & 7.222\e{-5} & 3.623  & 6.939\e{-6} & 4.217 \\
		\hline
		\end{tabular}
	\end{center}
	\caption{\label{2D_ACCURATE_RK45}Comparisons of convergence rate for 2D flow over a cylinder (here we use entropy to evaluate the error).}
\end{table}

\begin{figure}[!htb!]
\centering
     \begin{tabular}{|c|c|c|c|}
     \hline
     & Mesh Level 1  & Mesh Level 2 & Mesh Level 3 \\ 
     \hline
     \begin{sideways}\hspace*{20mm} Mesh  \end{sideways}
     & \includegraphics[width=0.28\textwidth, clip=, keepaspectratio]{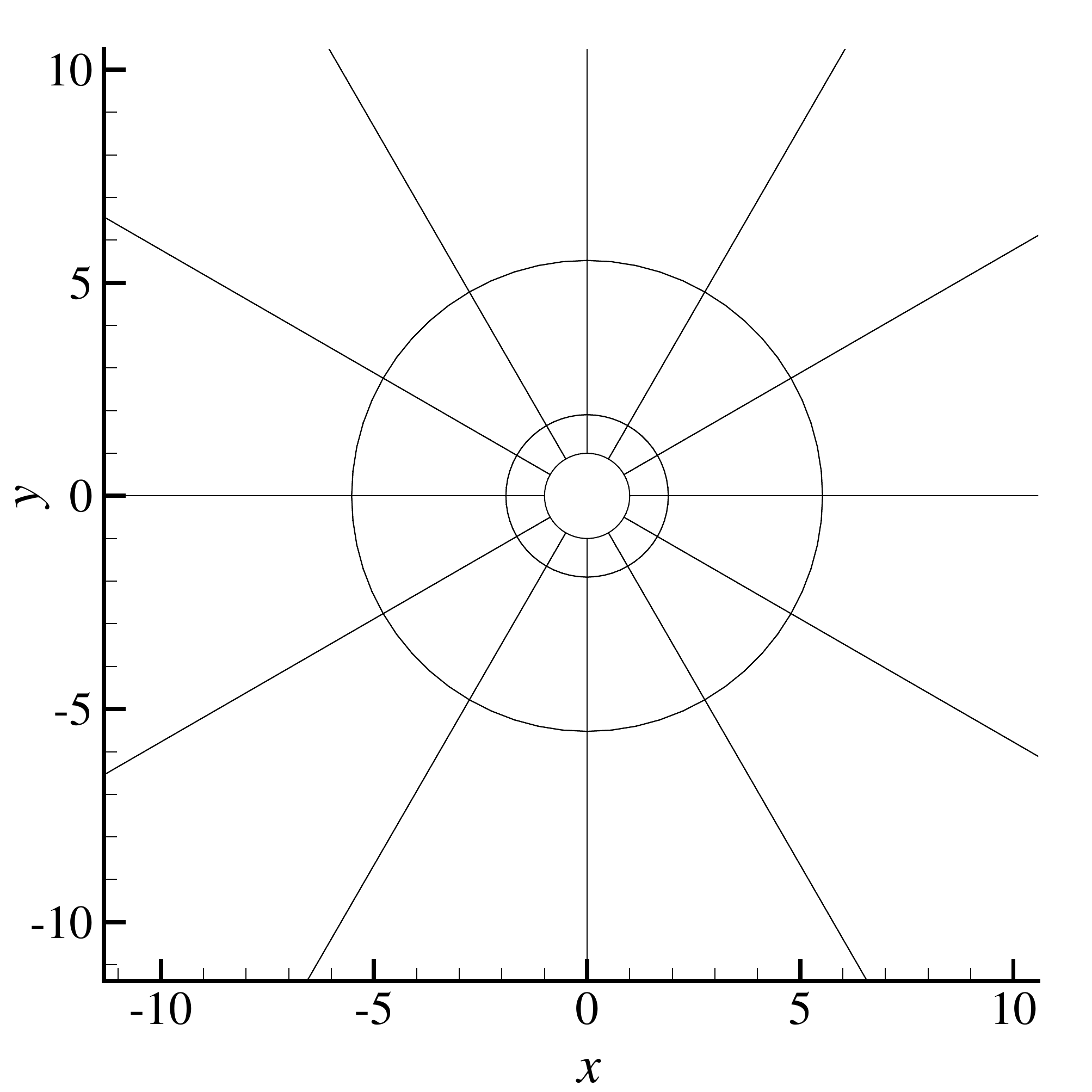}
     & \includegraphics[width=0.28\textwidth, clip=, keepaspectratio]{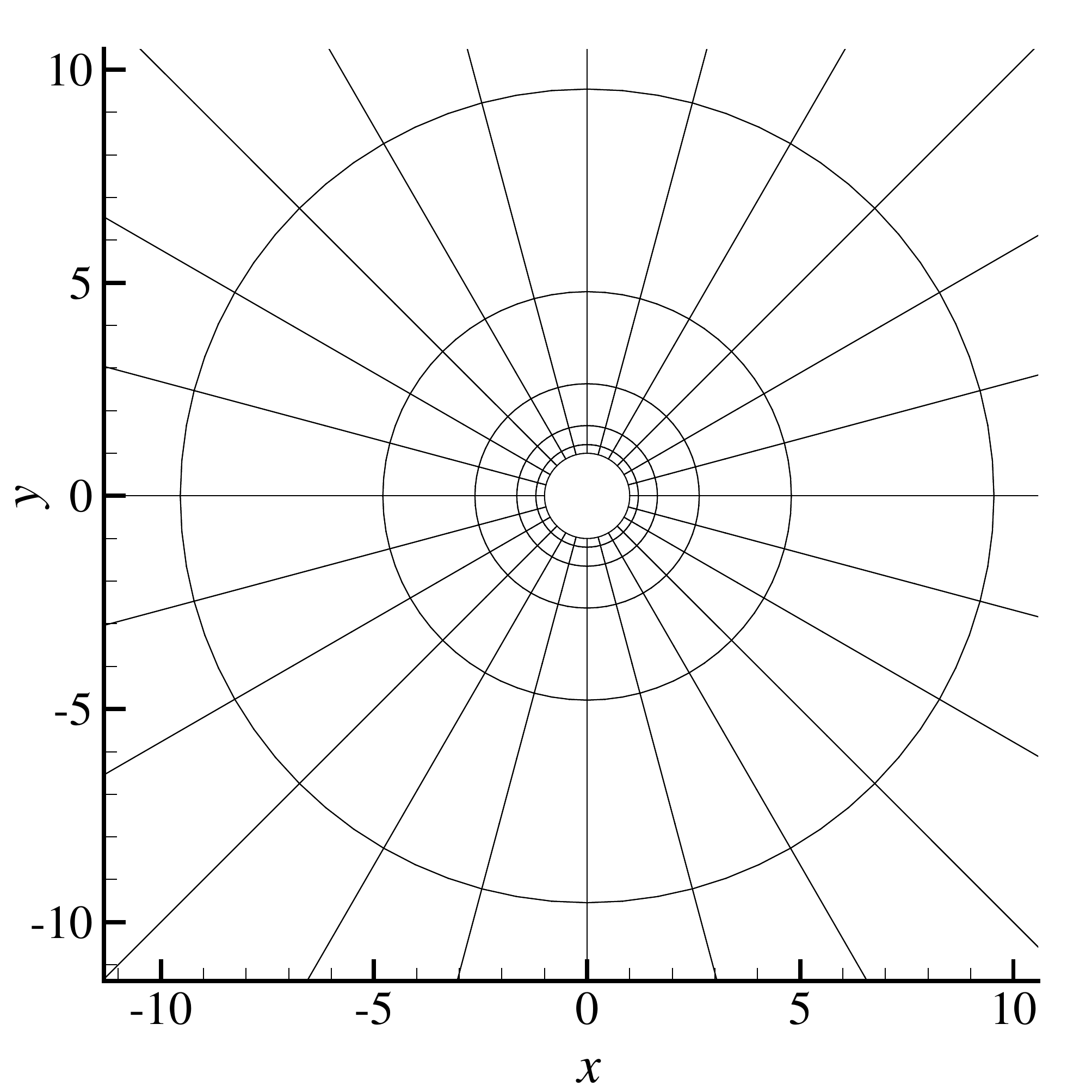}
     & \includegraphics[width=0.28\textwidth, clip=, keepaspectratio]{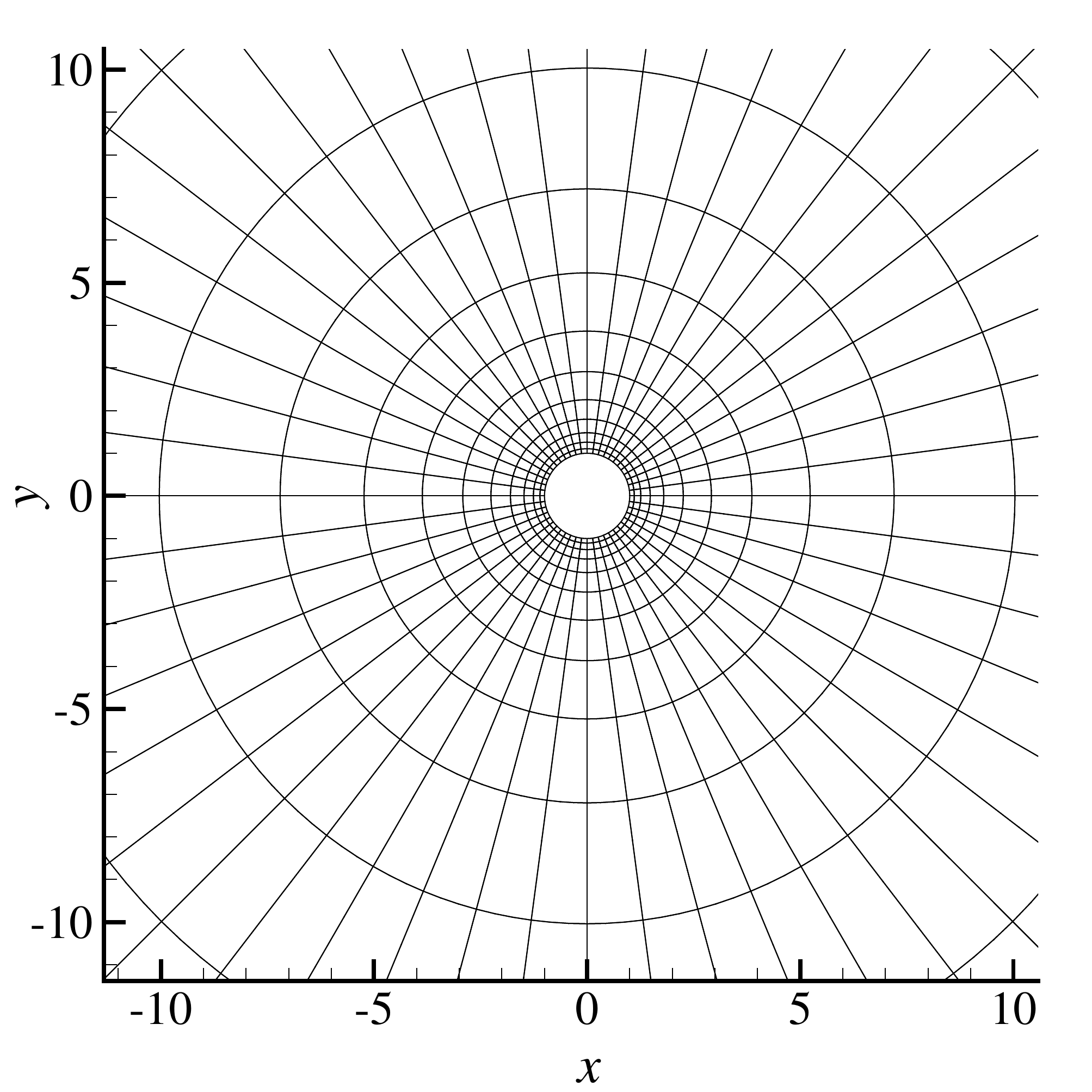} \\
     \hline
      \begin{sideways}\hspace*{14mm} DGP2 solutions  \end{sideways}
     & \includegraphics[width=0.28\textwidth, clip=, keepaspectratio]{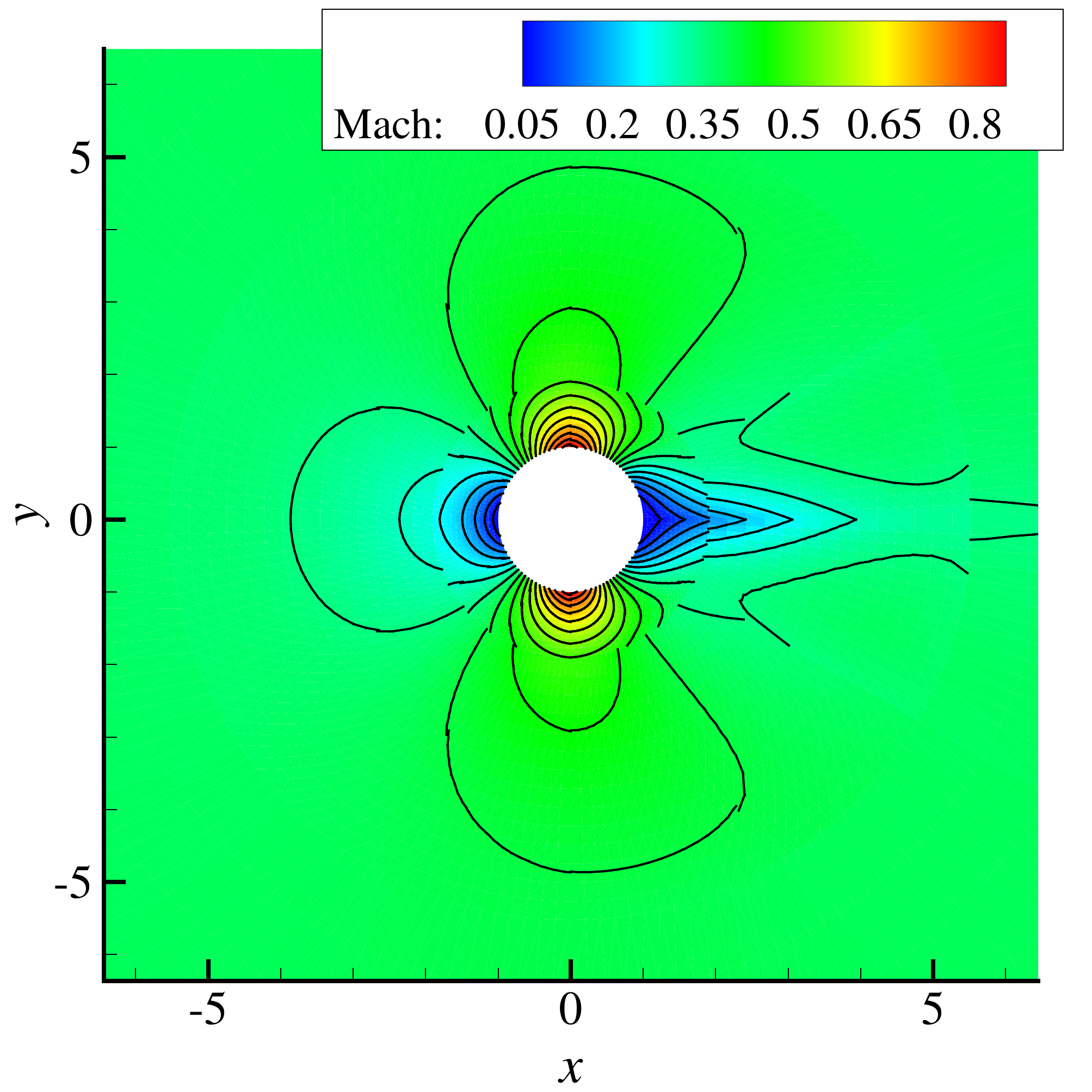}
     & \includegraphics[width=0.28\textwidth, clip=, keepaspectratio]{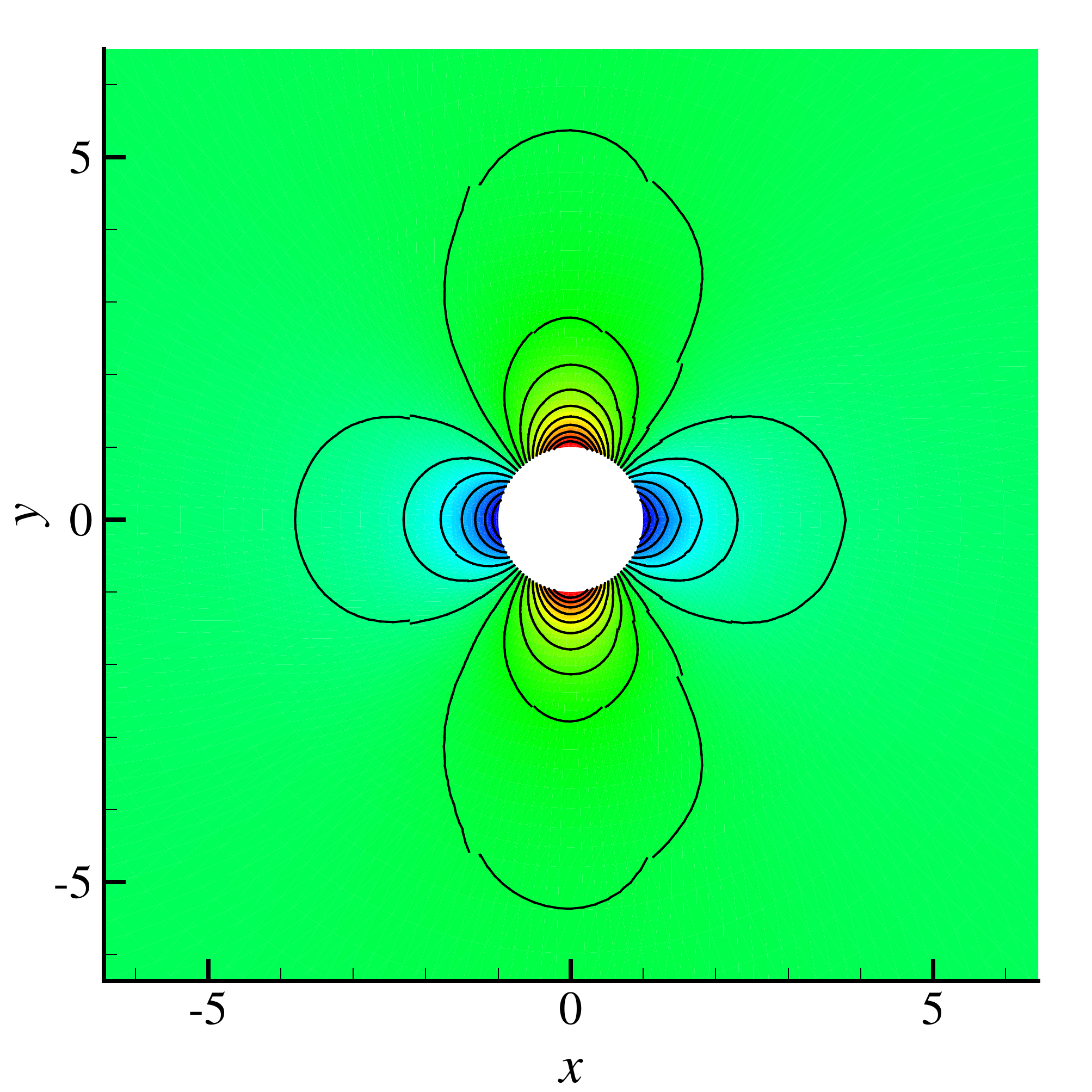}
     & \includegraphics[width=0.28\textwidth, clip=, keepaspectratio]{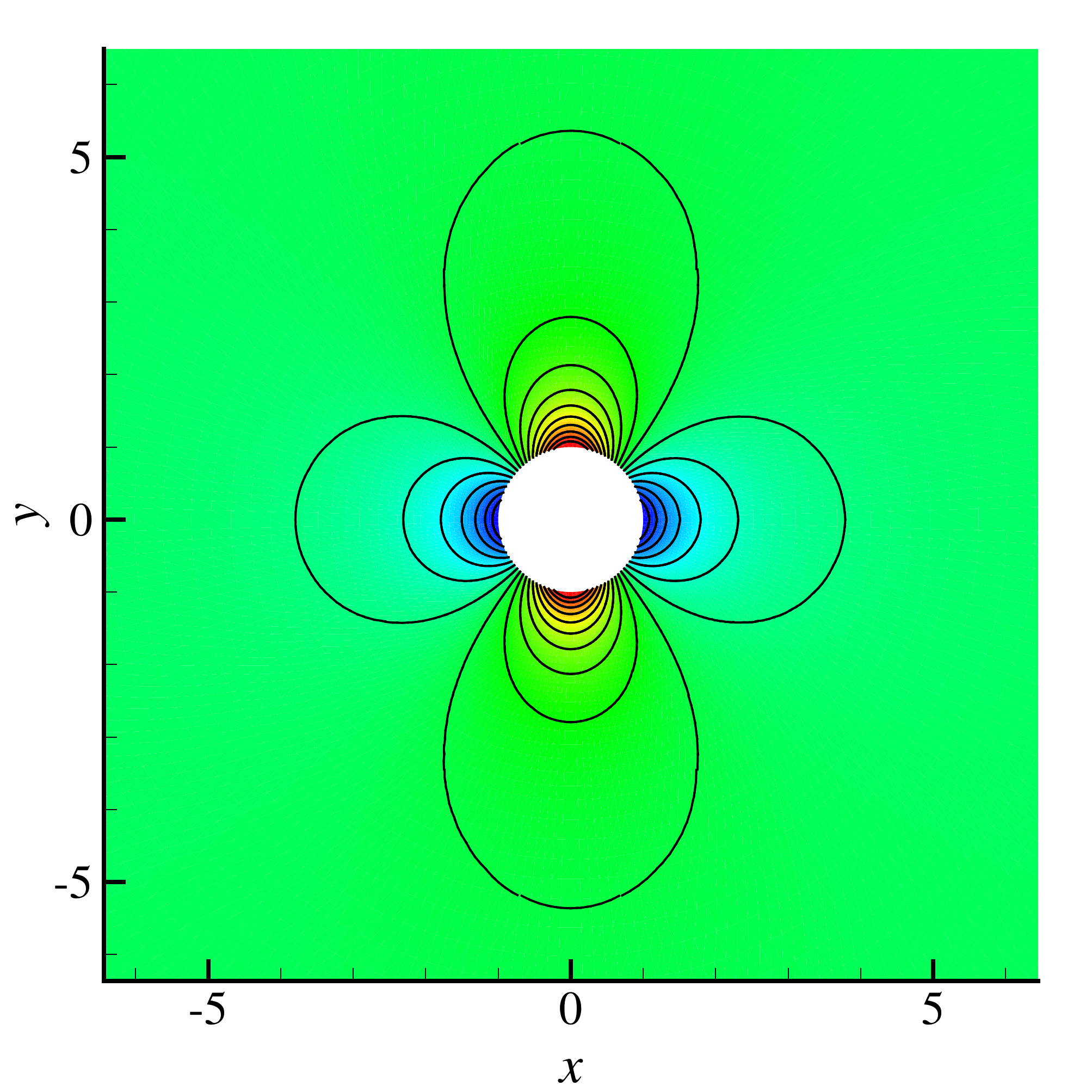} \\
     \hline
     \begin{sideways}\hspace*{14mm} Convergence  \end{sideways}
     & \includegraphics[width=0.28\textwidth, clip=, keepaspectratio]{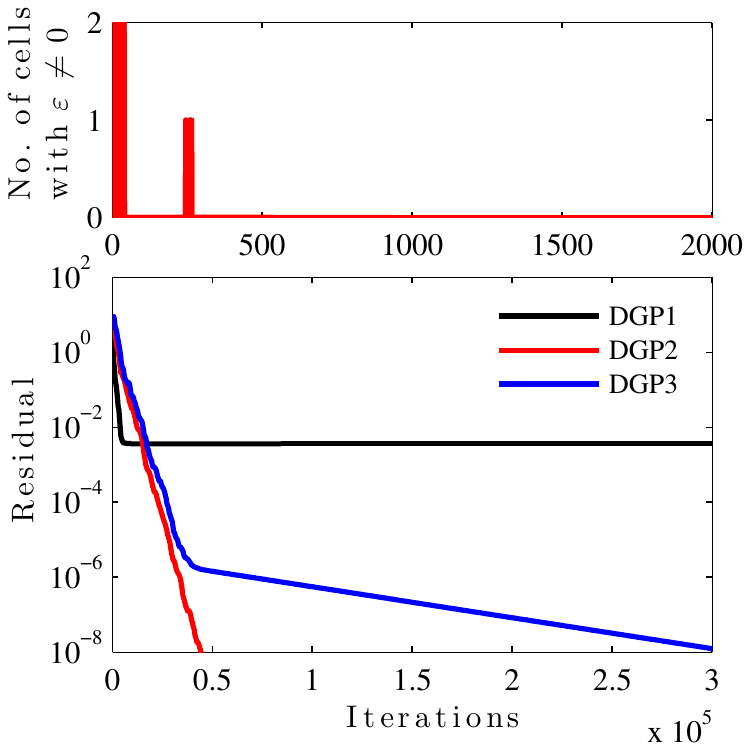}
     & \includegraphics[width=0.28\textwidth, clip=, keepaspectratio]{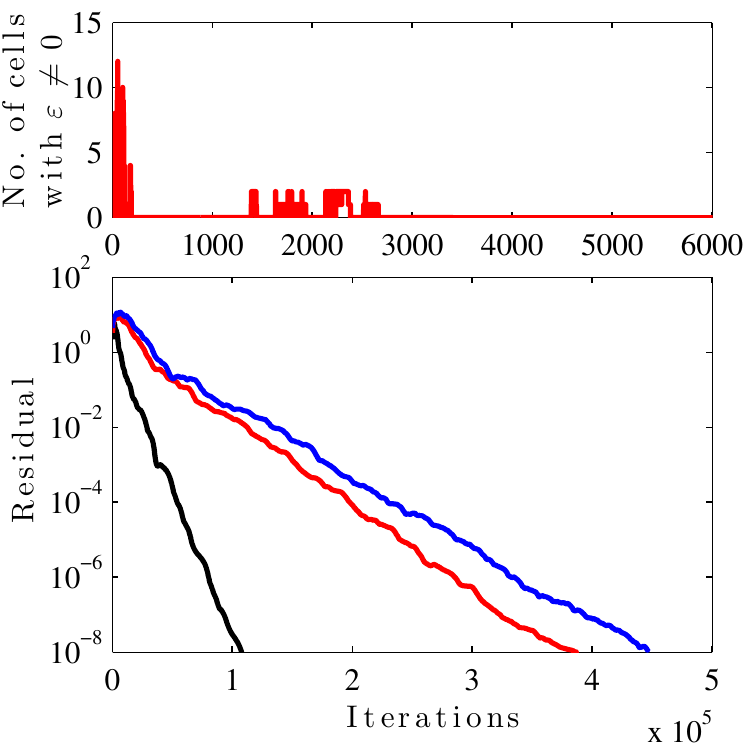}
     & \includegraphics[width=0.28\textwidth, clip=, keepaspectratio]{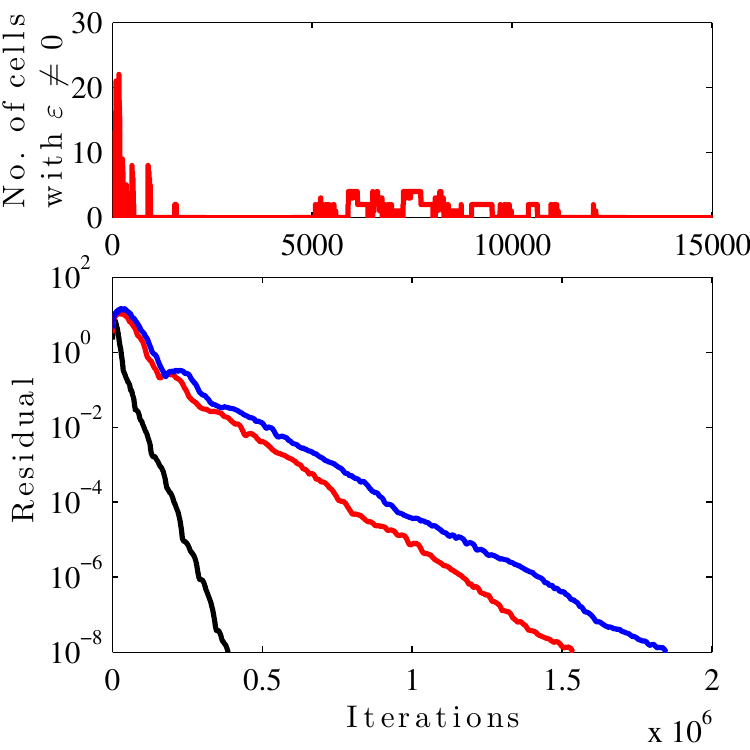} \\
     \hline
    \end{tabular}
    \caption{\label{FLOW_OVER_CIRCLE_QUAD}EBDG-solution of flow over a cylinder on curved quadrilateral meshes with
    three different refinement levels; top: computational mesh in near-field of the cylinder; middle: Mach number; bottom: convergence history and activation of entropy bounding as a function of iteration.}
\end{figure}

\begin{figure}[!htb!]
\centering
     \begin{tabular}{|c|c|c|c|}
     \hline
     & Mesh Level 1  & Mesh Level 2 & Mesh Level 3 \\ 
     \hline
     \begin{sideways}\hspace*{20mm} Mesh  \end{sideways}
     & \includegraphics[width=0.28\textwidth, clip=, keepaspectratio]{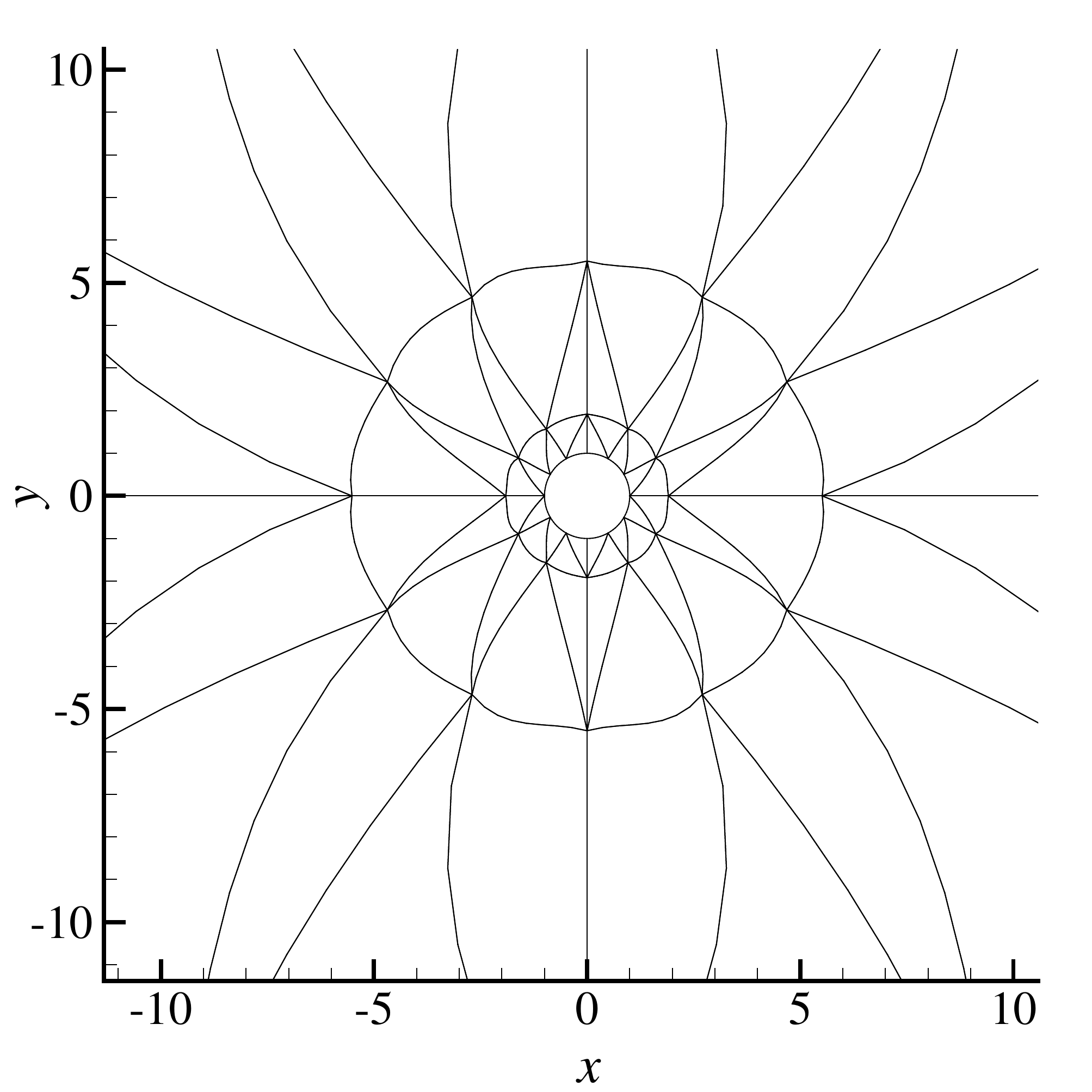}
     & \includegraphics[width=0.28\textwidth, clip=, keepaspectratio]{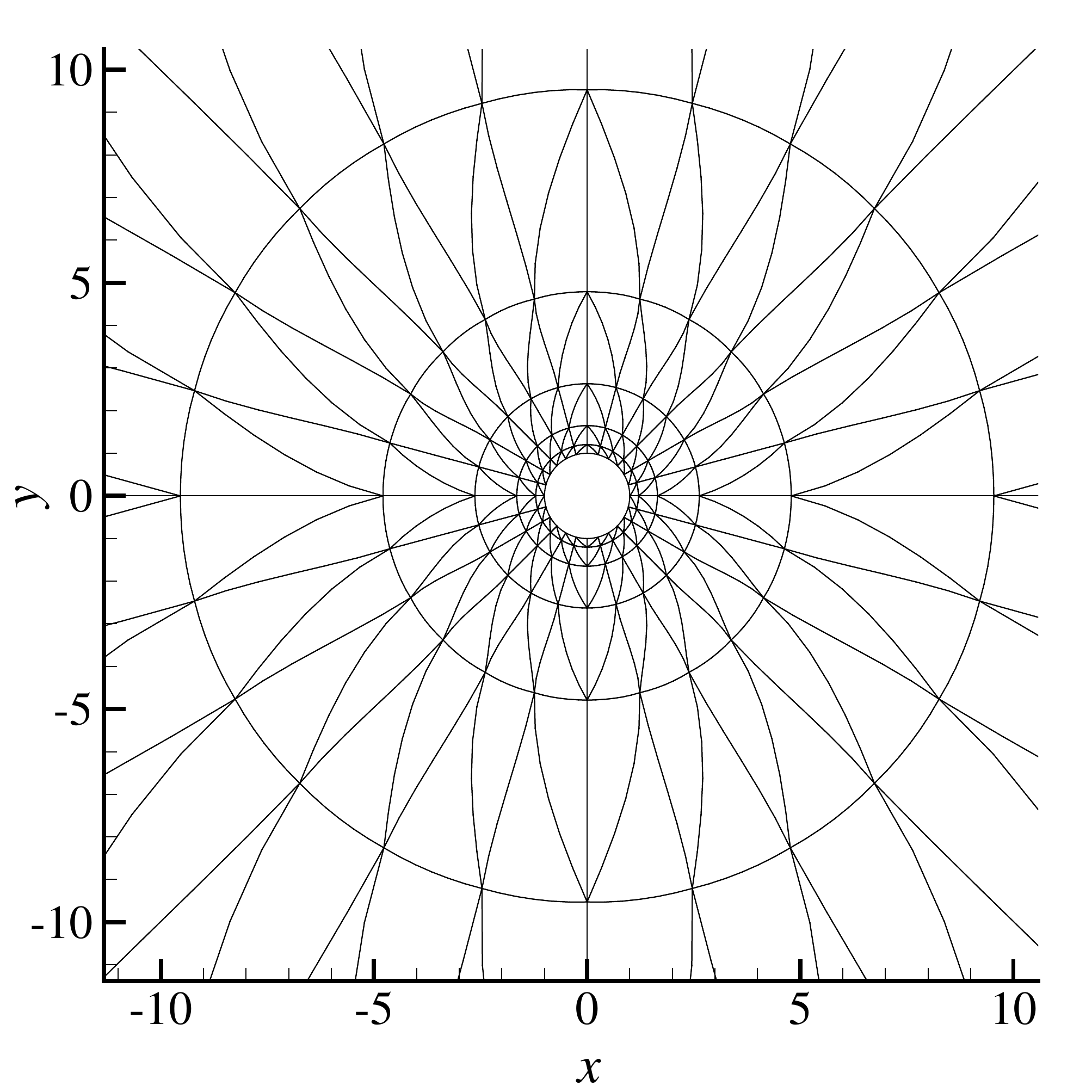}
     & \includegraphics[width=0.28\textwidth, clip=, keepaspectratio]{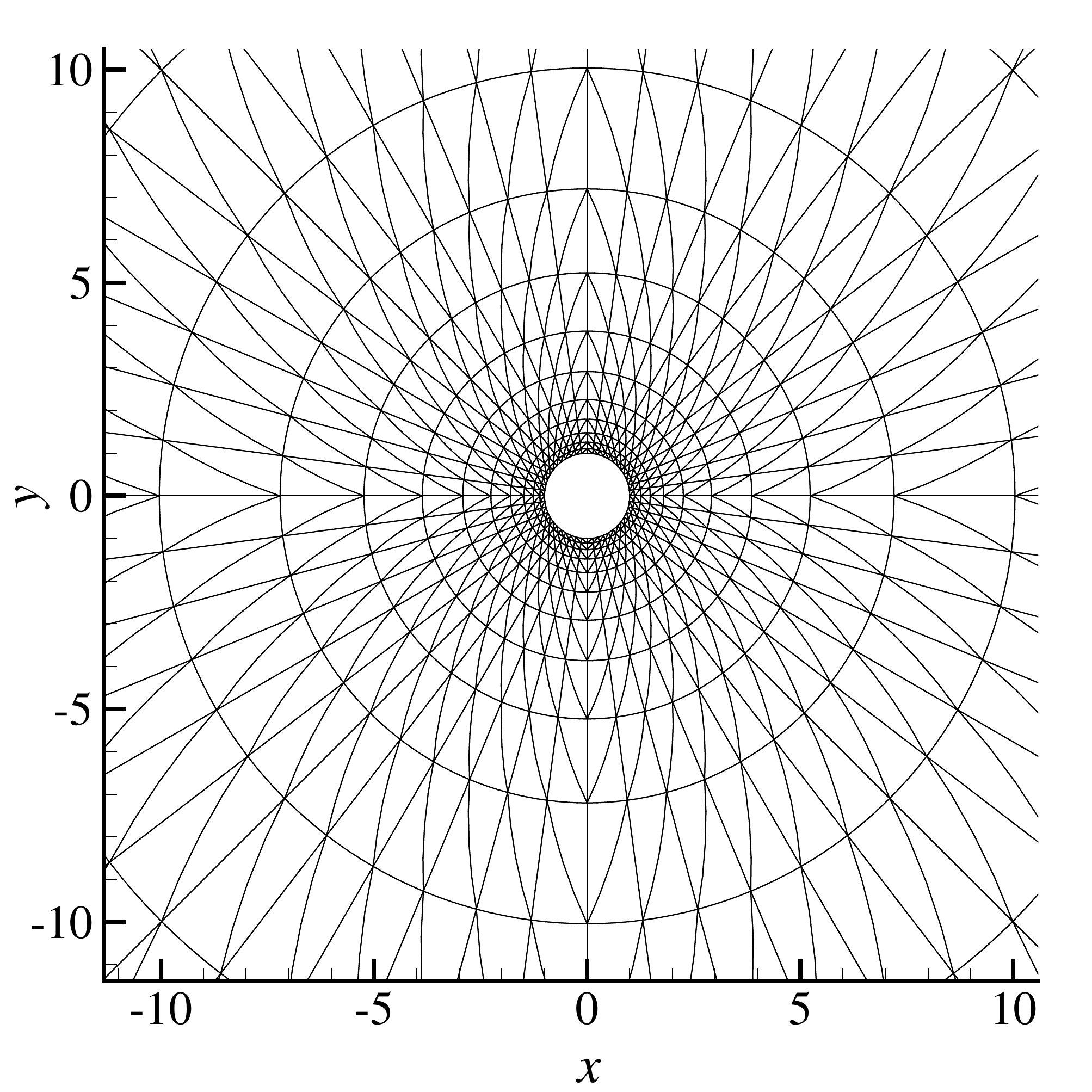} \\
     \hline
      \begin{sideways}\hspace*{14mm} DGP2 solutions  \end{sideways}
     & \includegraphics[width=0.28\textwidth, clip=, keepaspectratio]{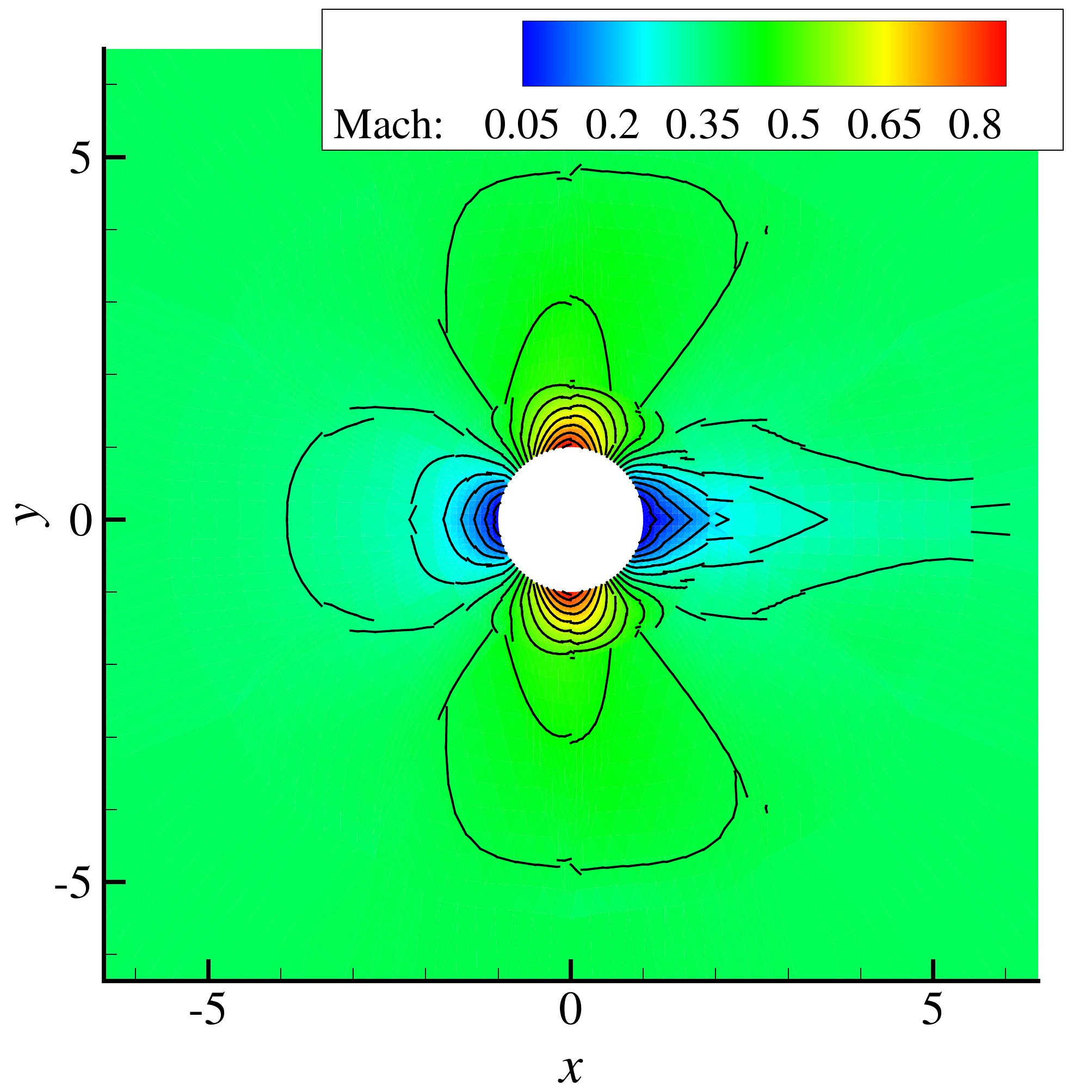}
     & \includegraphics[width=0.28\textwidth, clip=, keepaspectratio]{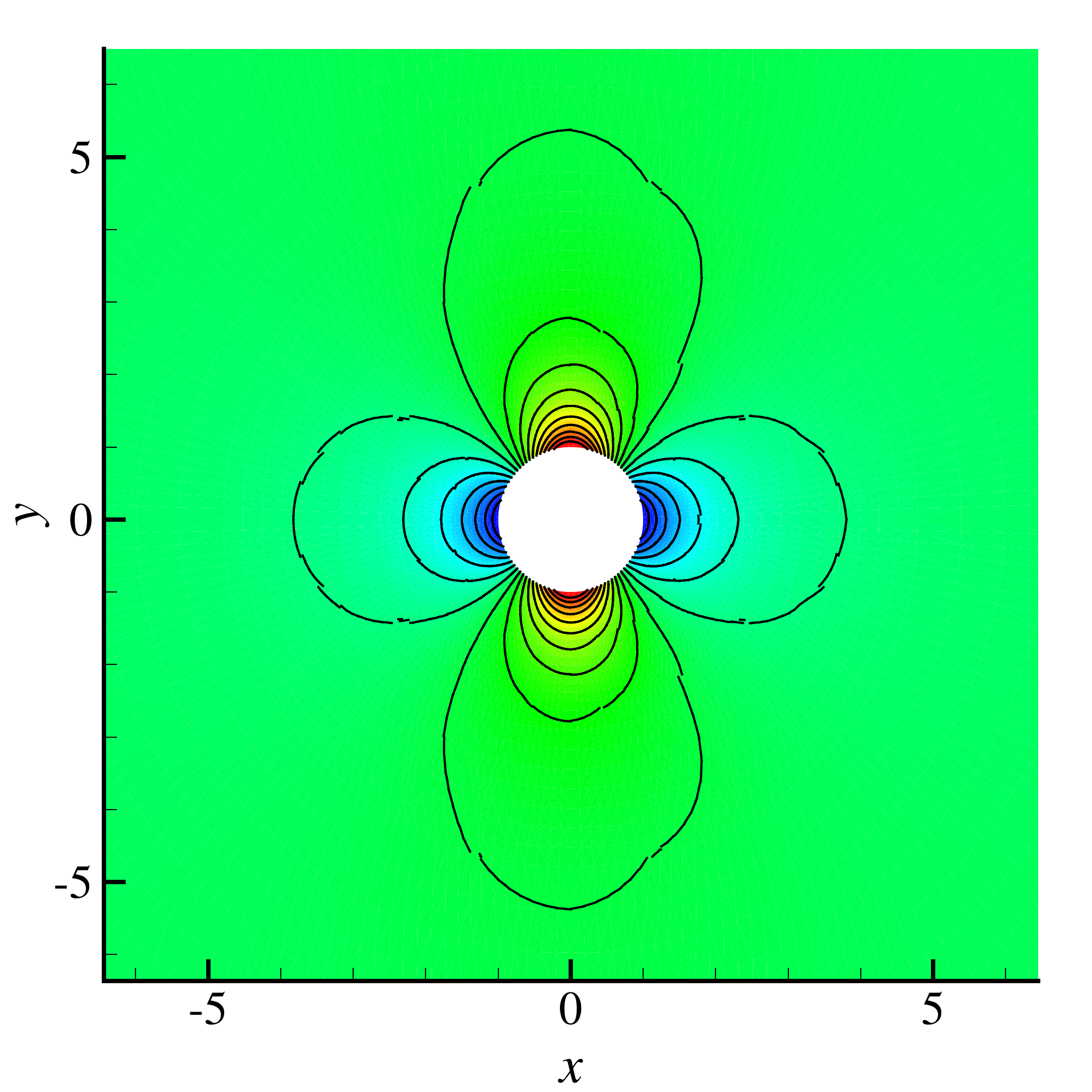}
     & \includegraphics[width=0.28\textwidth, clip=, keepaspectratio]{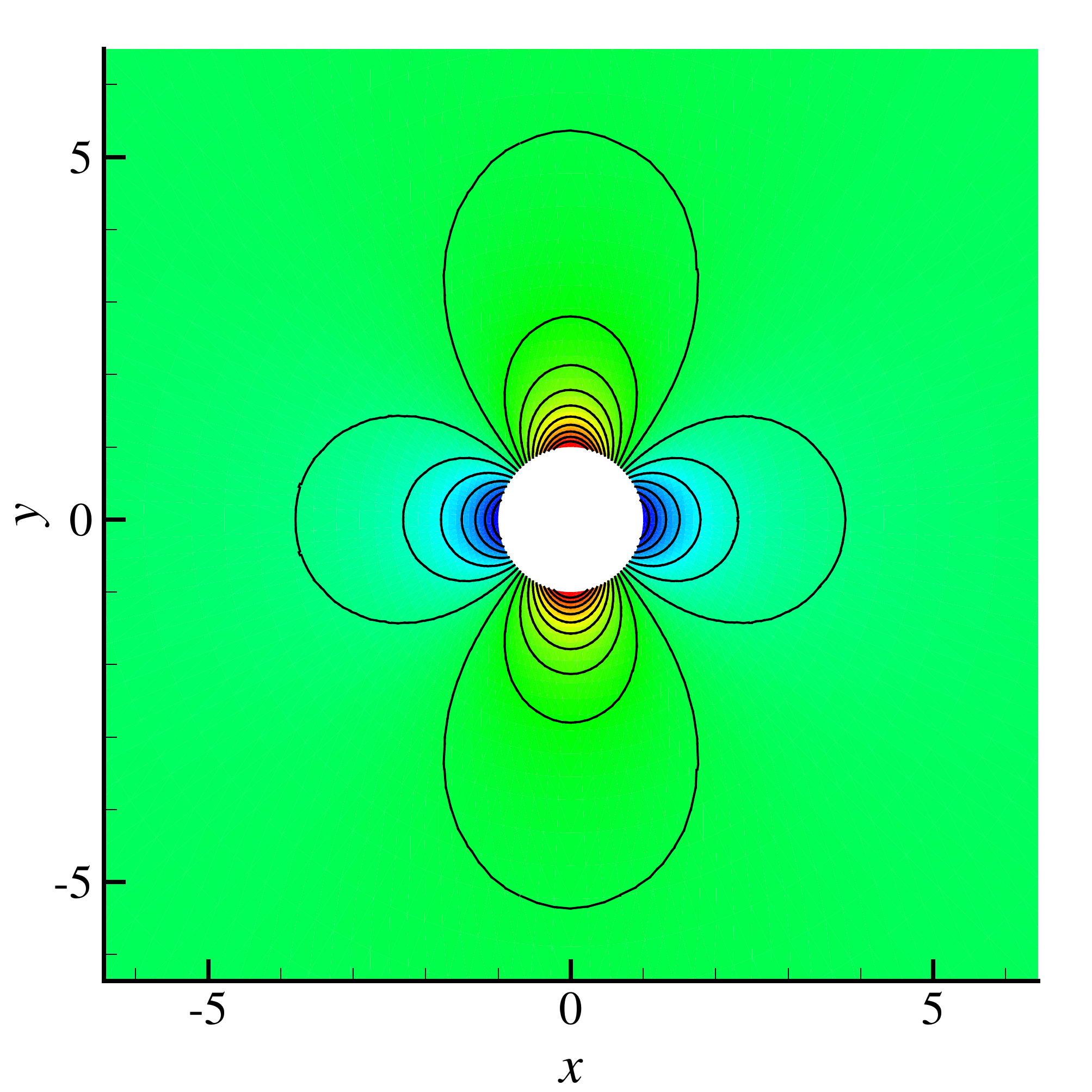} \\
     \hline
     \begin{sideways}\hspace*{14mm} Convergence  \end{sideways}
     & \includegraphics[width=0.28\textwidth, clip=, keepaspectratio]{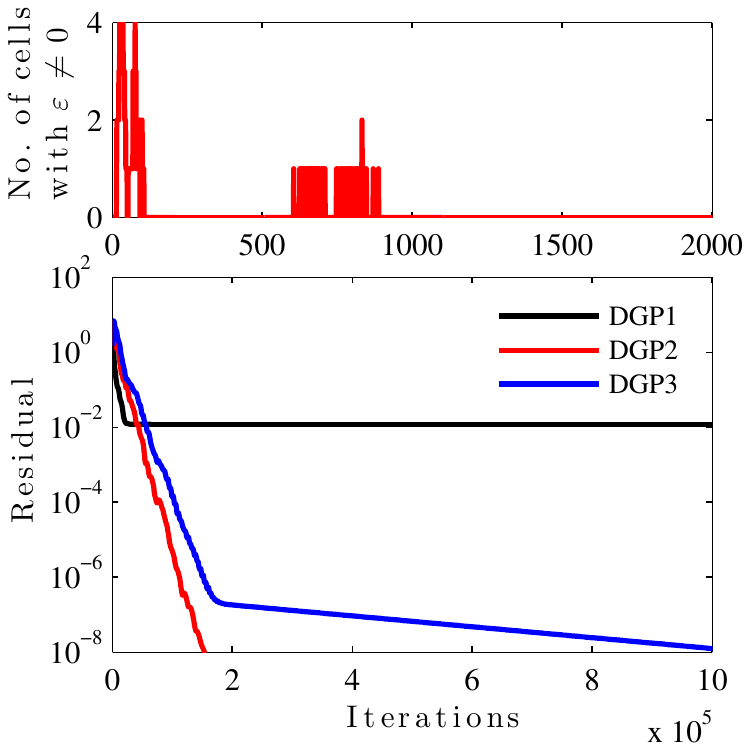}
     & \includegraphics[width=0.28\textwidth, clip=, keepaspectratio]{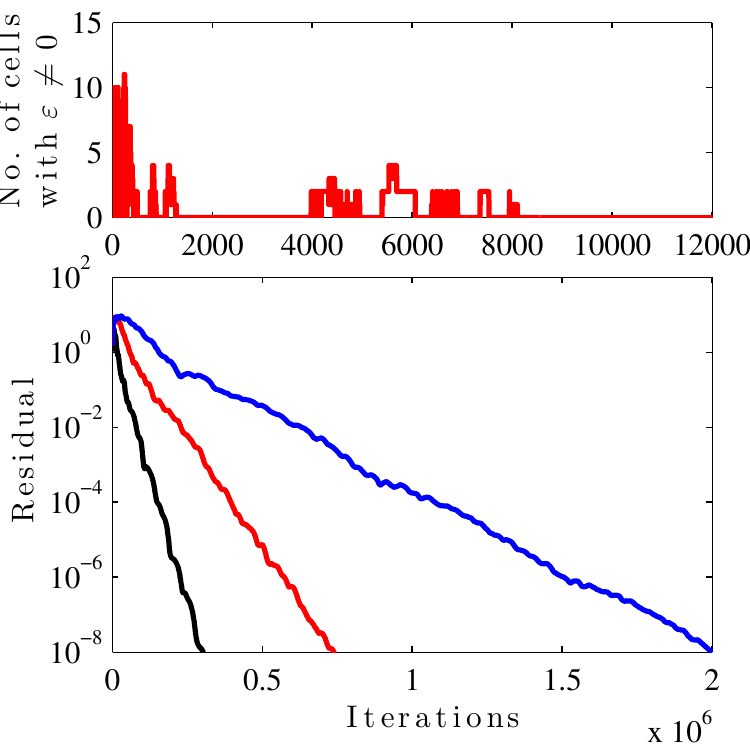}
     & \includegraphics[width=0.28\textwidth, clip=, keepaspectratio]{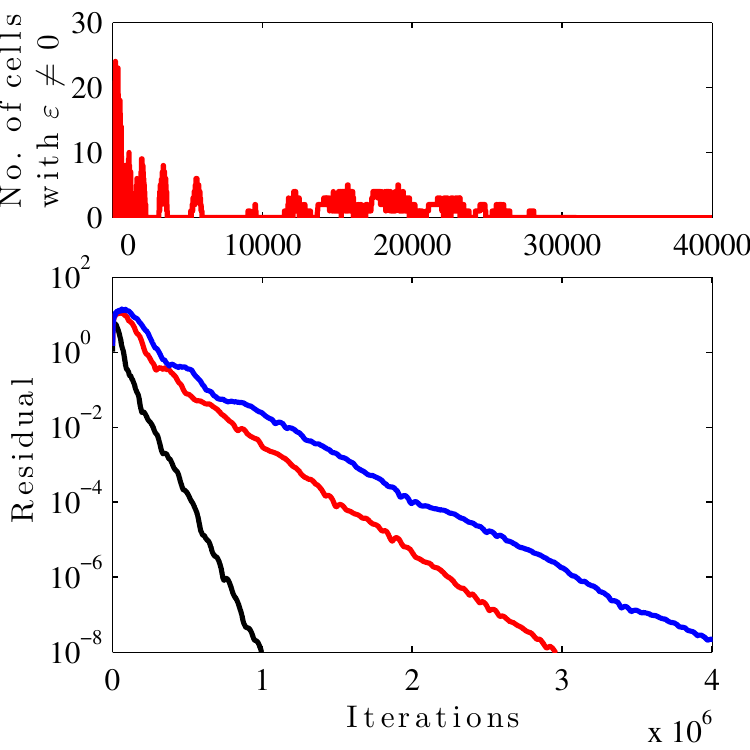} \\
     \hline
    \end{tabular}
    \caption{\label{FLOW_OVER_CIRCLE_TRI}DG-solution of flow over a cylinder on curved triangular meshes with
    three different refinement levels; top: computational mesh in the near-field of the cylinder; middle: Mach number; bottom: convergence history and activation of entropy bounding as a function of iteration.}
\end{figure}

Comparisons of the computational meshes, simulation results, and convergence properties are presented in Figs.~\ref{FLOW_OVER_CIRCLE_QUAD}  and \ref{FLOW_OVER_CIRCLE_TRI}. It is evident that the solution is improved by increasing the mesh resolution. The convergence history of the residual, provided in the last row of both figures, shows that entropy bounding is mostly activated during the start-up phase of the simulation to suppress numerical oscillations and ensure stability. It is interesting to note that the number of elements that require bounding is restricted to the region near the stagnation point upstream of the cylinder, and is confined to less than 8\% of the total number of elements. As the solution converges to the steady-state condition, entropy-bounding remains deactivated, retaining the high-order accuracy. Since the solution is smooth the physical entropy production is zero. Therefore, the convergence rates are measured in terms of entropy error using the discrete $L_2$-norm. A comparison of the convergence rates are presented in Table~\ref{2D_ACCURATE_RK45}, confirming that the optimal convergence rate is preserved even for complex geometries with curved elements.

\subsection{Two-dimensional double Mach reflection}
This test case is designed to assess the performance of EBDG for simulation of flows with strong shocks and wave structures. The numerical setup follows this of  Woodward and Colella \cite{DMR_1984}, representing a Mach 10 shock over a $30^\circ$-wedge. All quantities are non-dimensionalized, and the computational domain is $[0, 4] \times [0, 1]$. In the present study, we consider two different mesh-discretizations, consisting of a Cartesian mesh with quadratic elements ($L_e=h=0.02$) and a mesh with triangular elements ($L_e\approx0.02$). The pre-shock state is the same as that in Sec.~\ref{1D_MOV_SHOCK} and hence $s_0$ is set to 0.620. The CFL number is prescribed from Table \ref{1D_ACCURATE_SSPRK34} using a safety factor of 0.8.

\begin{figure}
     \begin{subfigmatrix}{2}
        \subfigure[EBDGP1, quadrilateral mesh with $h=0.02.$]{\includegraphics[width=0.475\textwidth, clip=, keepaspectratio]{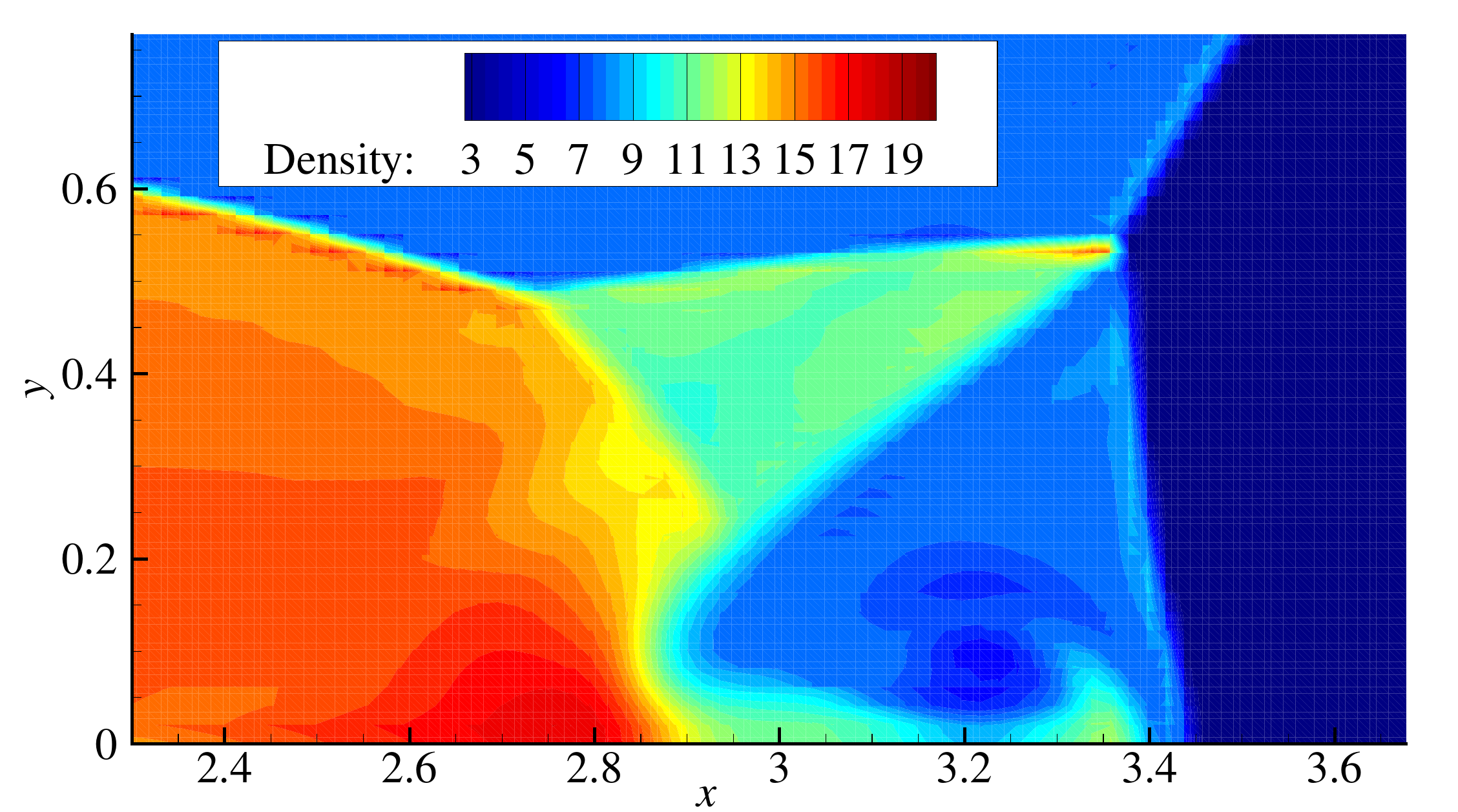}}
        \subfigure[EBDGP1, triangular mesh with $h=0.02.$]{\includegraphics[width=0.475\textwidth, clip=, keepaspectratio]{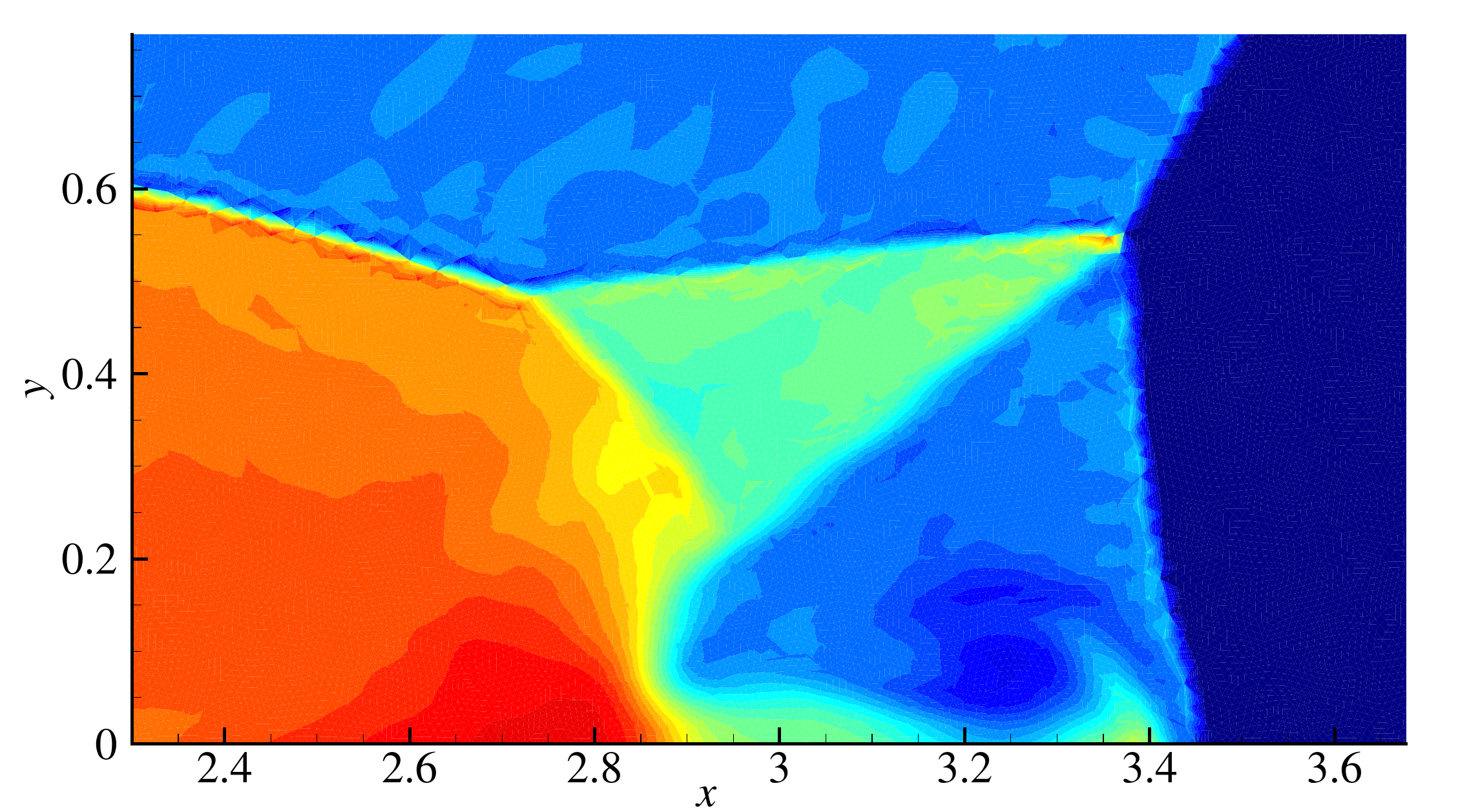}}
        \subfigure[EBDGP2, quadrilateral mesh with $h=0.02.$]{\includegraphics[width=0.475\textwidth, clip=, keepaspectratio]{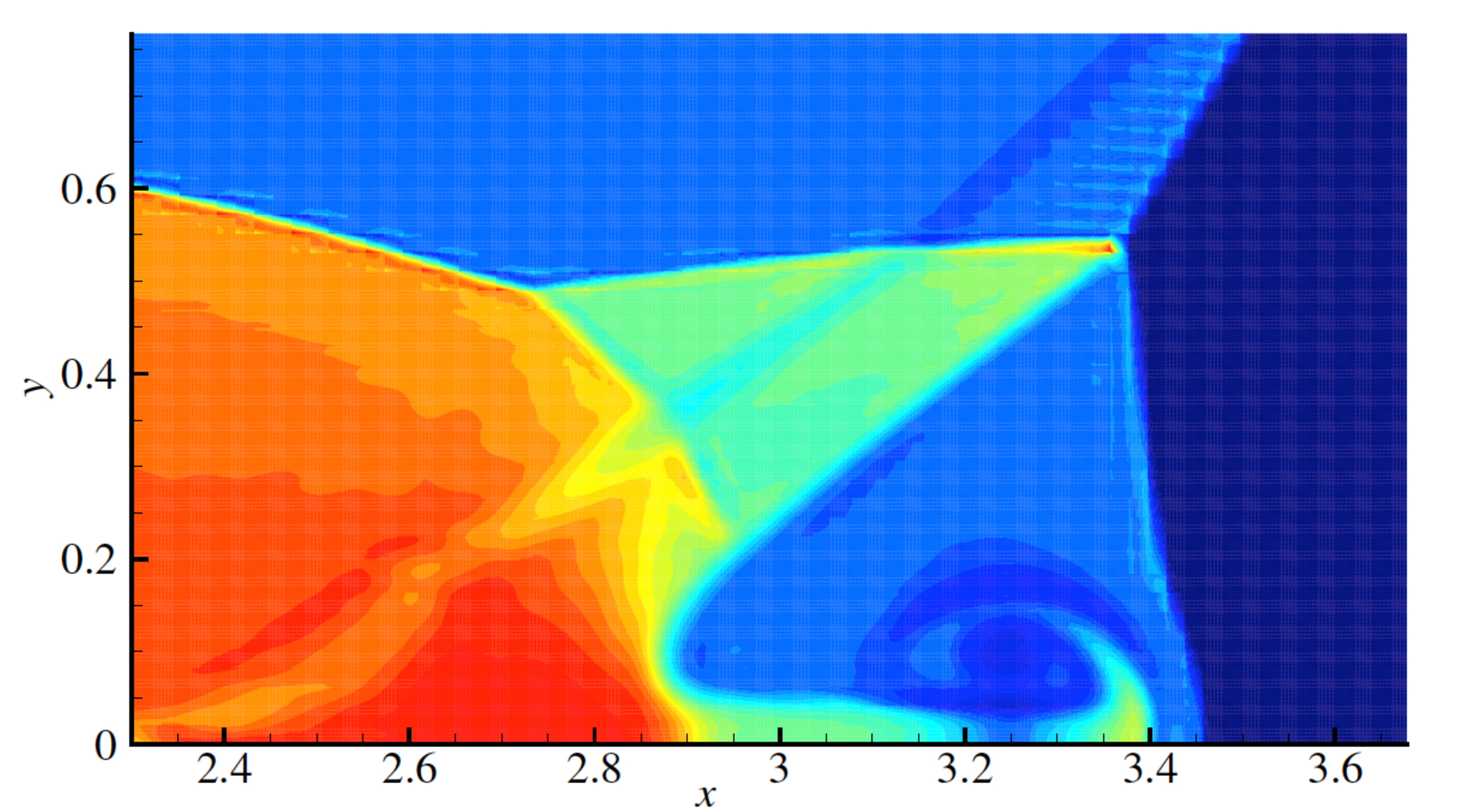}}
        \subfigure[EBDGP2, triangular mesh with $h=0.02.$]{\includegraphics[width=0.475\textwidth, clip=, keepaspectratio]{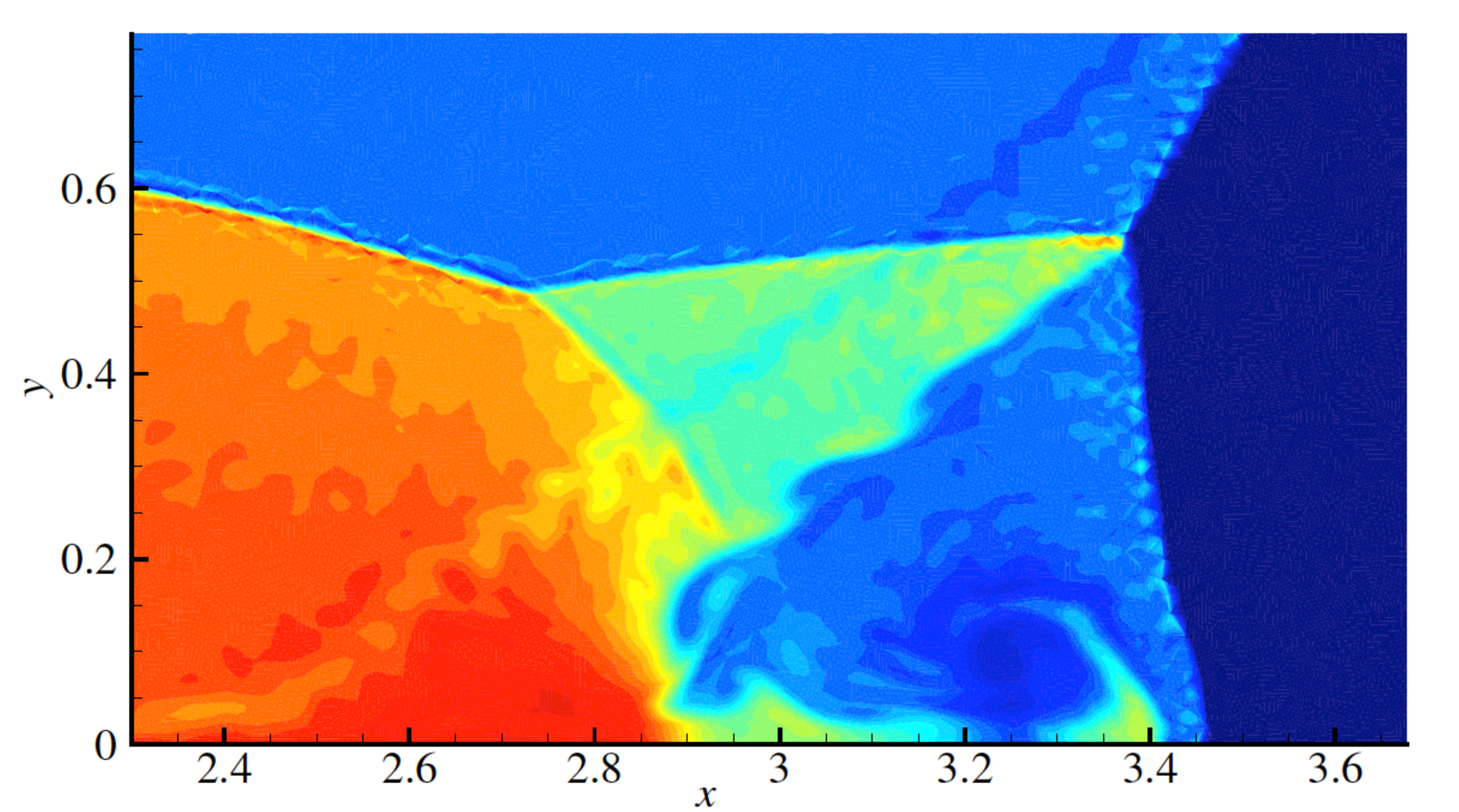}}
        \subfigure[\label{DMR_SOLUTIONS_WENO5}WENO5, quadrilateral mesh with $h=0.0067.$]{\includegraphics[width=0.475\textwidth, clip=, keepaspectratio]{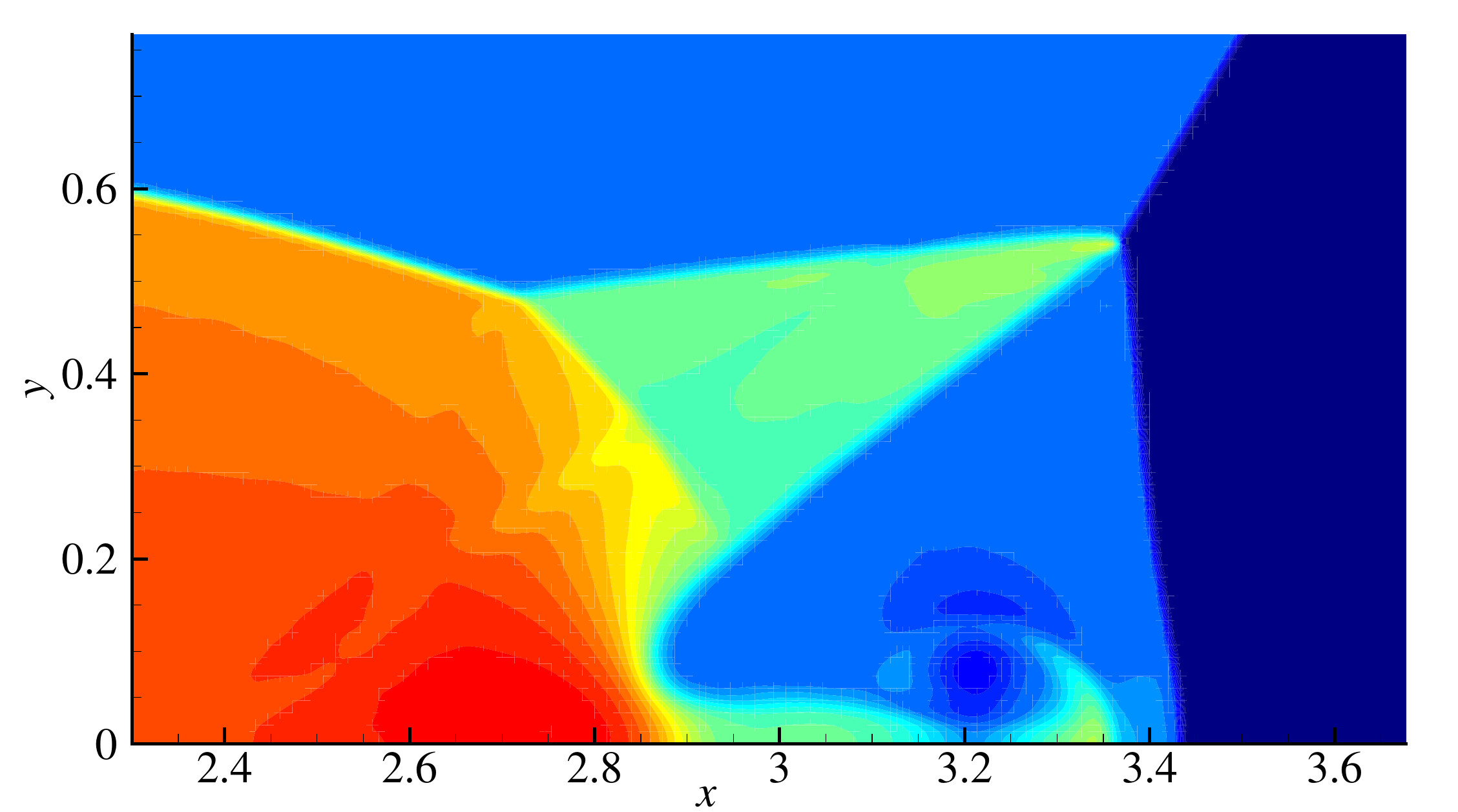}}
        \subfigure[\label{DMR_SOLUTIONS_DG_WENO_LIMITER}DGP2 with WENO limiter~\cite{WENOLIM_ZHONG}, quadrilateral mesh with $h=0.02$.]{\includegraphics[width=0.475\textwidth, clip=, keepaspectratio]{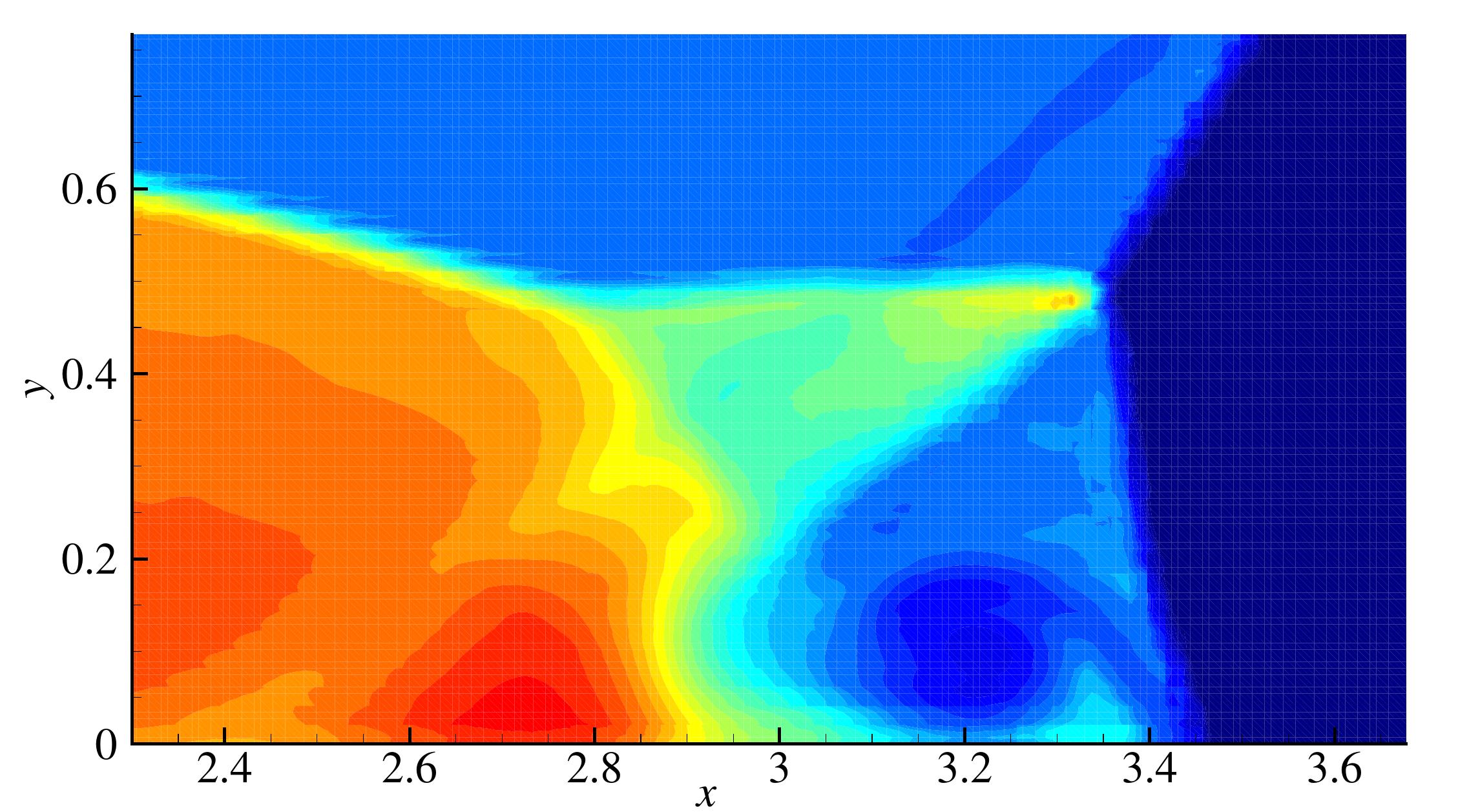}}
    \end{subfigmatrix}
   \begin{subfigmatrix}{1}
    \subfigure[Instantaneous snapshot of $\varepsilon$ for DGP2 on quadrilateral mesh.]{\includegraphics[width=0.85\textwidth, clip=, keepaspectratio]{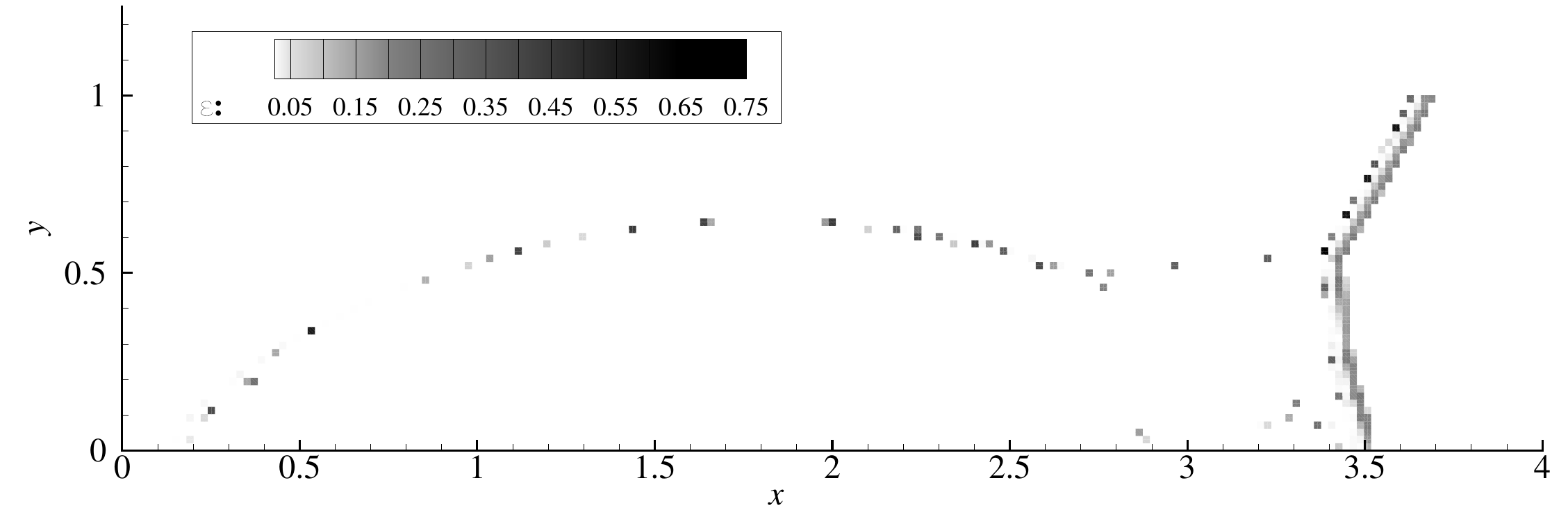}}
   \end{subfigmatrix}
    \caption{\label{DMR_SOLUTIONS}(Color online) Simulation results of double Mach reflection over a 30$^{\rm o}$-wedge.}
\end{figure}

Simulation results for density contours at time $t=0.25$ are shown in Fig.~\ref{DMR_SOLUTIONS}. The proposed EBDG-method captures all wave-features, and it is found that without enforcing the entropy constraint the solution diverges in the first iteration for these strong shock conditions. For comparison, a reference solution obtained using a fifth-order WENO-scheme is shown in Fig.~\ref{DMR_SOLUTIONS_WENO5}, and results from a DGP2-simulation using a WENO-limiter~\cite{WENOLIM_ZHONG} are presented in Fig.~\ref{DMR_SOLUTIONS_DG_WENO_LIMITER}. Comparisons between EBDGP1 and EBDGP2 results show the benefit of the high-order scheme in providing improved representations of the shock-wave structure. At the same degrees of freedom, the DGP2-solution provides comparable predictions to that of the fifth-order WENO scheme, except for the small oscillations that cannot be removed by the linear scaling procedure. Compared to the DG-simulation with WENO-limiter (Fig.~\ref{DMR_SOLUTIONS_DG_WENO_LIMITER}), EBDG effectively avoids introducing excessive numerical dissipation since the solution is only entropy-constrained in regions in which the entropy condition is violated. 

\subsection{Three-dimensional supersonic flow over a sphere}
This test case extends the evaluation of the EBDG-method to three-dimensional configurations with complex geometries. Currently, robust approaches for capturing strong shocks in three-dimensional curved elements are still subject to investigation. This test case considers a flow at a Mach number of 6.8 over a sphere. The radius of the sphere is $R=1$. Due to the geometric symmetry, the computational domain considers only an eighth section of the domain, and it extends to $3R$ in radial direction. Symmetry boundary conditions are imposed at the planes $y=0$ and $z=0$, and outflow boundary conditions are prescribed at $x=0.$ Normal velocity inflow is prescribed at the outer shell with the following specification:
\ben
\rho_\infty &=& 1.4\;,\\
u_\infty &=& -6.80\;,\\
v_\infty &=& 0.0\;,\\
w_\infty &=& 0.0\;,\\
p_\infty &=& 1\;.
\een
Slip-wall conditions are imposed at the surface of the sphere. The computational domain is discretized with quadratically curved hexahedron elements. For the initial mesh, the radial dimension is partitioned with 14 elements with a linear stretching factor of 1.1 while the azimuthal dimension of the plane at $x = 0$ is partitioned using 12 elements. DGP2 is applied for this case and the CFL number is 0.8CFL$^{\text{EB}}$ with CFL$^{\text{EB}}=0.056$.

\begin{figure}[!htb!]
     \begin{subfigmatrix}{2}
     \subfigure[\label{SPHERE_SOLU_MESH_BL}Baseline mesh.]{\includegraphics[width=0.48\textwidth, clip=, keepaspectratio]{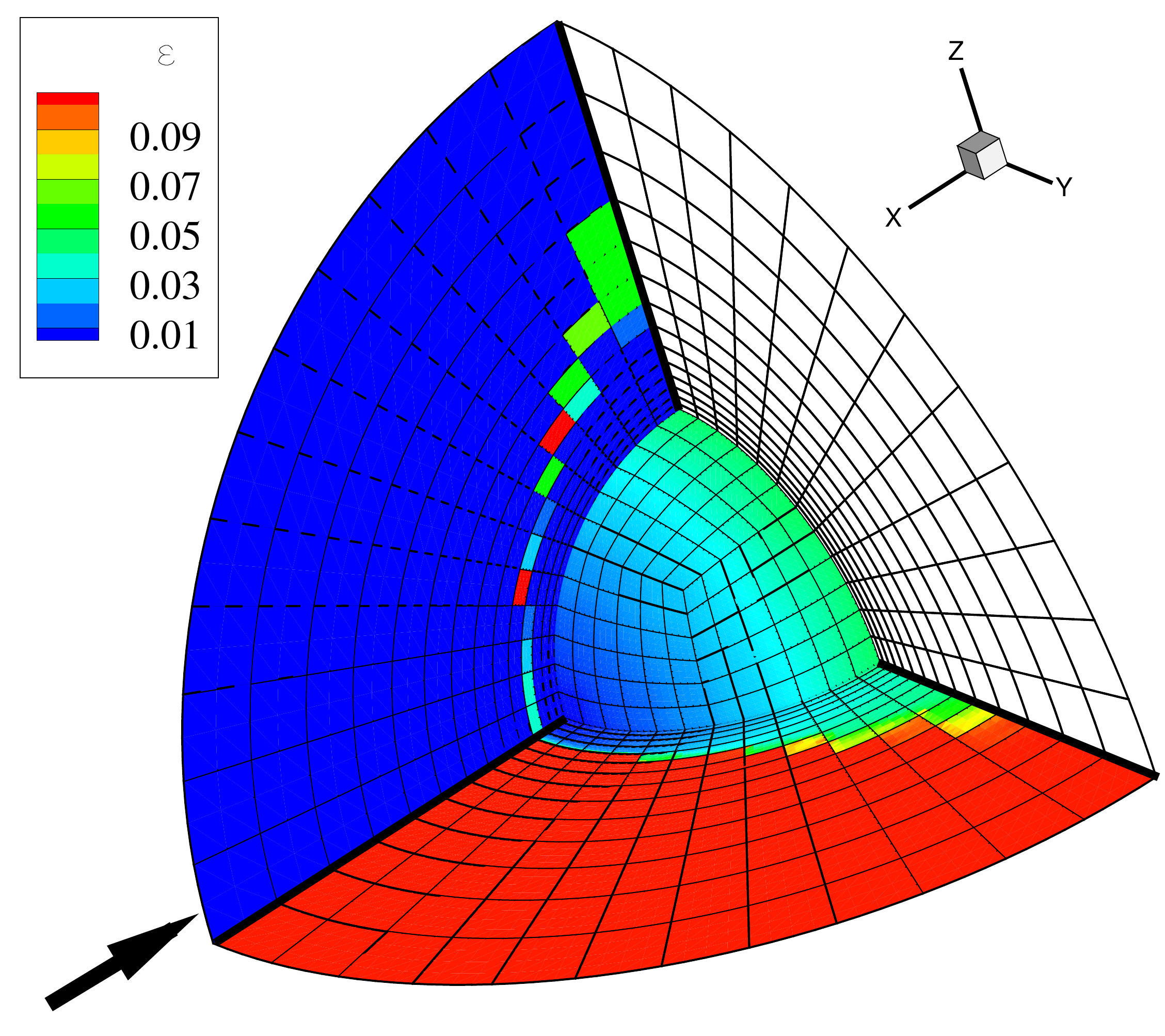}}
     \subfigure[\label{SPHERE_SOLU_MESH_L1}Local refinement at level one.]{\includegraphics[width=0.48\textwidth, clip=, keepaspectratio]{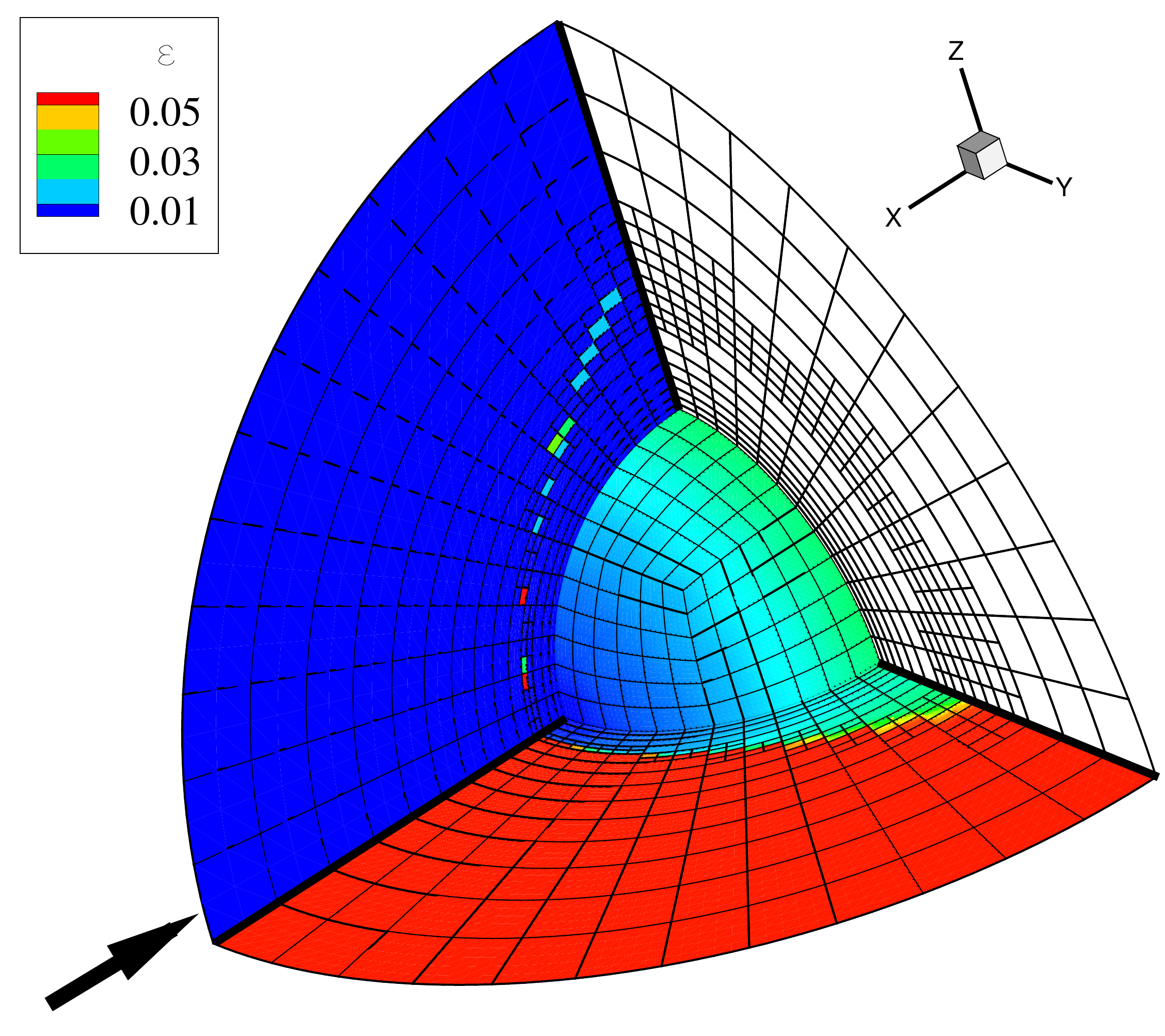}}
     \subfigure[\label{SPHERE_SOLU_MESH_L2}Local refinement at level two.]{\includegraphics[width=0.48\textwidth, clip=, keepaspectratio]{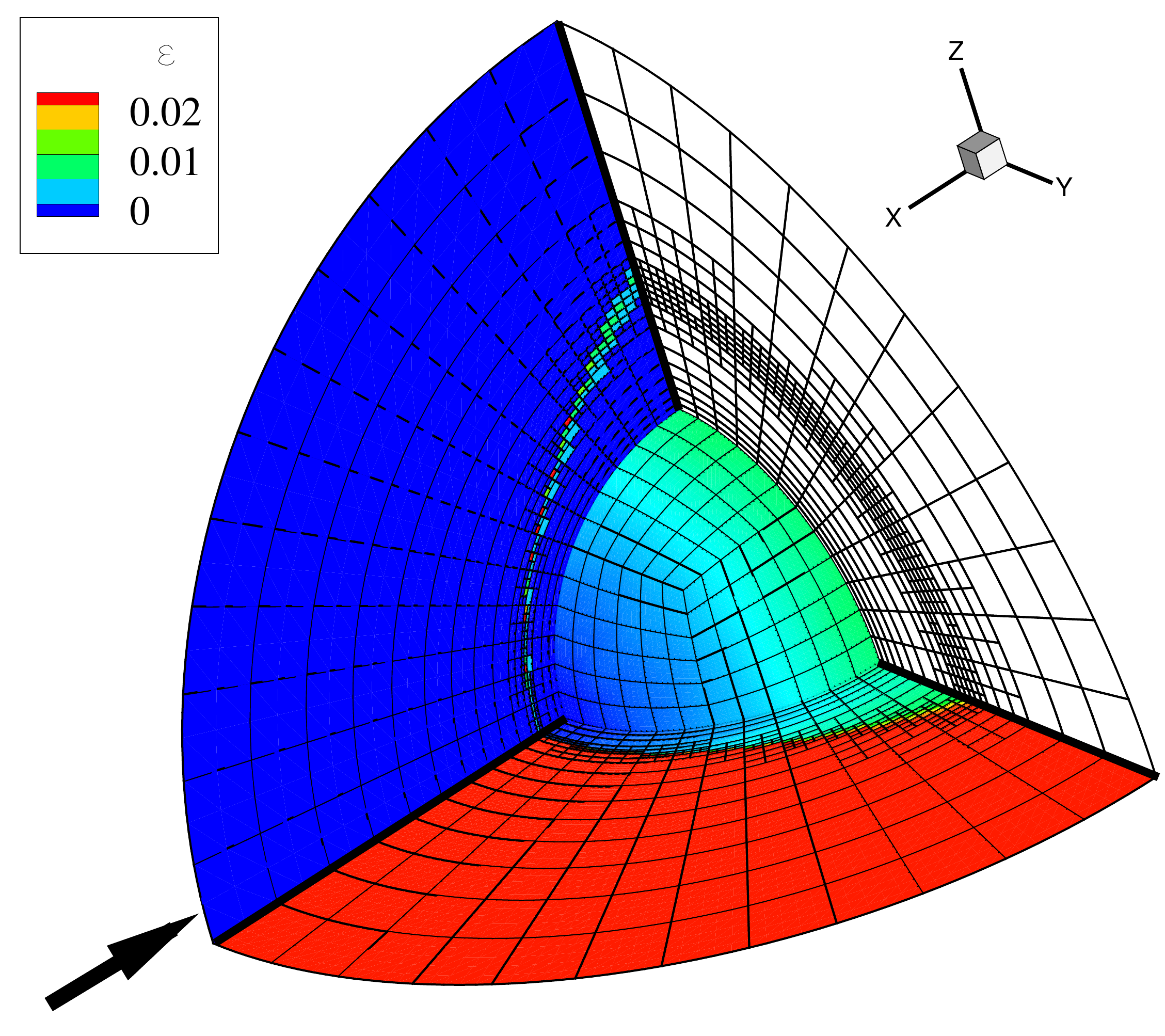}}
     \subfigure[\label{SPHERE_SOLU_MESH_L3}Comparison with measurement in~\cite{REFDATA_SHOCKSPHERE_1967}.]{\includegraphics[width=0.48\textwidth, clip=, keepaspectratio]{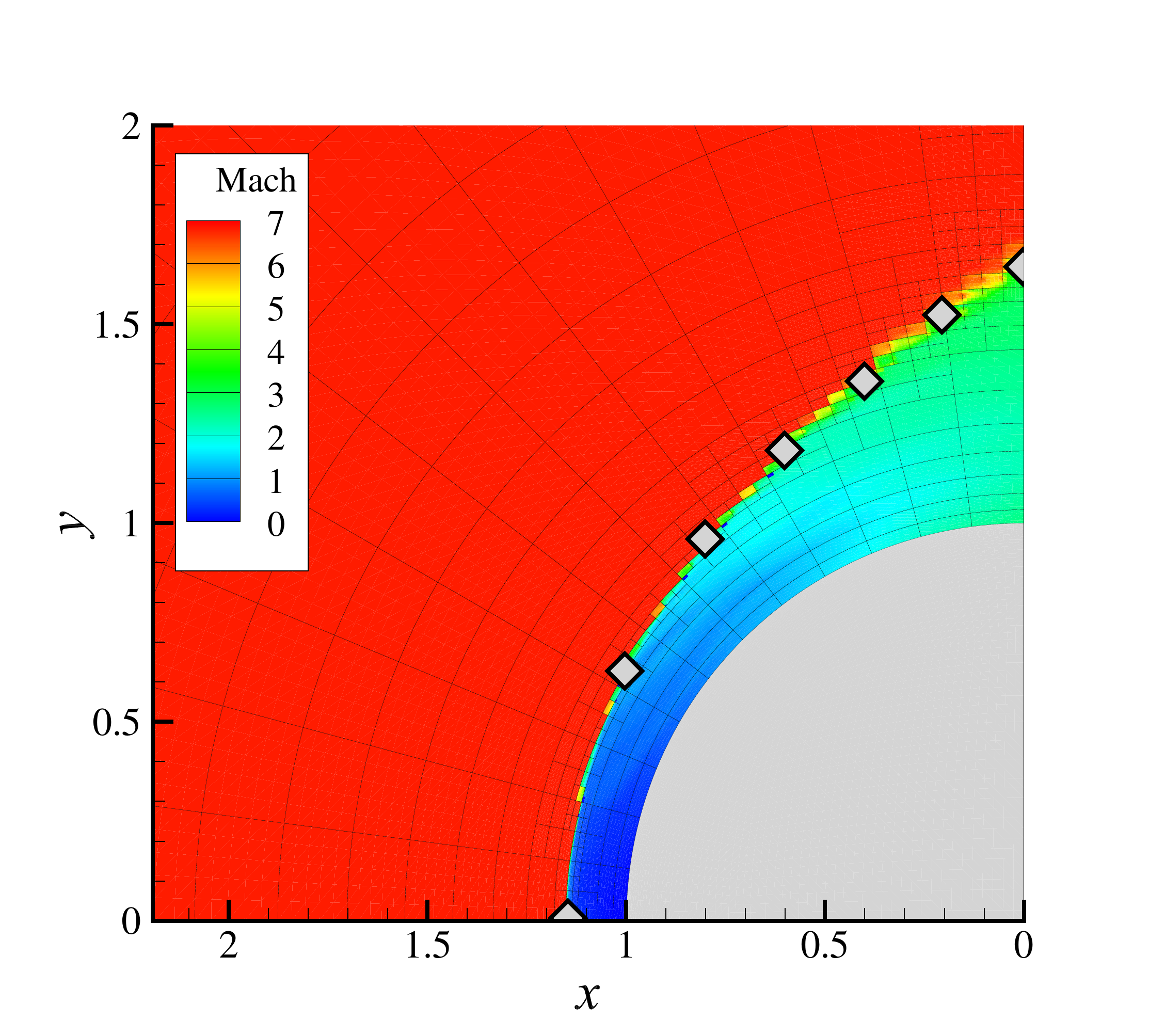}}
     \end{subfigmatrix}
    \caption{\label{SPHERE_SOLU}(Color online) Simulations of Mach 6.8 flow over a sphere showing the $\epsilon$ profile ($y=0$ plane) and Mach-number distribution ($z=0$ plane) on (a) baseline mesh, and simulation results with local refinement with (b) one refinement level and (c) two refinement levels. Comparisons of the shock location with measurements by Billig~\cite{REFDATA_SHOCKSPHERE_1967} are shown in (d).}
\end{figure}

Simulation results are illustrated in Fig.~\ref{SPHERE_SOLU_MESH_BL}, showing the surface mesh and isocontours of the Mach number. The bounding parameter $\varepsilon$ can be utilized as a indicator for local mesh refinement. We sample the elements with non-zero $\varepsilon$ values over few iterations, and then locally refine these elements. Results using one and two levels of refinement are shown in Figs.~\ref{SPHERE_SOLU_MESH_L1} and~\ref{SPHERE_SOLU_MESH_L2}, respectively. This direct comparison shows that the shock profiles become sharper with increasing resolution. Since the bounding parameter is sharp, the mesh-refinement is confined to a narrow region in the vicinity of the shock. The flow-field solution behind the shock is smooth, and no entropy bounding is applied in this region. To provide a quantitative analysis, simulation results from the EBDG-method are compared against measurements by Billig~\cite{REFDATA_SHOCKSPHERE_1967} in Fig. \ref{SPHERE_SOLU_MESH_L3}, showing good agreement between experiments and computations.

\section{Conclusions}
A regularization technique for the discontinuous Galerkin scheme was developed using the entropy principle. Motivated by the FV entropy solution, the high-order DG-scheme is stabilized by constraining the solution to obey the entropy condition. The implementation of the  resulting entropy-bounding discontinuous Galerkin scheme relies on two key components, namely a limiting operator and a CFL-constraint. These essential components were derived by considering first a one-dimensional setting and the subsequent extension to multi-dimensional configurations with arbitrary and curved elements. Specifically, utilizing the interpolation basis we were able to extend the entropy bounding (also including positivity preserving) to arbitrarily shaped elements independent of specific quadrature rules. The bounding procedure is obtained from algebraic operations, resulting in a computationally efficient and simple implementation. A sufficient CFL-condition was rigorously derived and proofed to ensure that the entropy constraint  can be enforced on different types and orders of elements. By considering different configurations, numerical tests were conducted to examine accuracy and stability of the entropy-bounding DG-scheme. These test cases confirm the efficacy in regularizing solutions in the vicinity of discontinuities, generated either by true flow physics or during the transient solution update. The added benefit of the entropy bounding method is its utilization as a refinement indicator. 

Since the herein proposed entropy bounding scheme relies on a linear scaling operator, it is not capable to remove shock-triggered oscillations of smaller magnitude, although it stabilizes the solution and prevents the solver from diverging. As a final remark, the derivation that was presented in this study is general and extendable to other discontinuous schemes with sub-cell solution representations, such as spectral finite volume schemes~\cite{ZJWANG_SPV_2006} and the flux reconstruction scheme~\cite{HUYNH_FR_2013}.  Therefore, entropy-bounding, as an idea, has the potential to improve the robustness of shock-capturing for these emerging high-order numerical methods.

%
%
%
\section*{Acknowledgment}
Financial support through NSF with Award No. CBET-0844587 is gratefully acknowledged. This work used the Extreme Science and Engineering Discovery Environment (XSEDE), which is supported by National Science Foundation grant number ASC130004. Helpful discussions with Yee Chee See on the mathematical analysis are appreciated.
\appendix
\section{\label{APP_COMB_RULE}Combination Rule}
In this section, we derive an estimate for the maximum characteristic speed for convex state solution. For this, we consider a state vector $U$ of Eq.~(\ref{EQ_EULER_STATE}), which is written in the following form:
\be
\label{COMBINATION_U_STAR_LAMBDA}
\mathsf{U} = \sum_{k} \beta_k \mathsf{U}_k\;,
\ee
where $\beta_k > 0 $ and $\sum_k \beta_k  =1 $. The maximum characteristic speed of $\mathsf{U}$ is:
\ben
\nu(\mathsf{U}) = |u(\mathsf{U})| + c(\mathsf{U}) = |u(\mathsf{U})| + \sqrt{\gamma(\gamma - 1) \left( e(\mathsf{U}) - \frac{1}{2} |u(\mathsf{U})|^2\right)}\;,
\een
in which $u$ and $E$ can be calculated according to Eq. (\ref{COMBINATION_U_STAR_LAMBDA}) as
\ben
u(\mathsf{U})  = \sum_k \alpha_k u(\mathsf{U}_k),~\qquad e(\mathsf{U}) = \sum_k \alpha_k e(\mathsf{U}_k)\;,
\een
and 
\ben
\alpha_k = \frac{\beta_k \rho(\mathsf{U}_k)}{\sum_k \beta_k \rho(\mathsf{U}_k) }\;
\een
is a new set of coefficients that is introduced to convert from conservative to primitive variables. Furthermore, because of
\ben
\gamma(\gamma -1) e(\mathsf{U}) &=& \sum_k \alpha_k \left(c^2(\mathsf{U}_k) + \frac{\gamma(\gamma-1)}{2}|u(\mathsf{U}_k)|^2\right)\;,\\
|u(\mathsf{U})| & = & \sqrt{|\sum_k\alpha_ku(\mathsf{U}_k)|^2} \;,\\
&\leq& \sqrt{\sum_k\alpha_k |u(\mathsf{U}_k)|^2}\;,
\een
we obtain 
\ben
\nu(\mathsf{U}) &=& |u(\mathsf{U})| + \sqrt{\sum_k \alpha_k c^2(\mathsf{U}_k) +\frac{\gamma(\gamma-1)}{2} \left( \sum_k \alpha_k |u(\mathsf{U}_k)|^2 - {|u(\mathsf{U})|^2}\right)} \;,\\
& \leq &  \sqrt{\sum_k\alpha_k |u(\mathsf{U}_k)|^2} + \sqrt{\sum_k \alpha_k c^2(\mathsf{U}_k) +\frac{\gamma(\gamma-1)}{2} \sum_k \alpha_k |u(\mathsf{U}_k)|^2 }\;,\\
& \leq & \sqrt{2 \sum_k \alpha_k c^2(\mathsf{U}_k) + (2 + \gamma(\gamma-1))\sum_k \alpha_k |u(\mathsf{U}_k)|^2}\;,\\
& \leq & \sqrt{2 + \gamma(\gamma-1)} \sqrt{\sum_k \alpha_k (c^2(\mathsf{U}_k) + |u(\mathsf{U}_k)|^2)}\;,\\
& \leq & \sqrt{2 + \gamma(\gamma-1)} \sqrt{\sum_k \alpha_k (c(\mathsf{U}_k) + |u(\mathsf{U}_k)|)^2}\;,\\
& \leq & \sqrt{2 + \gamma(\gamma-1)}\max_k \left\{c(\mathsf{U}_k) + |u(\mathsf{U}_k)|\right\}\;.
\een
We used this estimation for preselecting the dissipation coefficient $\lambda^*$ used in Lemma \ref{LEMMA_ENTROPY_MULTID}.

\section{Formulation of $\mathcal{J}^\partial_k $}
For a two-dimensional configuration, the curve of an edge is parameterized by $g_k \in \mathbb{R}$, and the surface Jacobian is written as $\mathcal{J}^\partial_k = \left[\frac{\partial x_1}{\partial g_k}, \frac{\partial x_2}{\partial g_k}\right]$; for a three-dimensional configuration, an edge surface is parameterized by $g_k=\left[g_{k}^{(1)}, g_{k}^{(2)}\right]^T \in \mathbb{R}^2$, and the Jacobian can be written as $\mathcal{J}^\partial_k = \left[\frac{\partial x_1}{\partial g_{k}^{(1)}}, \frac{\partial x_2}{\partial g_{k}^{(1)}}, \frac{\partial x_3}{\partial g_{k}^{(1)}}\right] \times \left[\frac{\partial x_1}{\partial g_{k}^{(2)}},\frac{\partial x_2}{\partial g_{k}^{(2)}},\frac{\partial x_3}{\partial g_{k}^{(2)}}\right]$.


 \end{document}